\documentclass[acmsmall,screen]{acmart}
\acmJournal{TELO}
\AtBeginDocument{%
  }

\setcopyright{acmlicensed}
\copyrightyear{20XX}
\acmYear{20XX}
\acmDOI{XXXXXXX.XXXXXXX}
\acmISBN{978-1-4503-XXXX-X/18/06}

\newcommand{\hide}[1]{}

\definecolor{lennartcolor}{rgb}{0.2,0.6,0.2}
\definecolor{pascalcolor}{rgb}{0.4,0.6,0.9}
\definecolor{todocolor}{rgb}{0.7,0.1,0.1}
\definecolor{changedcolor}{rgb}{0.42,0.27,0.57}
\definecolor{addedcolor}{rgb}{0.867,0.176,0.361}

\newcommand{\redacted}[1]{\emph{[anonymized for review]}}

\usepackage{cleveref}
\usepackage{subcaption}

\usepackage{algorithm}
\usepackage[noend]{algpseudocode}
\usepackage{algorithmicx}

\algblockdefx{Input}{EndInput}{\textbf{Input: }}{}
\algblockdefx{Output}{EndOutput}{\textbf{Output: }}{}
\algnotext{EndInput}
\algnotext{EndOutput}

\newcommand{\LU}{\text{LogUnif}}

\usepackage[normalem]{ulem}

\begin{document}

\title[BONO-Bench]{BONO-Bench: A Comprehensive Test Suite for Bi-objective Numerical Optimization with Traceable Pareto Sets}

\author{Lennart Schäpermeier}
\orcid{0000-0003-3929-7465}
\affiliation{%
  \institution{University of Münster}
  \city{Münster}
  \country{Germany}
}
\email{lennart.schaepermeier@uni-muenster.de}

\author{Pascal Kerschke}
\orcid{0000-0003-2862-1418}
\affiliation{%
  \institution{TU Dresden \& ScaDS.AI}
  \city{Dresden}
  \country{Germany}
}
\email{pascal.kerschke@tu-dresden.de}


\begin{abstract}
    The evaluation of heuristic optimizers on test problems, better known as \emph{benchmarking}, is a cornerstone of research in multi-objective optimization.
    However, most test problems used in benchmarking numerical multi-objective black-box optimizers come from one of two flawed approaches: On the one hand, problems are constructed manually, which result in problems with well-understood optimal solutions, but unrealistic properties and biases.
    On the other hand, more realistic and complex single-objective problems are composited into multi-objective problems, but with a lack of control and understanding of problem properties.

    This paper proposes an extensive problem generation approach for bi-objective numerical optimization problems consisting of the combination of theoretically well-understood convex-quadratic functions into unimodal and multimodal landscapes with and without global structure.
    It supports configuration of test problem properties, such as the number of decision variables, local optima, Pareto front shape, plateaus in the objective space, or degree of conditioning, while maintaining theoretical tractability: The optimal front can be approximated to an arbitrary degree of precision regarding Pareto-compliant performance indicators such as the hypervolume or the exact R2 indicator.
    To demonstrate the generator's capabilities, a test suite of 20 problem categories, called \emph{BONO-Bench}, is created and subsequently used as a basis of an illustrative benchmark study.
    Finally, the general approach underlying our proposed generator, together with the associated test suite, is publicly released in the Python package \texttt{bonobench} to facilitate reproducible benchmarking.
\end{abstract}

\ccsdesc[500]{Applied computing~Multi-criterion optimization and decision-making}
\ccsdesc[500]{Theory of computation~Continuous optimization}

\keywords{Multi-objective Optimization, Benchmarking, Performance Assessment, Pareto-optimal}

\received{01 June 2025}
\received[revised]{13 October 2025}
\received[revised]{14 January 2026}

\maketitle

\section{Introduction}

In multi-objective optimization, theoretical results guaranteeing the performance of optimizers are rare and limited \cite{rudolph2000convergence,doerr2021theoretical}.
Therefore, the creation, choice, and usage of test problems have been a critical driving force for this field, especially in its continuous branch.
For example, the DTLZ test suite \cite{deb2002scalable,deb2005scalable} was influential in the development of the NSGA-II algorithm \cite{deb2002fast}, with both having a lasting impact on multi-objective optimization to this day.

The seven test problems in DTLZ are constructed in a bottom-up fashion, leading to analytically well-understood optimal solutions, as well as known optimization difficulties (e.g., multimodality), which help to examine whether a given optimizer can overcome them.
However, DTLZ, and other test suites with similar construction mechanisms, suffer from a limited range of expressiveness, leading to certain biases within the test suites.
In particular, bottom-up approaches lack the capabilities to generate complex multimodal Pareto sets and fronts.
To mitigate one of the apparent biases, separable test problems can be rotated, as performed, e.g., by \citet{igel2007covariance}.
Another issue is that the single-objective problems created by these approaches often lack realistic characteristics.
For example, the DTLZ problems with $k>2$ objectives have unrealistic subproblems as any $k-1$ objectives can reach their optimal values simultaneously \cite{ishibuchi2016performance}, and their single-objective optima lay at the decision boundaries.

Another, more recent, stream of research tries to amend these issues by creating multi-objective test problems from combinations of well-established single-objective optimization problems to create more complex multi-objective test problems.
\citet{brockhoff2022using} create bi-objective test problems by a pairwise combination of problems from the black-box optimization benchmark (BBOB) suite, which is widely considered as gold standard in single-objective continuous black-box optimization.
Due to the diverse complexity of the underlying BBOB set, the structures of the resulting bi-objective problems are much more complex compared to problems constructed using a bottom-up mechanism like DTLZ.
However, optimal solutions can only be approximated, and the emerging properties of the multi-objective problems cannot be well controlled.
In another composition-based approach,
\citet{kerschke2016towards,kerschke2019search} utilized single-objective multiple peaks problems (MPM2) \cite{wessing2015multiple} to create multi-objective problems with complex multimodal structures.
Meanwhile, \citet{toure2019bi} and \citet{glasmachers2019challenges} conducted a more detailed analysis of a particularly specific yet fundamentally important subset thereof: bi-objective convex quadratic functions.
Building upon a deep understanding of this problem class, \citet{schaepermeier2023peak} then showed that the optimal trade-off sets of the resulting bi-objective problems can be approximated to any given precision, leading to the first test problem generator that combines the strengths of both streams: guarantees about its optimal solutions and the power to create complex problem landscapes with and without global structure.

In this work, we propose a problem generator, which builds upon this more recent stream of research and advances the design of scalable test functions for continuous bi-objective optimization.
It enables fine-grained control over the following seven problem properties while preserving approximation guarantees for optimal trade-off solutions:

\begin{itemize}
    \item \emph{Type of multimodality}: Unimodal, or multimodal with/without global structure.
    \item \emph{Degree of multimodality}:  Number of local optima in each objective.
    \item \emph{Pareto set shape}: Axis-aligned, linear, curved, as well as multimodal variants.
    \item \emph{Pareto front shape}: Convex, linear, or concave for both, individual local optima as well as the global trade-off set of non-dominated solutions.
    \item \emph{Conditioning}: Sensitivity to certain directions in search space.
    \item \emph{Discretization}: Plateaus induced by rounding values in the objective space.
    \item \emph{Decision space dimensionality}: Any number of decision variables $n \geq 2$ can be used.
\end{itemize}

Building on this methodology, we propose a test suite for \underline{b}i-\underline{o}bjective \underline{n}umerical \underline{o}ptimization \underline{bench}marking -- BONO-Bench -- consisting of $20$ different problem classes that enable the study of research questions on a standardized subset of problems.
We package this with a standardized evaluation pipeline that (i) automatically generates targets for optimal R2 and hypervolume indicator values (i.e., the only two Pareto-compliant performance metrics that can be computed without a reference set), and (ii) produces runtime profiles based on the optimization progress of benchmarked algorithms w.r.t.~the number of function evaluations.

This work is structured as follows:
\Cref{sec:background} provides the relevant background on multi-objective optimization and the design of test functions.
In \Cref{sec:construction}, our problem generator is presented, starting with its core components -- transformed variants of convex-quadratic functions -- followed by an explanation of how these single-objective building blocks are combined into bi-objective test problems with controllable properties.
The generator section is completed with the description of a mechanism enabling the approximation of the Pareto front to a user-defined quality threshold.
Next, \Cref{sec:bono} introduces BONO-Bench, a suite of 20 bi-objective test problems generated with our proposed method, and designed to cover a broad spectrum of structural characteristics.
A brief experimental study in \Cref{sec:study} demonstrates the properties and practical usability of BONO-Bench.
Finally, \Cref{sec:conclusions} summarizes this work and offers some concluding remarks.

\section{Background} \label{sec:background}

We begin by giving an introduction to (numerical) multi-objective optimization, with a particular focus on different construction techniques for test problems as well as the Pareto sets of bi-objective convex-quadratic problems.

\subsection{Numerical Multi-objective Optimization}

In the most general sense, multi-objective optimization deals with the simultaneous optimization of a decision vector $x \in \mathcal X$ on more than one objective function of interest \cite{miettinen1999nonlinear}.
W.l.o.g., we assume that all objectives are to be minimized:
\begin{equation}
    \begin{split}
        \underset{x \in \mathcal X}{\text{minimize}} \quad & F(x) = (f_1(x),\dots,f_k(x)), \\
        \quad &f_i(x): \mathcal X \rightarrow \mathbb R.
    \end{split}
\end{equation}
More specifically, in this paper, we focus on bi-objective problems, i.e., $k=2$, with $d$-dimensional numerical search spaces, i.e., $\mathcal X \subseteq \mathbb R^d$.
If not mentioned otherwise, we assume that $\mathcal X$ is box-constrained and use, w.l.o.g., $\mathcal X = [-5,5]^d$.

In multi-objective optimization, the comparison of solutions is usually done by using the \emph{dominance} relation: a solution $x_1$ dominates a solution $x_2$, written as $x_1 \prec x_2$, iff $\forall i \in \{ 1,\dots,k\}: f_i(x_1) \leq f_i(x_2)$ and $\exists i \in \{1, \ldots, k\}: f_i(x_1) < f_i(x_2)$.
If all relations are strict, $x_1$ is said to \emph{strongly dominate} $x_2$, otherwise $x_1$ \emph{weakly dominates} $x_2$.
Due to the multiple objectives, there is generally not one singular optimal solution, but rather a trade-off set of optimal solutions which reflect different compromises between the objectives.
The optimal solution set contains all nondominated solutions $X^* = \{x \in \mathcal X | \nexists x' \in \mathcal X: x' \prec x\}$ and is regularly referred to as the \emph{Pareto set}, and its image under $F$ is known as the \emph{Pareto front} $Y^* = F(X^*)$.
Furthermore, in the objective space, the point containing the unconstrained optima of all objectives is called the \emph{ideal} point $Y_I=(\min f_1,\dots,\min f_k)$.
It can also be viewed as the worst point (weakly) dominating all other Pareto-optimal solutions.
Conversely, the best objective vector that is (weakly) dominated by all Pareto-optimal solutions is known as the \emph{nadir} point $Y_N$.
In bi-objective optimization, if the single-objective optimum $x_i^* \in \mathbb R^d$ of each objective $f_i$ is unique, $Y_N=(f_1(x_2^*),f_2(x_1^*))$.

\paragraph{Approaches for Optimizing Multi-objective Problems}

When solving multi-objective optimization problems, the desired outcome is (a specific subset of) the Pareto set.
Some methods, such as the $\varepsilon$-constraint technique \cite{haimes1971bicriterion}, modify the original optimization problem by converting all but one objective into additional constraints $f_i(x) \leq \varepsilon_i$ and solving the resulting single-objective optimization problem.
Similarly, scalarization techniques combine the original objective functions into one common objective, e.g., by performing a linear combination of the objectives.
While these kinds of approaches then generate a singular output, they require extensive prior knowledge about the problem to be solved to produce meaningful scalarizations, and may not be able to reach all areas of the Pareto front, e.g., in concave regions.

Therefore, it is often more desirable to approximate the full Pareto front at first -- requiring multi-objective approaches -- where the approximation itself is either the desired outcome or may then be inspected by a decision-maker to select the most preferable trade-off. 
Due to the difficulties humans face when imagining landscapes with more than three dimensions, multi-objective problems -- typically consisting of at least two dimensions in both the search and objective spaces -- are often approached as if they were black-box problems.
In this setting, i.e., optimizing multi-objective problems from a black-box perspective, multi-objective evolutionary algorithms (MOEAs) -- such as NSGA-II \cite{deb2002fast}, SMS-EMOA \cite{beume2007sms}, or MO-CMA-ES \cite{igel2007covariance} -- are the methodology of choice due to their high flexibility.
In recent years, however, this topic has also gained importance in adjacent domains like machine learning.
In particular, Bayesian approaches like SMAC3 \cite{lindauer2022smac3} are emerging, given their effectiveness when optimizing with limited budgets. 
Other research groups complemented the aforementioned methods by proposing specialized strategies, such as techniques based on direct search and gradient ascent or descent, respectively \cite{audet2010mesh,wang2017hypervolume,grimme2019multimodality,schapermeier2022mole}.

\paragraph{Performance Assessment}

When comparing multiple Pareto front approximations -- e.g., to examine which algorithm works best in a specific configuration for a given problem -- it is necessary to summarize the quality of an approximation set in a scalar performance measure.
The most desirable performance measures fulfill a property called \emph{Pareto compliance}, which requires that the performance indicator is guaranteed to improve when a new nondominated or dominating solution is added to the approximation set.
The hypervolume metric (also known as the $\mathcal S$ metric) is likely the most widely applied such performance measure.
It computes the dominated hypervolume (for $k=2$: area) of the approximation set w.r.t. an anti-optimal reference point, fulfilling Pareto compliance for the area dominated by the reference point \cite{zitzler1998multiobjective,beume2009complexity,guerreiro2021hypervolume}.
More recently, the R2 indicator \cite{hansen1998evaluating,brockhoff2012properties} has also been shown to be Pareto-compliant when computed exactly \cite{jaszkiewicz2024exact,schaepermeier2024reinvestigating,schaepermeier2025r2}.
It requires an ideal/utopian reference point and is Pareto-compliant w.r.t. all solutions dominated by it.
Its main advantage is the much more frequent availability of a usable ideal point -- e.g., a value of zero for loss functions in machine learning, or a cost of zero in real-world applications -- where the attainment of good anti-optimal reference points for the computation of the hypervolume metric is often more difficult.

In addition to the aforementioned Pareto-compliant performance indicators, several other metrics -- in particular, IGD \cite{coello2005solving}, IGD$^+$ \cite{ishibuchi2015modified}, and the Hausdorff distance $\Delta_p$ \citep{schutze2012using} -- are widely used and well-established in the field \cite{audet2021performance,grimme2021peeking}. However, they rely on knowledge of the true Pareto front or a high-quality reference set, which could limit their suitability for realistic scenarios.

\paragraph{Visualizing Multi-objective Problems}

\begin{figure}
    \centering
    \begin{subfigure}{\linewidth}
        \includegraphics[width=0.325\linewidth]{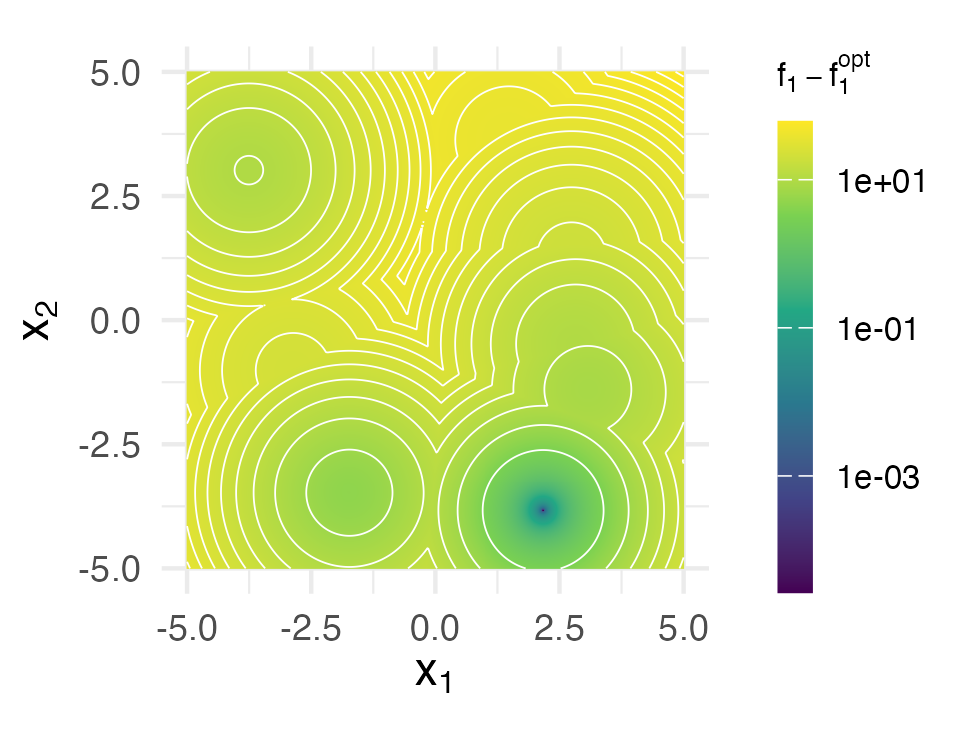}
        \includegraphics[width=0.325\linewidth]{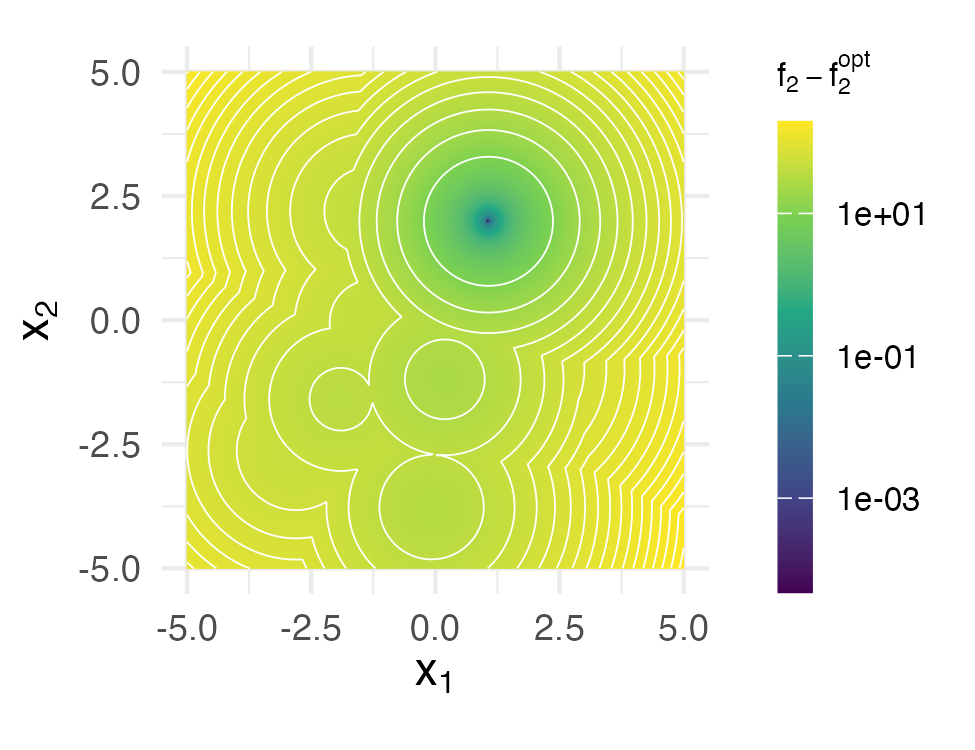}
        \includegraphics[width=0.325\linewidth]{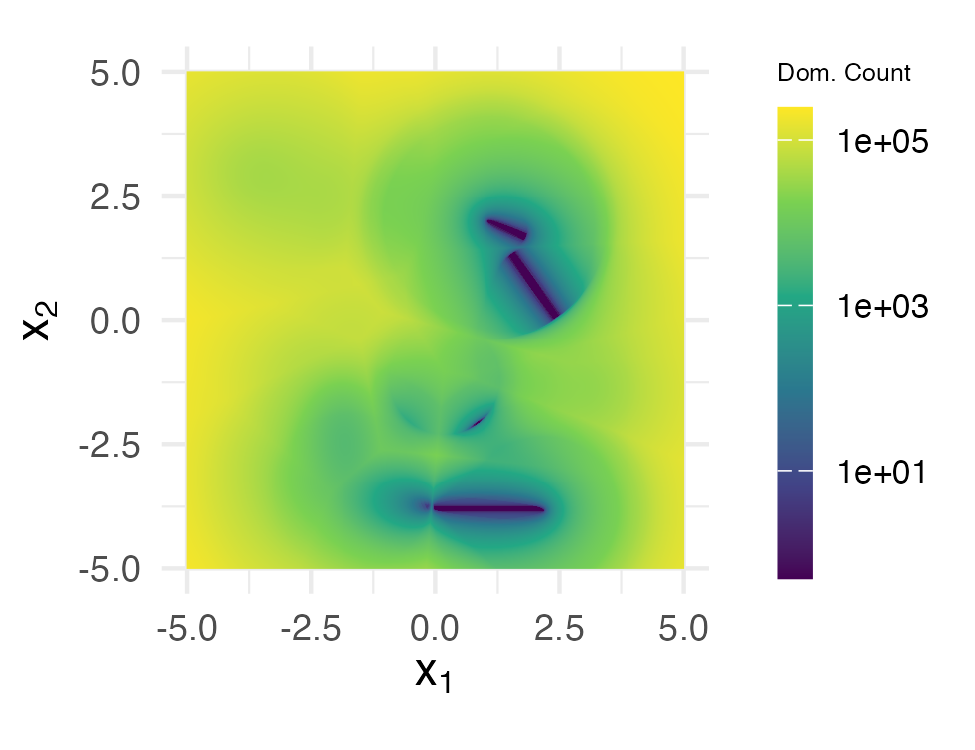}
        \caption{\emph{Left:} Heatmap of $f_1$ values. \emph{Middle:} Heatmap of $f_2$ values. \emph{Right:} Heatmap of domination counts. While the domination counts can reveal the position and distribution of globally optimal solutions, other landscape features stay missing. All heatmaps use log-scaled colormaps.}
        \label{fig:plot_heatmaps}
    \end{subfigure}
    \begin{subfigure}{\linewidth}
        \includegraphics[width=0.27\linewidth]{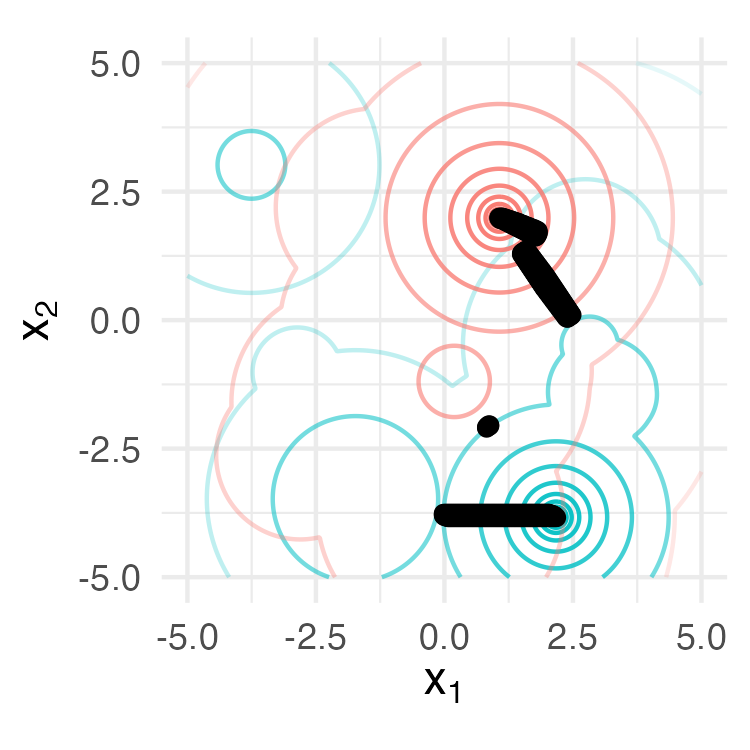}
        \includegraphics[width=0.36\linewidth]{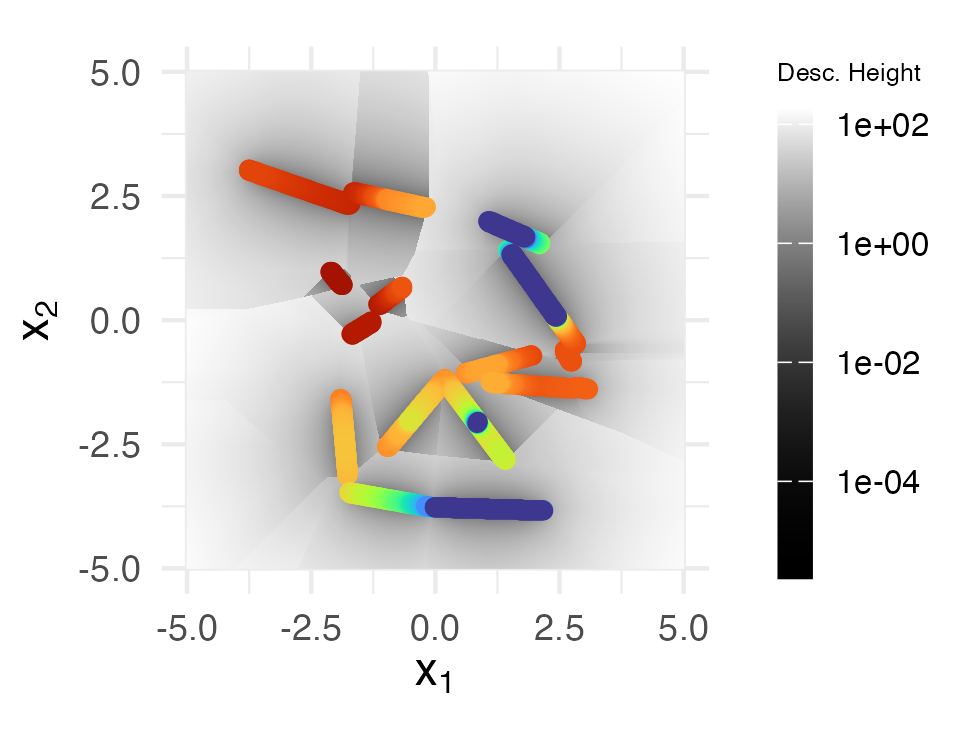}
        \includegraphics[width=0.36\linewidth]{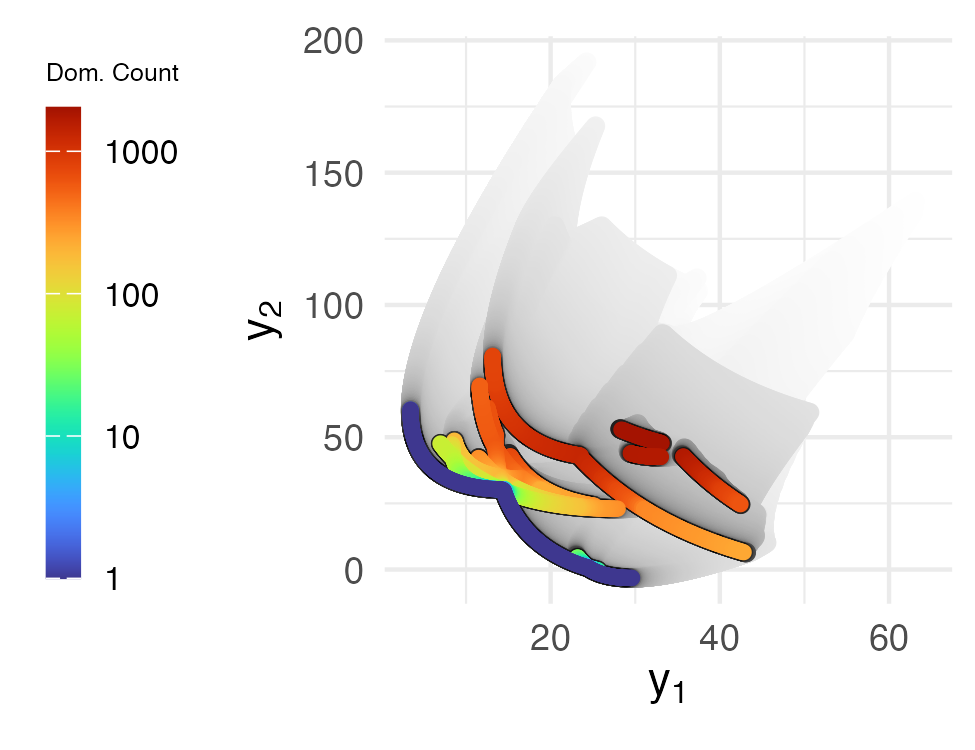}
        \caption{\emph{Left:} Bi-contour plot with nondominated points overlaid. \emph{Middle:} PLOT visualization in decision space. \emph{Right:} PLOT visualization in objective space. The PLOT visualization reveals all locally efficient sets, including the globally efficient sets, as well as attraction basins and interactions between them.}
        \label{fig:plot_worms}
    \end{subfigure}
    \caption{Different landscape visualizations of the same bi-objective test problem (instance 10 of the two-dimensional FewSpheres problem, cf.~\Cref{sec:bono_multi_without}).}
    \label{fig:plot}
\end{figure}

With $k \geq 2$ objectives and $d \geq 2$ decision variables, multi-objective problem landscapes are inherently hard to visualize, requiring pre-processing and aggregations to visualize this at least four-dimensional space in a useful way.
For an overview of recent visualization techniques, we refer to \citet{schapermeier2022plotting}.
An overview of the landscape visualizations discussed below (and used later in this paper) is presented in \Cref{fig:plot}.

With the exception of simple plots that depict each objective individually (left and center image of \Cref{fig:plot_heatmaps}), 
the domination count-based heatmaps by \citet{fonseca1995multiobjective} (right image of \Cref{fig:plot_heatmaps}) are the first visualizations of multi-objective problems.
Similarly simple is a plot combining the contours of both functions with the Pareto optimal solutions (left image in \Cref{fig:plot_worms}).
While these give an insight into the global structure of the single- and multi-objective functions, they are rather basic and cannot give insight into truly multi-objective landscape features such as locally efficient sets and their associated attraction basins.

For problem visualization in this paper, we apply the Plot of the Landscape with Optimal Trade-offs (PLOT), originally introduced by \citet{schaepermeier2020plot}.
It performs an approximate multi-objective gradient descent on the sample used for the visualization, and uses the outcomes to produce a gray-scale gradient field heatmap (depicting proximity to a locally optimal solution), as well as a separate colored visualization of the domination counts of locally optimal solutions, ranging from optimal (blue) to worst (red) locally-optimal solution (center image of \Cref{fig:plot_worms}). Using the same coloring, the problem's distinctive characteristics can also be depicted in its objective space (right image of \Cref{fig:plot_worms}).

The shown problem (instance 10 of the two-dimensional BONO15 problem, i.e., a FewSpheres problem, cf.~\Cref{sec:bono_multi_without}) features multiple locally efficient sets -- making it a \emph{multilocal} problem -- with the images of some of these locally efficient sets (i.e., the local fronts) close to the Pareto front.
This is to distinguish from \emph{multiglobal} problems, for which all locally efficient sets are also globally efficient, i.e., the images of all efficient sets contribute to the Pareto front.
For a differentiated view on multimodality in (numerical) multi-objective optimization, we refer to \citet{grimme2021peeking}.

\subsection{Test Problem Construction Techniques}

In literature, there is a wide array of test function construction techniques for numerical multi-objective test problems.
Despite their multitude, existing problem generators can be divided into two main streams: \emph{bottom-up} and \emph{composite} construction.

\subsubsection{Bottom-up Construction}

Early collections of multi-objective benchmarking test problems, such as the widely used ZDT~\cite{zitzler2000comparison} and DTLZ~\cite{deb2005scalable}, often used a bottom-up construction technique.
That is, when creating the test problem instances, properties like the location of the Pareto set, rudimentary multimodality (often only utilizing shifted copies of the Pareto set), and Pareto front shape were predefined.
This is often achieved by segmenting the decision variables into position variables, whose variation produces locally mutually nondominated solutions, and distance variables, whose variation results in locally dominated solutions.
A bottom-up construction is also used in more recent problem suites, such as the MMF problems which focus on multiglobal multi-objective optimization \cite{yue2019novel}, or the more flexible DBMOPP generator for multimodal distance-based problems \cite{fieldsend2021visualizable}.

This technique often leads to simple problem properties, such as axis-parallel Pareto sets in predictable locations (e.g., in the center or along the boundary of the decision space), separability, and lack of complex multimodality, which are easily exploited by suitable optimizers.
Also, the single-objective subproblems created by bottom-up approaches often lack recognizable properties of usual single-objective test problems.
On the positive side, the exact knowledge of the Pareto set and front allows to generate precise reference sets and indicator values to reliably assess whether a given optimizer has successfully solved the problem to a specified quality.
Bottom-up generators also allow for precise control of Pareto front geometries by construction.

\subsubsection{Composite Construction}

Recent test problem generators more frequently compose multiple separately established single-objective optimization problems into a multi-objective one.
The multi-objective problems then have emerging problem characteristics, which are not completely known up-front.
Probably the most well-known benchmark suite utilizing this technique is the bi-objective BBOB~\cite{brockhoff2022using}, though other studies use this technique, e.g., in conjunction with single-objective (multiple) peaks functions \cite{kerschke2019search,glasmachers2019challenges,schaepermeier2023peak}.

The main advantage of the composite approach is that single-objective problems are constructed in a more principled manner, which makes them more representative of real-world problems.
They can also result in more complex Pareto sets and fronts, which increases the expressiveness of problems created using this approach.
On the flip side, this technique also comes with a trade-off between problem complexity and the knowledge and control about properties of the resulting multi-objective landscapes: While \citet{glasmachers2019challenges} can only produce bi-ellipsoid problems with linear Pareto sets, they can precisely control the positioning of such set and the resulting Pareto front shape.
In the other extreme, bi-objective BBOB \cite{brockhoff2022using} comprises a large variety of single-objective challenges, but neither multi-objective properties nor Pareto set or front are known or can be approximated with any kind of guarantees.

\subsection{Pareto Sets of Convex-quadratic Problems} \label{sec:convex_quadratic}

Convex-quadratic problems are a subset of numerical bi-objective optimization problems with good theoretical understanding in literature \cite{toure2019bi,glasmachers2019challenges}.
In a convex-quadratic bi-objective optimization problem, both objectives $f_i$, $i\in\{1,2\}$, have the following form:
\begin{equation}\label{eq:convex_quadratic}
    f_i(x) = \frac 1 2 (x - x_i^*)^T H_i (x - x_i^*) + y_i^*,
\end{equation}
where, $x_i^* \in \mathbb R^d$ denotes the decision space position of the optimal solution of $f_i$, while $y_i^* \in \mathbb R$ defines the desired objective value at $x_i^*$.
Finally, $H_i \in \mathbb R ^{d \times d}$ is the positive-definite and symmetric Hessian matrix that consists of all second derivatives of $f_i$.
In the context of optimization, the condition number $\kappa$ of the Hessian matrix is an important property that describes the relative sensitivity of the steepest and most shallow descent directions.
For a Hessian matrix $H$, it is computed by the ratio of its largest and smallest eigenvalues:
\begin{equation}
    \kappa(H) = \frac {\lambda_\text{max}(H)} {\lambda_\text{min}(H)}.
\end{equation}
Some of the existing black-box optimization benchmark suites use condition numbers up to $\kappa = 10^6$ to produce test problems with high conditioning \cite{hansen2021coco,glasmachers2019challenges}.

Due to the smoothness and convexity in objective space of the resulting multi-objective problem $F$, we can reliably find all points within the Pareto set using the optima of the linear interpolations of $f_1$ and $f_2$:
\begin{align}\label{eq:interpolation}
    f_t(x) &= (1 - t) f_1(x) + t f_2(x), \ t \in [0,1].
\end{align}
$f_t$ is optimized at the singular point $x_t^*\in \mathbb R^d$ at which its gradient vanishes:
\begin{align}
\begin{split}
     & \nabla f_t(x_t^*) = (1-t) H_1 (x_t^* - x_1^*) + t H_2 (x_t^* - x_2^*) = 0 \\
     & \Rightarrow [(1 - t) H_1 + t H_2] x_t^* - [(1 - t) H_1 x_1^* + tH_2 x_2^*] = 0 \\
     & \Rightarrow x_t^* = [(1 - t) H_1  + t H_2]^{-1} [(1 - t) H_1 x_1^* + tH_2 x_2^*].
\end{split}
\end{align}
That is, using $H_t = (1-t)H_1 + tH_2$, the resulting Pareto set (PS) of $F$ is parametrized as:
\begin{equation}
    PS(F) = \left\{H_t^{-1} [(1 - t) H_1 x_1^* + tH_2 x_2^*] \ | \ t \in [0,1]\right\}.
\end{equation}
Depending on the relative positioning of the Hessian matrices, the Pareto set may be linear or curved \cite{glasmachers2019challenges}.
In particular, if $H_1=H_2$, the resulting Pareto set is always linear.
Furthermore, \citet{toure2019bi} observed that the position of the Pareto set does not change under strong monotone transformations of the objectives.
This enables the creation of different Pareto front shapes, such as convex, linear or concave fronts with various curvatures \cite{glasmachers2019challenges}.
This also allows to apply the parametrization of the Pareto sets of the individual subfunctions of multimodal multi-objective problems as long as those are monotonically transformed convex-quadratic problems as is the case with the MPM2 generator \cite{wessing2015multiple,kerschke2019search,schaepermeier2023peak}.
We address this in more detail in \Cref{sec:approximation}.

If we also specify $y_t^* = (1-t)f_1(x_t^*) + tf_2(x_t^*)$, we can derive the full canonical quadratic form of the resulting interpolated problem:
\begin{align}
    f_t(x) &= \frac 1 2 (x - x_t^*)H_t(x-x_t^*) + y_t^*.
\end{align}
Later, we will need the sum of two convex-quadratic problems, which can be viewed as a special case of the linear interpolation:
\begin{align}\label{eq:quadratic_addition}
\begin{split}
    f_1(x) + f_2(x) &= 2 f_{0.5}(x) \\
    &= \frac 1 2 (x - x^*)^T H (x - x^*) + y^*,
\end{split}
\end{align}
where $H = H_1 + H_2$, $x^* = H^{-1}(H_1 x_1^* + H_2 x_2^*)$, and $y^* = f_1(x^*) + f_2(x^*)$.

\section{A Generator for Bi-objective Test Problems with Known Optima}
\label{sec:construction}

In this section, we introduce the components of our test problem generator, which is capable of creating diverse bi-objective numerical test problems with traceable Pareto sets.
In Section~\ref{sec:peaks}, we summarize the fundamental building block of our test problem generator, i.e., single-objective peak functions along with various transformation operators to configure the problem's complexity.
Next, we describe how those functions can be combined to produce bi-objective problems (Section~\ref{sec:biobjective}) that are either unimodal, or multimodal with or without global structure.
At last, we describe how we can ensure Pareto front approximation up to a user-specified threshold (Section~\ref{sec:approximation}).

\subsection{Single-objective Peak Problems}\label{sec:peaks}

Convex-quadratic functions as described in \Cref{sec:convex_quadratic} also serve as the fundamental building block for all test problems constructed in our approach.
In particular, we apply the following strictly monotone transformations to the convex-quadratic problems: 
\begin{equation} \label{eq:peak}
    f(x) = s \cdot \left( \frac 1 2 (x - x^*)^T H (x - x^*) \right)^{p/2} + y^*.
\end{equation}
Here, $p > 0$ determines that the distance to $x^*$ behaves similarly to a corresponding $L_p$-norm, and $s > 0$ influences the scaling of objective values.
Using the same Hessian in both objectives, $p$ directly controls the front shape: $p < 1$ in both objectives produces a concave Pareto front, $p=1$ creates a linear front and $p>1$ results in a convex Pareto front.
Note that all of these transformations are strictly monotonic, so that the location of level sets, and thus the Pareto set does not change compared with the original quadratic problem.

Finally, we optionally make use of a discretization operator, which rounds the objective values to the nearest multiple of a given stepping value $h > 0$.
We denote the discretization operation as $\lfloor y \rfloor_h = h \cdot \left \lfloor \frac y h \right \rfloor$.
While this transformation is not strictly monotonic, its weak monotonicity still allows to perform high-quality approximations of Pareto sets (cf.~\Cref{sec:approximation}).
Concluding the specification of single-objective peak functions, we arrive at the following formulation for discretized problems:
\begin{equation} \label{eq:stepped_peak}
    f_h(x) = \left\lfloor s \cdot \left( \frac 1 2 (x - x^*)^T H (x - x^*) \right)^{p/2} \right\rfloor_h + y^*,
\end{equation}
with $h>0$ being the rounding step size.
As $h \rightarrow 0$, $f_h$ more and more closely approximates a corresponding function $f$ without discretization, so we can consider $h=0$ as the special case from \Cref{eq:peak}, where no discretization is applied.

\subsection{Generating Bi-objective Problems with Varying Properties} \label{sec:biobjective}

In the following, we combine the unimodal peak functions from \Cref{sec:peaks} to produce bi-objective problems of varying complexity. In particular, we elaborate on how the different parameters of \Cref{eq:peak,eq:stepped_peak}, respectively, need to be configured to produce problems that are either unimodal (\Cref{sec:unimodal}), multimodal with (\Cref{sec:multi_with}) or without (\Cref{sec:multi_without}) global structure.

\subsubsection{Bi-objective Unimodal Problems} \label{sec:unimodal}

For all bi-objective unimodal problems, we create problems following the form of \Cref{eq:peak,eq:stepped_peak}.
While each problem category is based on a different approach, many properties are shared between the specific generators.
The single-objective optima $x^*_1$ and $x^*_2$ are drawn uniformly at random from $[-4,4]^d$ to ensure that most resulting Pareto sets are automatically located within the decision space boundaries of $[-5,5]^d$. This is necessary, as nonlinear connections between two points from a bounding box $B$ can exceed this initial bounding box.
We use rejection sampling to ensure that they have a distance of at least $|| x^*_1 - x^*_2 || \geq 2$.
Inspired by \citet{glasmachers2019challenges}, the scale parameter is sampled log-uniformly between $10^0$ and $10^6$, i.e., $s_i \sim \LU (10^0,10^6)$, and the optimal objective value for each problem is sampled uniformly from $y_i^* \sim \text{Unif}(-s,s)$ to produce a wide array of objective value ranges.
If not further restricted, the distance parameter $p$ is sampled log-uniformly $p \sim \LU (\frac 1 3, 3)$ and is set identically in both objectives to control the Pareto front shapes to be convex and concave with equal likelihood.
Additionally, in most cases, we create random Hessian matrices with conditioning $\kappa$ by sampling a diagonal matrix with eigenvalues $D = \text{diag}(1, \kappa,\lambda_3,\dots,\lambda_d)$ where $\lambda_i \sim \LU (1,\kappa)$ for $i\geq 3$ and then performing a random rotation $H = R^T D R$, ensuring that $\kappa(H) = \kappa$.

Even with just this unimodal setup, we can already create a wide array of different unimodal test problems (cf.~\Cref{sec:bono_unimodal}).
As discussed by \citet{glasmachers2019challenges}, just covering all cases of bi-objective unimodal problems with linear Pareto sets -- which can be varied regarding the alignment of the Hessians and the decision space grid, alignment of the Hessians with each other, and Pareto front shape -- results in $54$ different problem classes.
In addition to such cases with linear Pareto sets, we also support the creation of more general, curved Pareto sets and problems with discretization.
In our test suite (see Section~\ref{sec:bono}), we distinguish between $5$ cases with linear Pareto sets, as well as $2$ more general cases that result in curved Pareto front shapes, one of which is discretized in the objective space.

\subsubsection{Multimodal Problems with Global Structure} \label{sec:multi_with}

\begin{figure}
    \centering
    \begin{subfigure}{0.49\textwidth}
        \includegraphics[width=\linewidth]{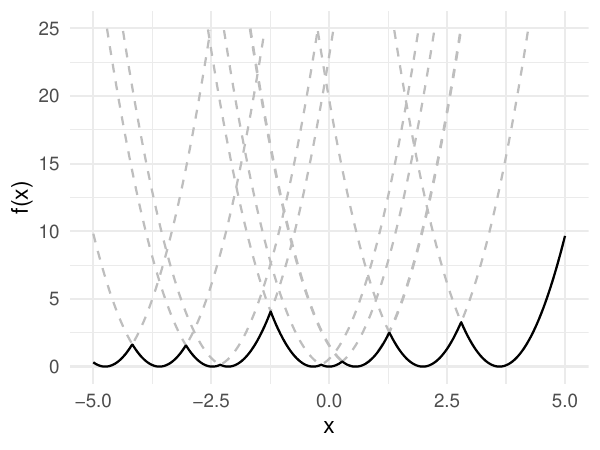}
        \caption{Ten quadratic perturbation functions (dashed) and their minimum (solid).}
    \end{subfigure}
    \hfill
    \begin{subfigure}{0.49\textwidth}
        \includegraphics[width=\linewidth]{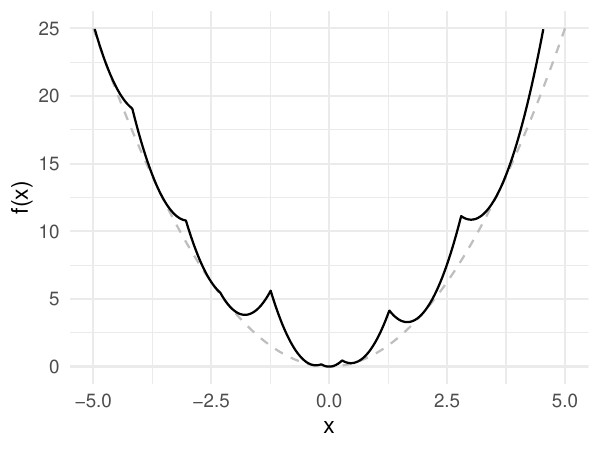}
        \caption{Approximation of the dashed quadratic function with perturbations, resulting in the solid line.}
    \end{subfigure}
    \caption{Illustration of the perturbation of a quadratic function with other quadratic functions. The more perturbation functions are used, the smaller the difference becomes.}
    \label{fig:perturbation}
\end{figure}

In this first category of multimodal problems, we create problems \emph{with global structure} based on the unimodal problems from above.
To achieve this, we perturb the inner quadratic function with a concatenation of further quadratic problems.
Here, the individual objectives can be formulated and simplified as follows:
\begin{align} \label{eq:perturbed_problems}
\begin{split}
    f(x) &= s \cdot \left( \frac 1 2 (x - x^*)^T H (x - x^*) + \min_{j \in \{1, \dots, J\}} \left\{\frac 1 2 (x - x_j^*)^T H_j (x - x_j^*)\right\} \right)^{p/2} + y^* \\
        &= s \cdot \left(\min_{j \in \{1, \dots, J\}} \left\{\frac 1 2 (x - x_j'^*)^T H_j' (x - x_j'^*) + y_j'^*\right\} \right)^{p/2} + y^*,
\end{split}
\end{align}
where $H_j'$, $x_j'^*$ and $y_j'^*$ are the result of adding the original problem and the $j$-th perturbation function using \Cref{eq:quadratic_addition}.
This operation is illustrated in \Cref{fig:perturbation}.
Analogously, a stepped variant following \Cref{eq:stepped_peak} can be defined.
To ensure that the global optimum remains at its established position and the global shape of the function does not change too much, the first perturbation function is identical to the original one, i.e., $H_1 = H$ and $x_1^*=x^*$.
The remaining perturbation functions $H_j$ and $x^*_j$ for $j>1$ are randomly sampled with $x^*_j \sim \mathcal U ([-4,4]^d)$ and $\kappa(H_j) = \kappa(H)$.
Note that the ideal point is not affected by the perturbation, but depending on the approximation quality the nadir point may move further away from it.
For this study, we set the number of perturbation functions for all such generators to $J=500$, which presents a reasonable trade-off between the degree of multimodality, function complexity, and perturbation of the original function.
Setting $J=500$ also limits the amount of potential peak combinations to $J^2=250\,000$, which has a large impact on the runtime of the Pareto front approximation (see below).

\subsubsection{Multimodal Problems without Global Structure} \label{sec:multi_without}

In this final category of multimodal problems, we focus on such without any global structure, i.e., where local optima of the objective functions and the Pareto set do not follow any global trends.
The resulting problem formulation essentially equals the multimodal problems with global structure in \Cref{eq:perturbed_problems}, but comes without a peak function determining the global shape:
\begin{align} \label{eq:random_problems}
    f(x) = s \cdot \left(\min_{j \in \{1,\dots,J\}} \left\{\frac 1 2 (x - x_j^*)^T H_j (x - x_j^*) + y_j^*\right\} \right)^{p/2} + y^*.
\end{align}
The problem creation only differs in how the internal quadratic functions are sampled: The first one is determined to be the global optimum ($y_1^*=0$), while the other internal optima for $j>1$ are sampled from $y^*_j \sim \mathcal U(1,10)$.
This is inspired by the MPM2 generator \cite{wessing2015multiple}, which also introduces a small gap between the global and non-global local optima.
The global optimum $x_1^*$, as well as the other parameters $s$, $p$, and $y^*$ are sampled as usual if not mentioned otherwise.

\subsection{Pareto Front Approximation} \label{sec:approximation}

The key advantage of our problem construction approach is that each objective is composed of one or more monotonically transformed convex-quadratic basis functions, which retain some theoretical tractability in the final problem definition.
This design is essential for enabling the approximation of the problem's Pareto front.

\subsubsection{Composition of the Pareto Front}

For problems following the construction in \Cref{eq:perturbed_problems,eq:random_problems}, note that we can specify the Pareto set (and front) of each combination of convex-quadratic subproblems in the constituent single-objective problem.
As the Pareto set of the composite problem $F$ is made up of the minimum of all of these convex-quadratic subproblems, we can use the Pareto sets of all potential combination of subproblems (or: peaks) to restrict our search for the Pareto front of $F$.
Thus, we can specify the union of all Pareto \emph{sets} (UPS) of each potential peak combination by:
\begin{align}
\begin{split}
    \text{UPS}(F) &= \bigcup_{j\in\{1,\dots,J\},k \in \{1,\dots,K\}} PS(F_{j,k})
\end{split}
\end{align}
where $F_{j,k}$ denotes a multi-objective function that consist only of the $j$-th convex-quadratic component of $f_1$, as well as the $k$-th convex-quadratic component of $f_2$.
We can then use $\text{UPS}(F)$ to define the Pareto set of $F$ as:
\begin{equation}
    \text{PS}(F) = \left\{x \in \text{UPS}(F) \ | \nexists \ x' \in \text{UPS}(F): x' \prec x \right\},
\end{equation}
with the Pareto front $\text{PF}(F) = F(\text{PS}(F))$.
This construction still finds a complete Pareto front when we apply discretization using $\lfloor\cdot\rfloor_h$: No points that are dominated when no discretization is applied will be required to specify the complete Pareto front in the presence of discretization or any other non-strict monotone transformation, for that matter:
Note, however, that in this case the Pareto set $\text{PS}(F)$ may include additional points that also become non-dominated by rounding.
That is, with $x_m$ and $x'_m$ denoting the weakly monotonically transformed vectors, respectively:
\begin{equation}
    x \preceq x' \Rightarrow x_m \preceq x'_m, \text{but } x_m \preceq x'_m \nRightarrow x \preceq x'
\end{equation}
For BONO-Bench, we use this to perform an approximation of $\text{PF}(F)$ for strongly and weakly monotonically transformed functions up to a user-defined accuracy threshold, as we will describe now.

\begin{algorithm}[!t]
    \caption{Pareto Front Approximation}\label{alg:approximation}
    \begin{algorithmic}[1]
        
        \Input{Multiple peak problems $f_1$, $f_2$, indicator $I$, approximation tolerance $\delta_I$}
        \EndInput
        
        \Output{Pareto front approximation $\hat Y$}
        \EndOutput
        
        \Procedure{ApproximateFront}{$f_1, f_2, I, \delta_I$}
            
            \State $\text{front} \gets \textsc{NondominatedArchive}()$ \Comment{Approximated Pareto front}
            \State $\text{queue} \gets \textsc{PriorityQueue}()$ \Comment{Automatically sorts entries by $\varepsilon$}
            \State $\varepsilon_{total} \gets 0$ \Comment{Total approximation error}

            \For{$p_1$ in $\text{peaks}(f_1)$} \label{alg:approximation:peaks} \Comment{Setup sorted priority queue for all peak combinations}
                \For{$p_2$ in $\text{peaks}(f_2)$}
                    \State $t_l \gets 0$
                    \State $x_l \gets \text{xopt}(p_1)$ \label{alg:approximation:opt1}
                    \State $y_l \gets \{p_1(x_l), p_2(x_l)\}$
                    
                    \State $t_r \gets 1$
                    \State $x_r \gets \text{xopt}(p_2)$ \label{alg:approximation:opt2} 
                    \State $y_r \gets \{p_1(x_r), p_2(x_r)\}$

                    \If{$\text{front} \npreceq \text{ideal}(y_l, y_r)$} \Comment{Skip dominated local fronts}
                        \State \text{front}.add($y_l$)
                        \State \text{front}.add($y_r$)
                        \State $\varepsilon \gets \varepsilon_I(y_l, y_r)$
                        \State $\varepsilon_{total} \gets \varepsilon_{total} + \varepsilon$
                        \State \text{queue}.add($\{\varepsilon, p_1, p_2, t_l, t_r\}$) \label{alg:approximation:add_queue}
                    \EndIf
                \EndFor
            \EndFor

            \While{\emph{$\varepsilon_{total} > \delta_I$}} \label{alg:approximation:opt} \Comment{Iterate until approximation tolerance satisfied}
                \State $\varepsilon, p_1, p_2, t_l, t_r \gets \text{queue.pop()}$ \Comment{Get and remove element with largest error from queue}
                \State $\varepsilon_{total} \gets \varepsilon_{total} - \varepsilon$
                \State $t_m \gets (t_l + t_r) / 2$
                \State $x_m \gets \text{xopt}((1 - t_m) p_1 + t_m p_2)$ using \Cref{eq:quadratic_addition}
                \State $y_m \gets \{p_1(x_m), p_2(x_m)\}$
                \State \text{front}.add($y_m$)

                \If{$\text{front} \npreceq \text{ideal}(y_l, y_m)$} \Comment{Add \emph{left} subdivision if not dominated}
                \label{alg:addfront_l}
                    \State $\varepsilon \gets \varepsilon_I(y_l, y_m)$
                    \State $\varepsilon_{total} \gets \varepsilon_{total} + \varepsilon$
                    \State \text{queue}.add($\{\varepsilon, p_1, p_2, t_l, t_m\}$)
                \EndIf

                \If{$\text{front} \npreceq \text{ideal}(y_m, y_r)$} \Comment{Add \emph{right} subdivision if not dominated}
                \label{alg:addfront_r}
                    \State $\varepsilon \gets \varepsilon_I(y_m, y_r)$
                    \State $\varepsilon_{total} \gets \varepsilon_{total} + \varepsilon$
                    \State \text{queue}.add($\{\varepsilon, p_1, p_2, t_m, t_r\}$)
                \EndIf
            \EndWhile
            
            \State \textbf{return} $\text{front}$
    \EndProcedure
    \end{algorithmic}
\end{algorithm}

\begin{figure}
    \centering
    \begin{subfigure}{0.45\textwidth}
        \includegraphics[width=\linewidth]{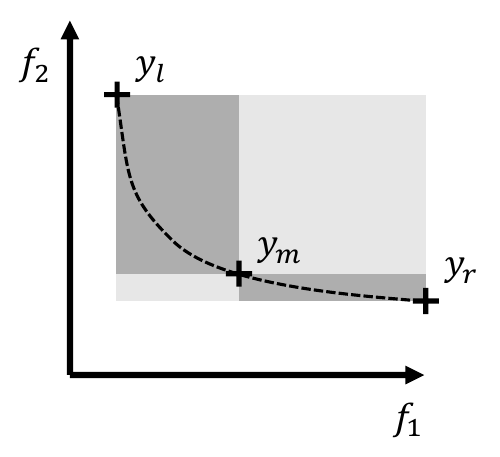}
        \caption{Single set approximation}
        \label{fig:single-set}
    \end{subfigure}
    \begin{subfigure}{0.4\textwidth}
        \includegraphics[width=\linewidth]{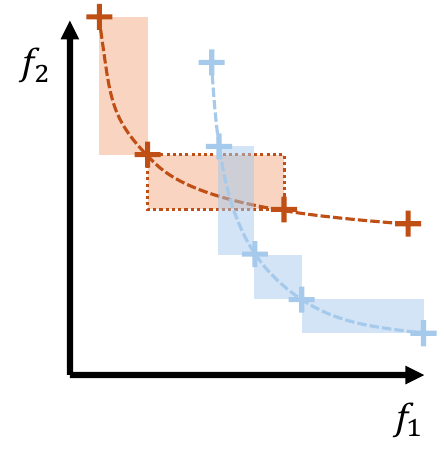}
        \caption{Multiple set approximation}
        \label{fig:multi-set}
    \end{subfigure}
    \caption{Illustration of our Pareto front approximation technique with the hypervolume indicator. The parametrized fronts of the peak combinations are illustrated as dashed lines, and actually known/evaluated points are represented by crosses. Subfigure (a) shows the process of subdividing the front between two evaluated points $y_l$ and $y_r$. The original approximation error $\varepsilon_I(y_l,y_r)$ is given by the lighter shaded area, while the remaining error after the addition of $y_m$ is given by the darker shaded area for $\varepsilon_I(y_l,y_m)$, and $\varepsilon_I(y_m,y_r)$, respectively. Subfigure (b) shows the interactions when multiple peaks contribute to the front: Approximation errors are computed for each front separately, though pairs of consecutive points whose ideal point is dominated by evaluated points are excluded from the further process. The next pair of points to subdivide is decided greedily by the largest gap, and illustrated here by the dotted area.}
    \label{fig:approx}
\end{figure}

\subsubsection{Approximation Procedure}

The complete \textsc{ApproximateFront} procedure to approximate the true Pareto front of a BONO-Bench problem is given in \Cref{alg:approximation}.
It takes two multiple peak problems $f_1$ and $f_2$ as input, as well as a performance indicator of choice, and an approximation threshold $\delta_I$.
The outcome is an approximation $\hat Y$ of the true Pareto front $Y$ for which it holds that $||I(Y) - I(\hat Y)|| < \delta_I$.
The procedure is inspired by \citet{schaepermeier2023peak}, with the main difference that we iteratively choose the next individual point to sample rather than increasing the sampling density for all peak combinations simultaneously.
Some relevant concepts for this procedure are illustrated in \Cref{fig:approx}.

After setting up some data structures to manage the approximation process, starting in line~\ref{alg:approximation:peaks}, we iterate over all pairwise peak combinations and assess the uncertainty of the indicator value between the corresponding ideal and nadir points by evaluating the respective single-objective optima (cf.~lines~\ref{alg:approximation:opt1} and \ref{alg:approximation:opt2}).
All segments are added to a priority queue which allows to manage the potentially Pareto optimal segments (cf.~line~\ref{alg:approximation:add_queue}).
By implementing it using a sorted list, we can perform updates on the queue in $\mathcal O (\log N)$ runtime, each.
In addition to the peaks and the associated indicator uncertainty, we also store the associated $t$-values (cf.~\Cref{eq:interpolation}), which are always $t_l=0$ and $t_r = 1$ during this first step of the initialization.
Further, we maintain a nondominated archive as the Pareto front approximation, which likewise is maintained in a sorted list with (amortized) logarithmic runtime per update.
Note that this step ensures that the endpoints of the true Pareto front are considered in the approximation and the indicator uncertainty can solely focus on the interior of the front.
In this first step, we obtain a worst case runtime of $\mathcal O(N \log N)$, where $N$ denotes the number of peak combinations.

In a second phase, beginning in line~\ref{alg:approximation:opt}, we iteratively remove the segment with the highest uncertainty from the queue and bisect it by interpolating the $t$ values associated with the endpoints, cf.~\Cref{fig:single-set}.
We can then derive $x_m$ from the interpolated peak problem using \Cref{eq:quadratic_addition}, and check whether the segments on the left (right) between $t_l$ and $t_m$ ($t_m$ and $t_r$) can still contribute to the Pareto front and, if so, add them to the queue (cf.~lines~\ref{alg:addfront_l} and \ref{alg:addfront_r}, respectively).
The Pareto front and indicator uncertainty are updated accordingly.
In each iteration, the uncertainty regarding the indicator value of the true Pareto front is reduced, and we terminate the procedure when $\varepsilon_{total}$ first passes below the requested approximation guarantee $\delta_I$.
While we cannot perform a precise runtime analysis for this second part of the approximation, for the hypervolume approximation we can reasonably assume that about $\frac 1 {\delta_{HV}}$ points on the Pareto front are required, which is the number of points required on a linear front.
In the worst case, this approximation has to be repeated for each peak combination, resulting in approximately $\frac N {\delta_{HV}}$ iterations, though practical runtimes are much more efficient in most cases, cf.~\Cref{sec:properties}.
For the exact R2 indicator, a similar asymptotic behavior can be observed.
Our implementation of \Cref{alg:approximation} supports both the normalized hypervolume and exact R2 indicators, which enables us to benchmark against both indicators.

\section{A Test Suite of 20 Bi-objective Test Problems with Known Optima}
\label{sec:bono}

Using the methods described in \Cref{sec:construction}, we now introduce BONO-Bench, a suite of 20 bi-objective test problems.
The suite is designed to cover a broad spectrum of problem characteristics. BONO1 to BONO7 are unimodal problems (\Cref{sec:bono_unimodal}), BONO8 to BONO14 are multimodal problems with global structure (\Cref{sec:bono_multi_with}), and BONO15 to BONO20 represent multimodal problems without global structure (\Cref{sec:bono_multi_without}).
Each problem group comes with its own generator in the \texttt{bonobench} package, which is parametrized according to the specific problem description, cf.~\Cref{tab:bono_unimodal,tab:bono_multi_without}. All parameters discussed in these tables can be modified to generate different problems.

To reveal the different structures of each BONO problem, we provide landscape visualizations for all two-dimensional problem instances with seed $0$ (cf.~\Cref{sec:setup}) in \Cref{fig:axis-aligned,fig:linear-ps,fig:free-ellipsoids,fig:mm-axis-aligned,fig:mm-linear-ps,fig:mm-free-ellipsoids}.
Note that paired problems are seeded identically to enable maximal comparability between, e.g., multimodal and unimodal variants or problems with different front shapes.
For each of the BONO problems, we show (from left to right) heatmaps of the two underlying single-objective peak functions, as well as the corresponding PLOT visualizations in the search and objective space, respectively.
Visualizations for fifteen total instances per two-dimensional problem are available at Zenodo \cite{zenodo2025bonobench}.

\subsection{Unimodal Problems} \label{sec:bono_unimodal}

\begin{table}[t]
    \centering
    \caption{Overview of the most important parameter variations within the unimodal BONO-Bench problems: Conditioning ($\kappa$), distance parameter ($p$), number of discretization steps ($N_h$), whether the Hessians are identical, and if the subproblems are separable. Problems with identical Hessians also have linear Pareto sets. In addition to being separable, BONO1 and BONO2 also have axis-aligned Pareto sets.}
    \label{tab:bono_unimodal}
    \begin{tabular}{cccccc}
        \toprule \bfseries Problem ID & $\boldsymbol{\kappa}$ & $\boldsymbol{p}$ & $\boldsymbol{N_h}$ & \bfseries $\boldsymbol{H_1}=\boldsymbol{H_2}$? & \bfseries Separable?\\
        \midrule
        BONO1  & $1$ & $2$ & - & Yes & Yes \\
        BONO2  & $\LU(10^5,10^6)$ & $2$ & - & Yes & Yes \\ \midrule
        BONO3  & $\LU(50,200)$ & $\LU(1.5,3)$ & - & Yes & No \\
        BONO4  & $\LU(50,200)$ & $1$ & - & Yes & No \\
        BONO5  & $\LU(50,200)$ & $\LU(\frac 1 3,\frac 2 3)$ & - & Yes & No \\ \midrule
        BONO6  & $\LU(50,200)$ & $\LU(\frac 1 3, 3)$ & - & No & No \\
        BONO7  & $\LU(50,200)$ & $\LU(\frac 1 3, 3)$ & $\lfloor\LU(50,201)\rfloor$ & No & No \\ \bottomrule
    \end{tabular}
\end{table}

An overview of the most important parameters for all unimodal BONO-Bench problems is given in \Cref{tab:bono_unimodal}.

\medskip

\begin{figure}[t!]
    \centering
    \begin{subfigure}{\textwidth}
        \centering
        \includegraphics[width=0.245\linewidth]{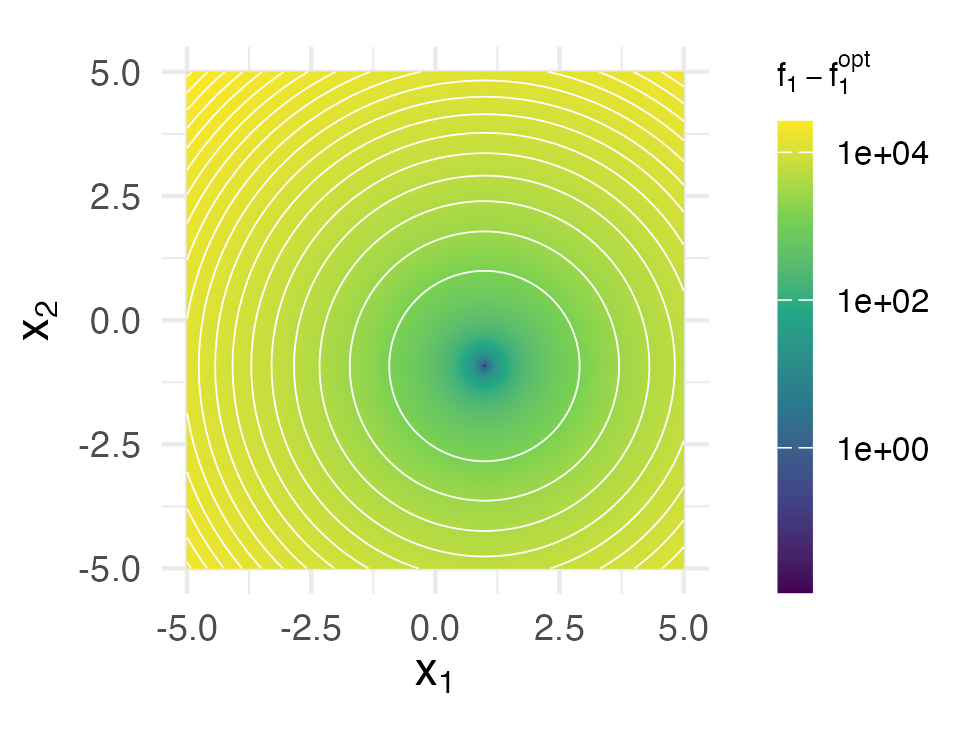}
        \includegraphics[width=0.245\linewidth]{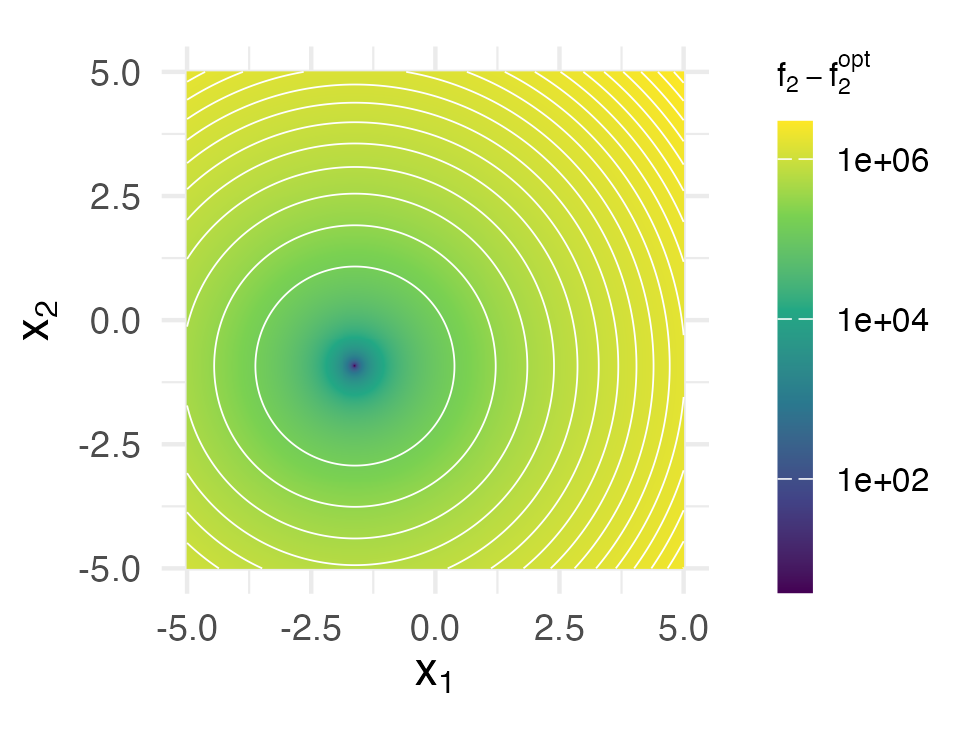}
        \includegraphics[width=0.245\linewidth]{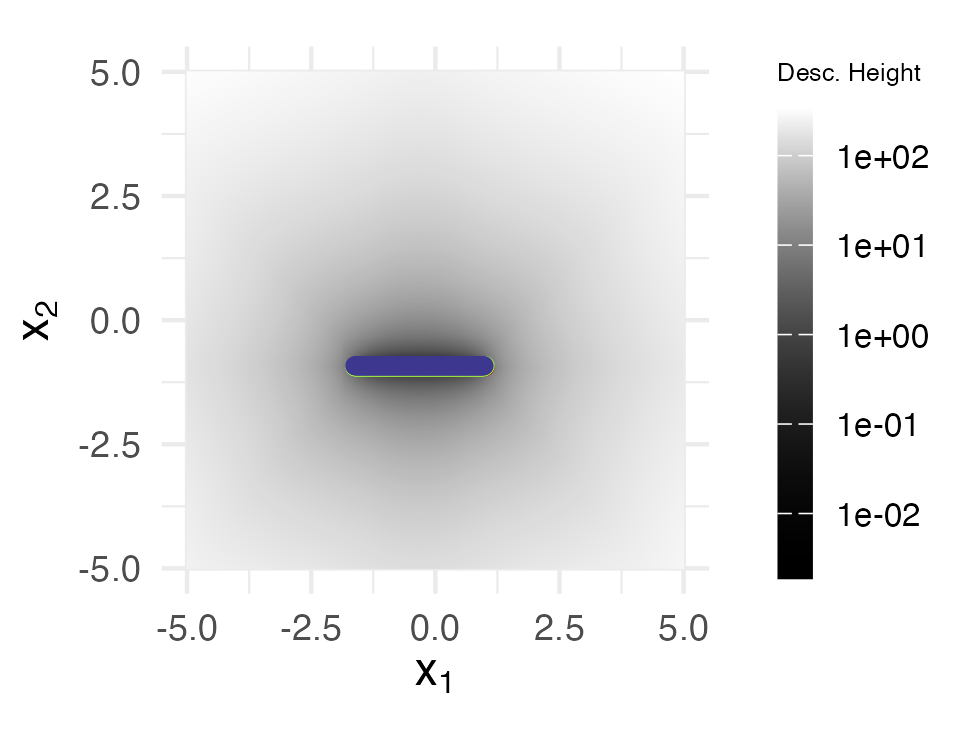}
        \includegraphics[width=0.245\linewidth]{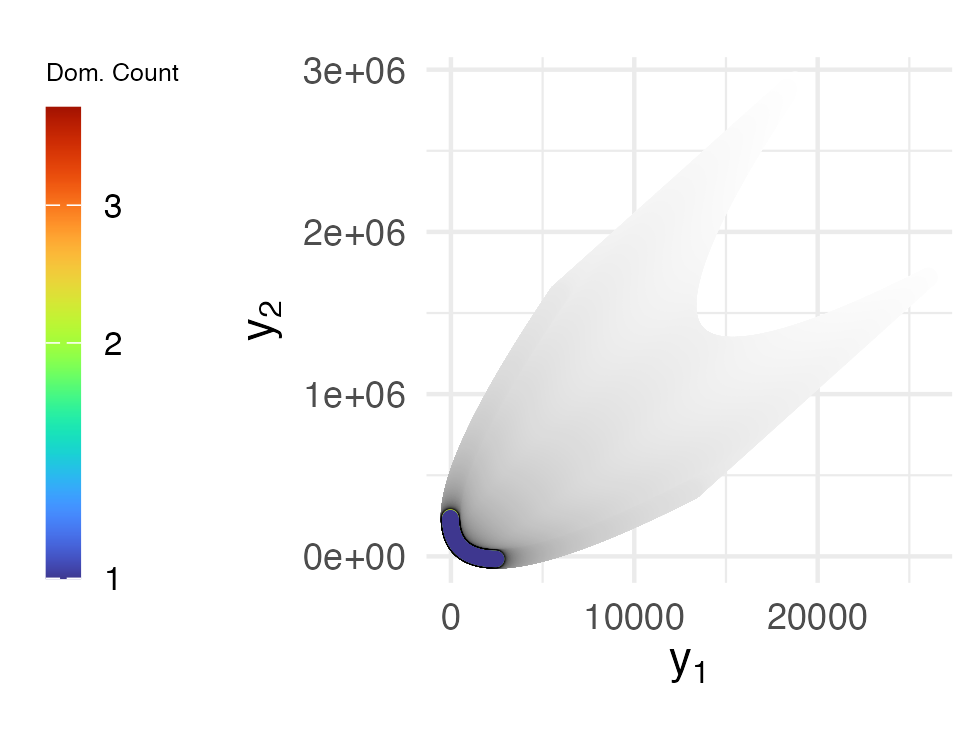}
        \caption{BONO1 instance: Axis-aligned spheres.}
    \end{subfigure}
    \begin{subfigure}{\textwidth}
        \centering
        \includegraphics[width=0.245\linewidth]{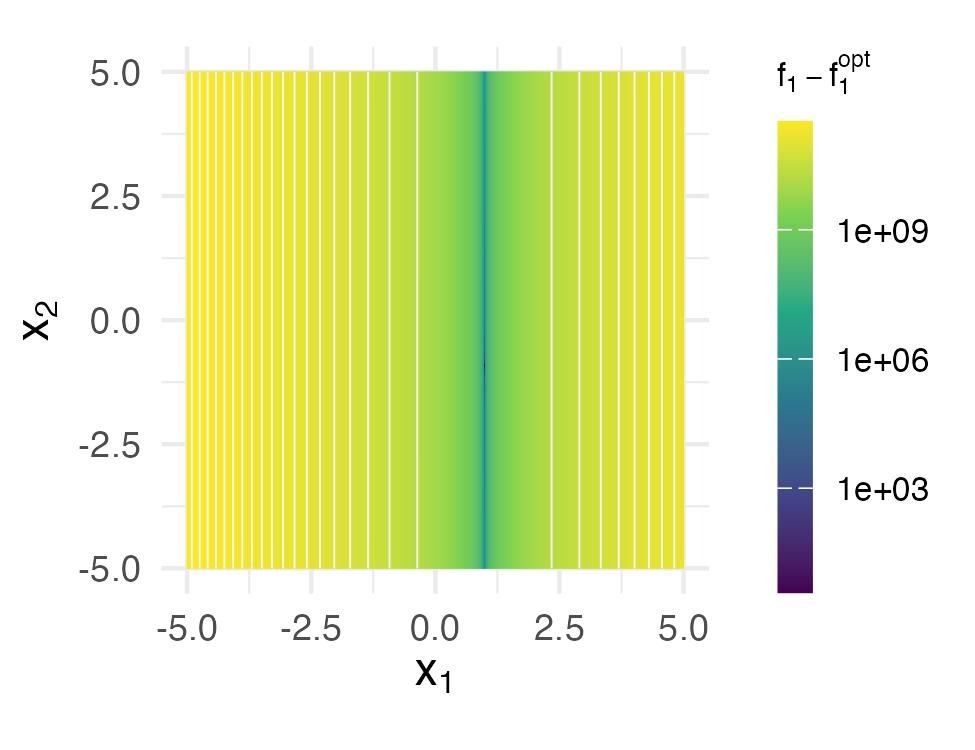}
        \includegraphics[width=0.245\linewidth]{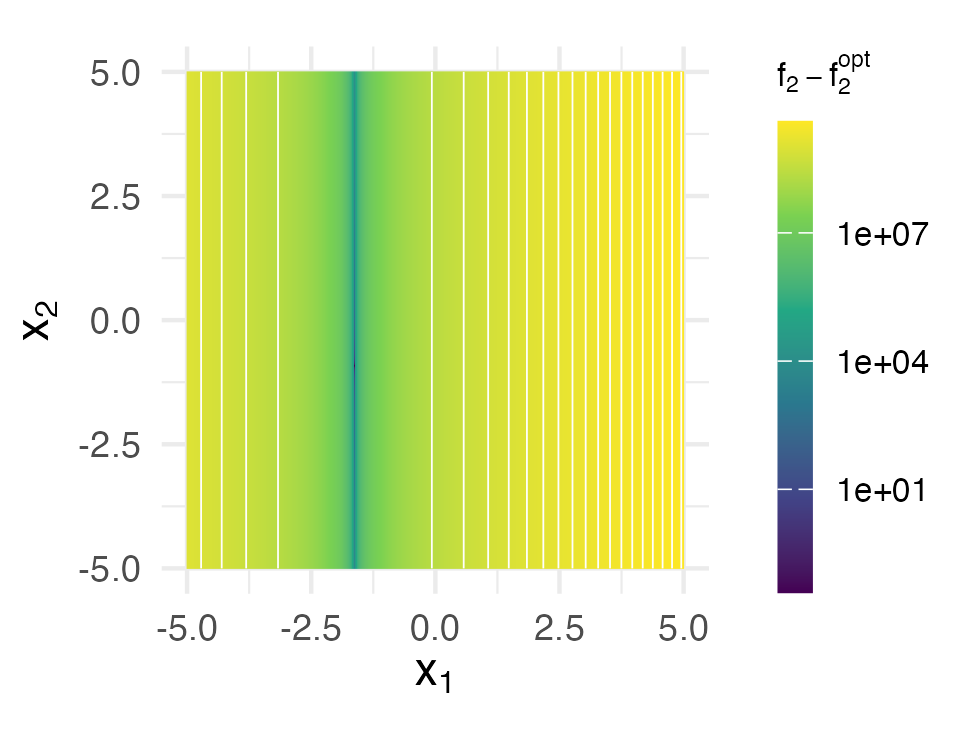}
        \includegraphics[width=0.245\linewidth]{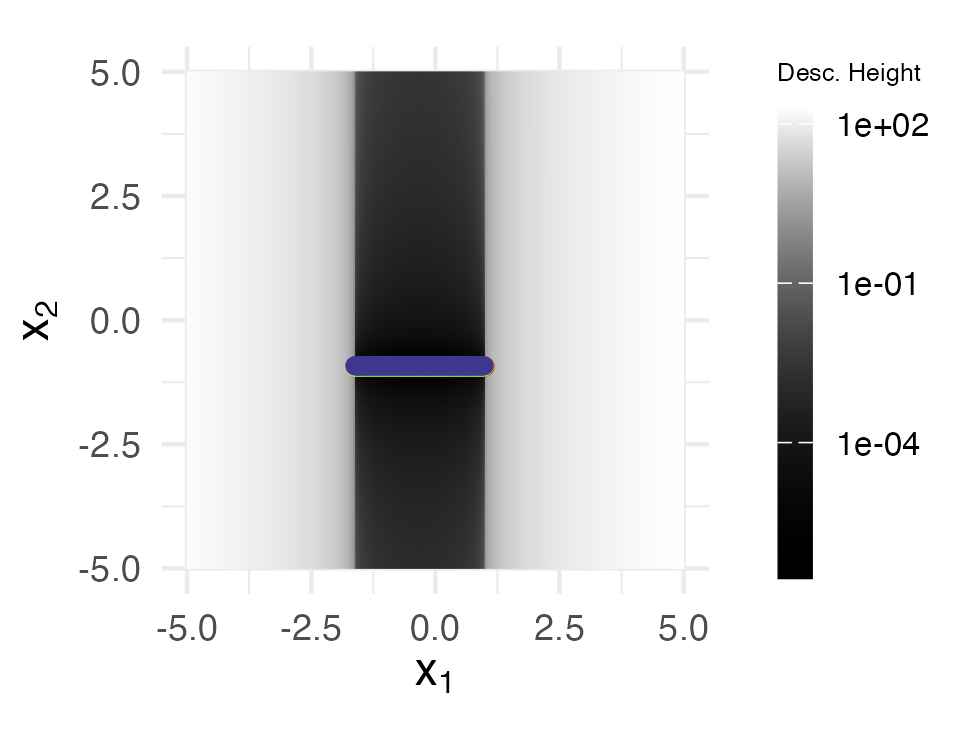}
        \includegraphics[width=0.245\linewidth]{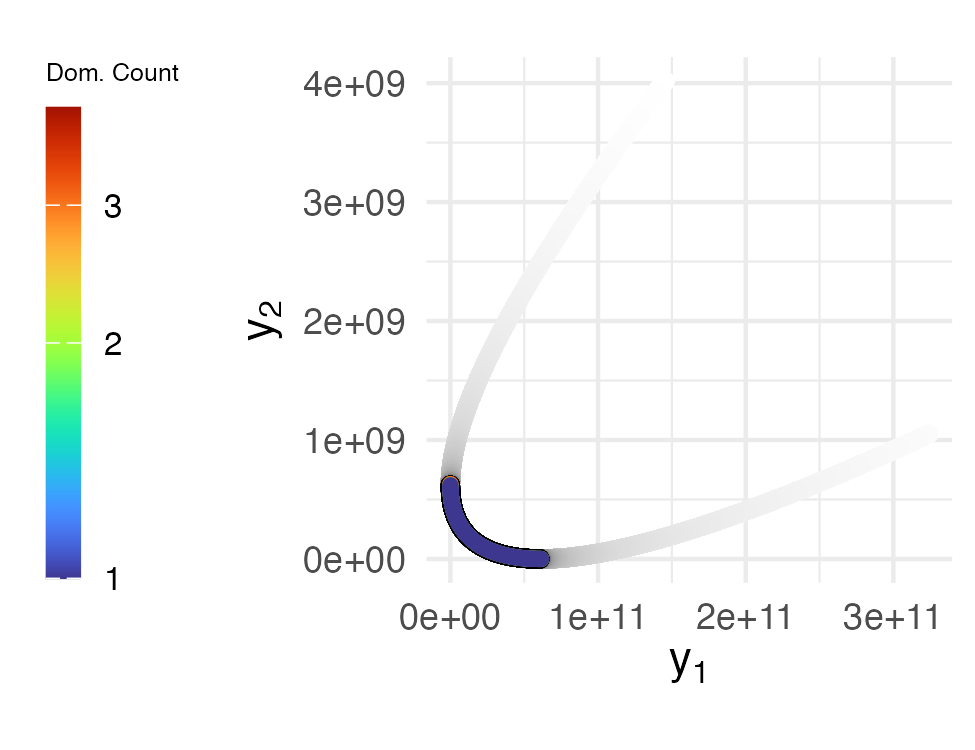}
        \caption{BONO2 instance: Axis-aligned ellipsoids.}
    \end{subfigure}
    
    \caption{Illustration of instances of the unimodal, axis-aligned problems from the BONO-Bench suite, which show log-scaled landscape visualizations of the individual objectives (left) as well as PLOT \cite{schaepermeier2020plot} visualizations of the decision and objective space (right). The upper (sphere) problem features no conditioning, while the lower (ellipsoid) problem is highly conditioned ($\kappa \approx 2.4 \cdot10^5$).}
    \label{fig:axis-aligned}
\end{figure}

The first two problems in BONO-Bench, i.e., BONO1 and BONO2, feature axis-aligned Pareto sets, and examples of both generators are illustrated in \Cref{fig:axis-aligned}.

\subsubsection{BONO1: Axis-aligned Spheres}

The axis-aligned spheres problem is the simplest optimization problem we consider.
It has $H_1 = H_2 = I_d$, a simple convex-quadratic front ($p=2$) and additionally $x^*_1$ and $x^*_2$ only differ on one axis, creating a Pareto set parallel to that axis.

\subsubsection{BONO2: Axis-aligned Ellipsoids}

In contrast to the axis-aligned sphere, BONO2 features a highly conditioned Hessian $H_1 = H_2 = \text{diag}(\lambda_1, \dots, \lambda_d)$ with $\kappa(H_1) \sim \LU(10^5,10^6)$.
Here, we use a permutation instead of a rotation matrix while sampling the Hessian to preserve separability.
It retains all other features of BONO1, namely the convex-quadratic Pareto front ($p=2$), an axis-parallel Pareto set, and separability of the individual objective functions.

\begin{figure}[t!]
    \centering
    \begin{subfigure}{\textwidth}
        \includegraphics[width=0.245\linewidth]{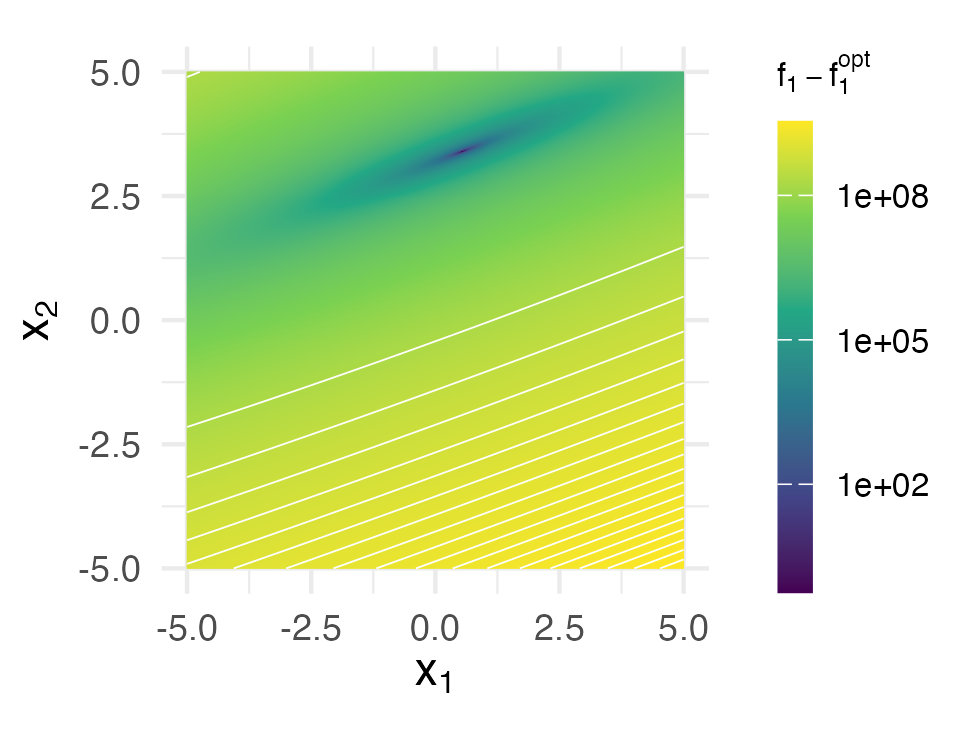}
        \includegraphics[width=0.245\linewidth]{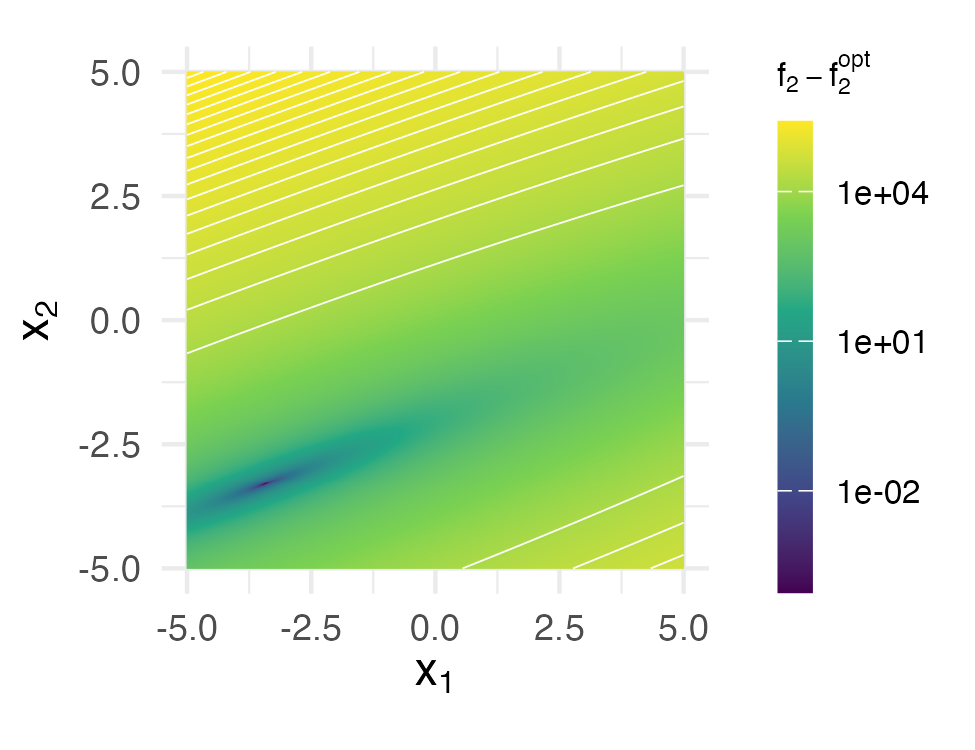}
        \includegraphics[width=0.245\linewidth]{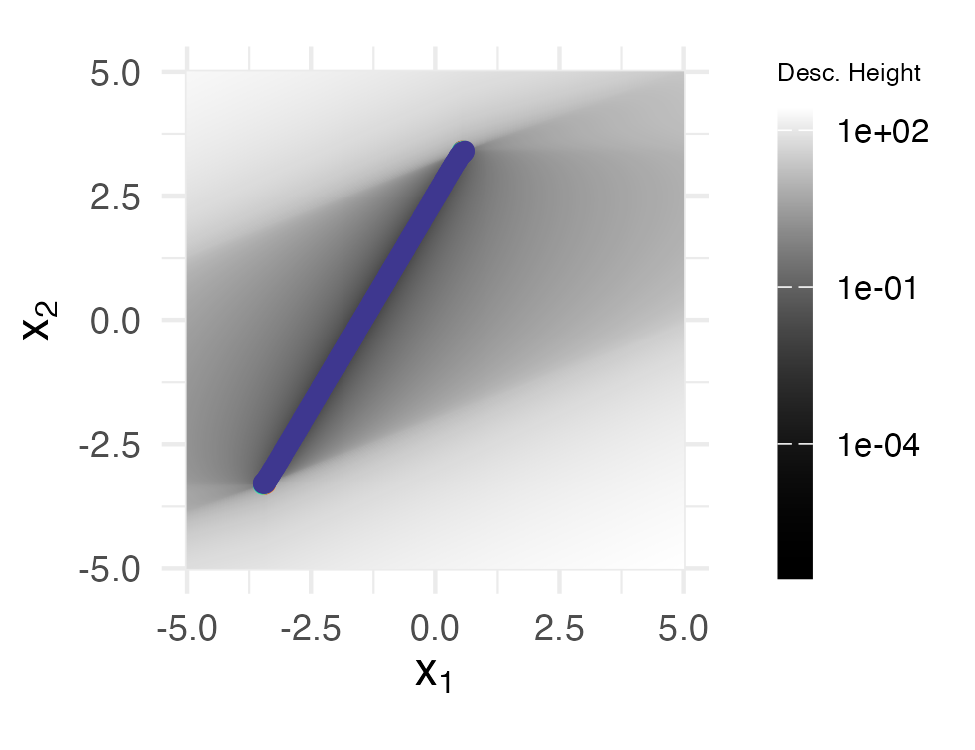}
        \includegraphics[width=0.245\linewidth]{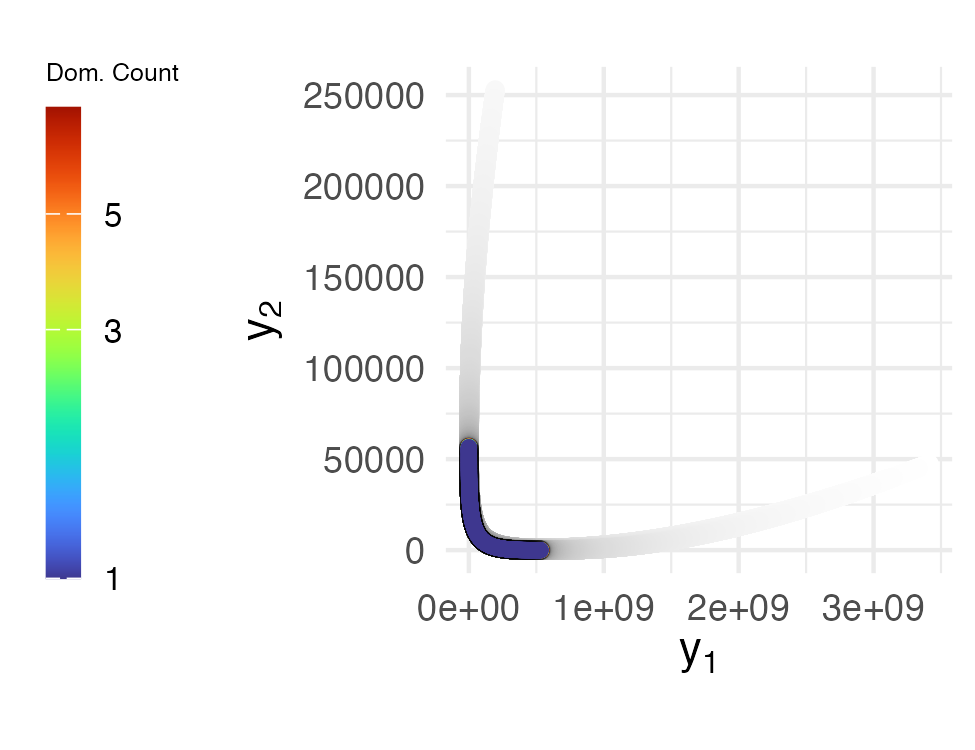}
        \caption{BONO3 instance: Convex ellipsoids.}
    \end{subfigure}
    \begin{subfigure}{\textwidth}
        \includegraphics[width=0.245\linewidth]{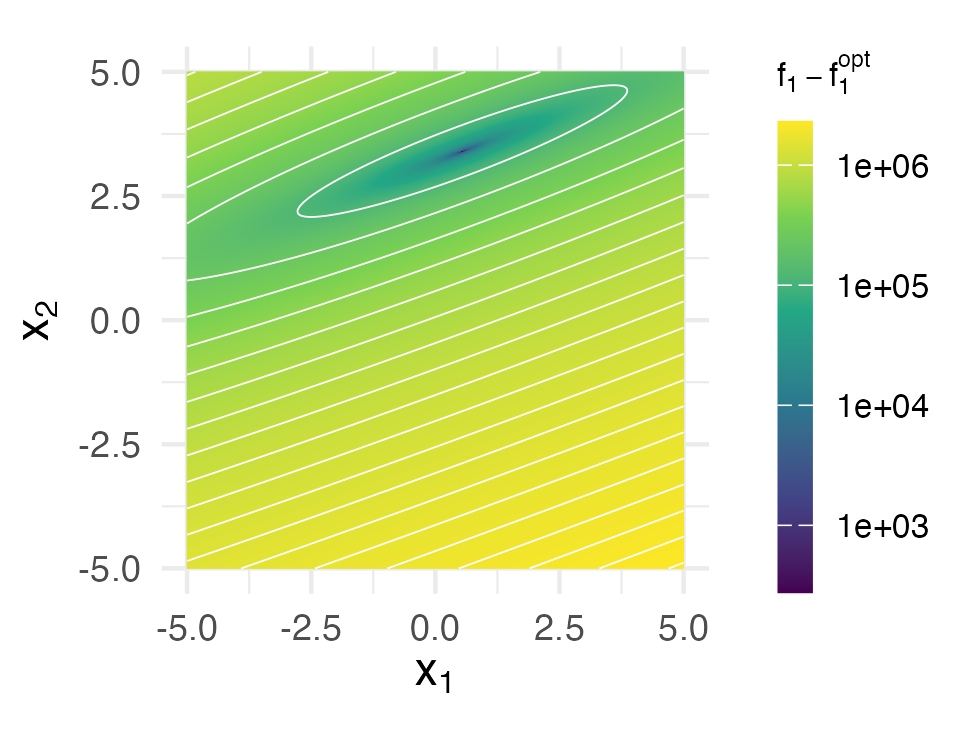}
        \includegraphics[width=0.245\linewidth]{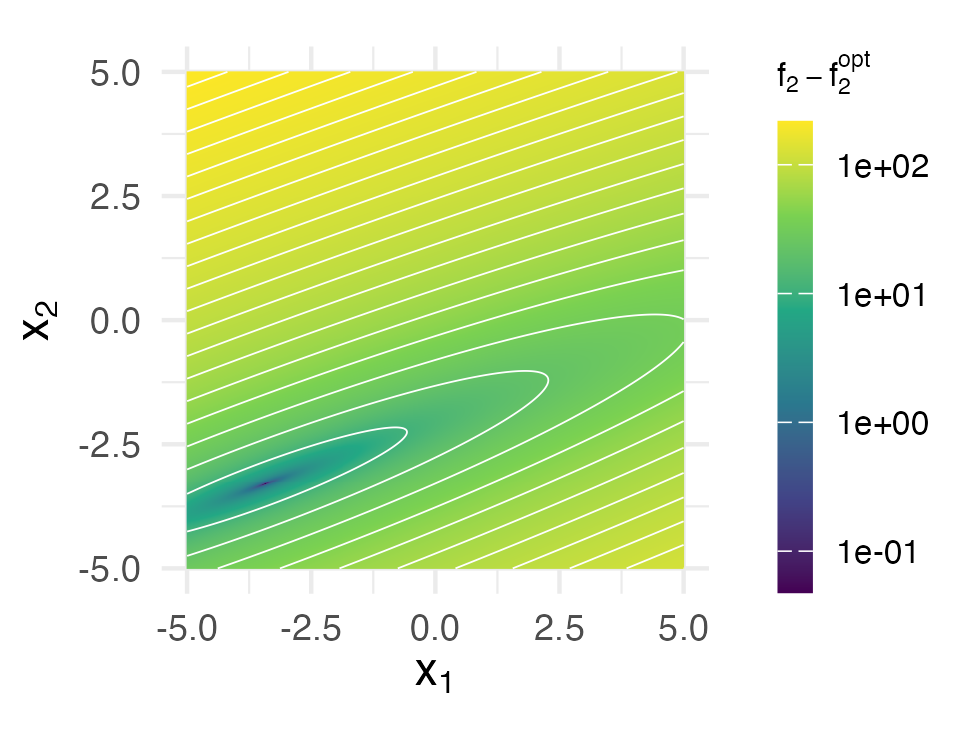}
        \includegraphics[width=0.245\linewidth]{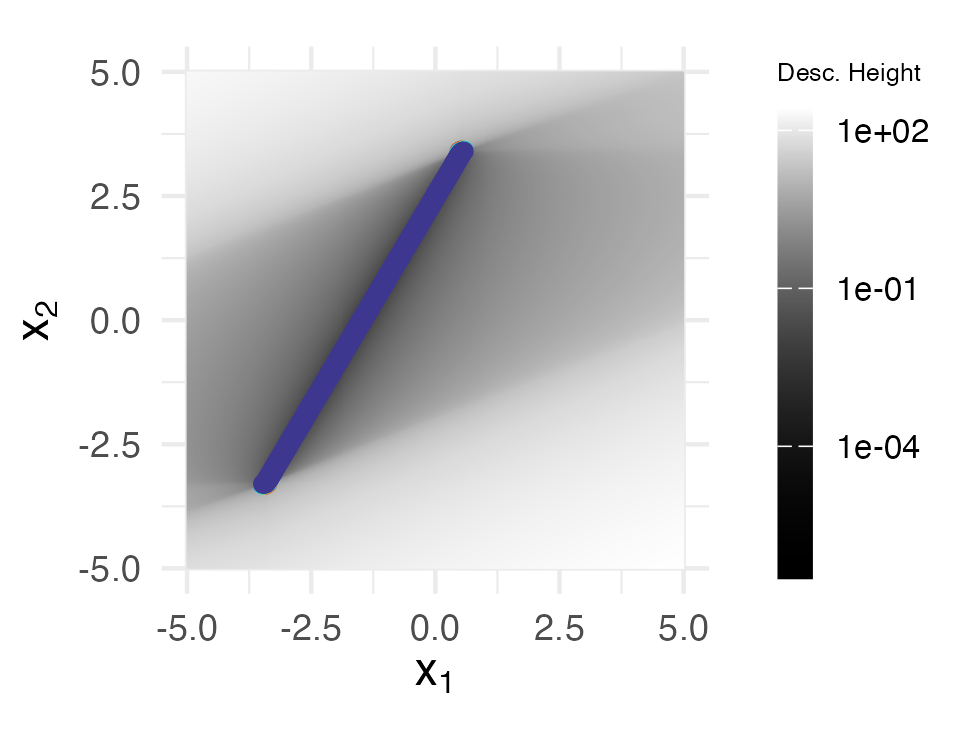}
        \includegraphics[width=0.245\linewidth]{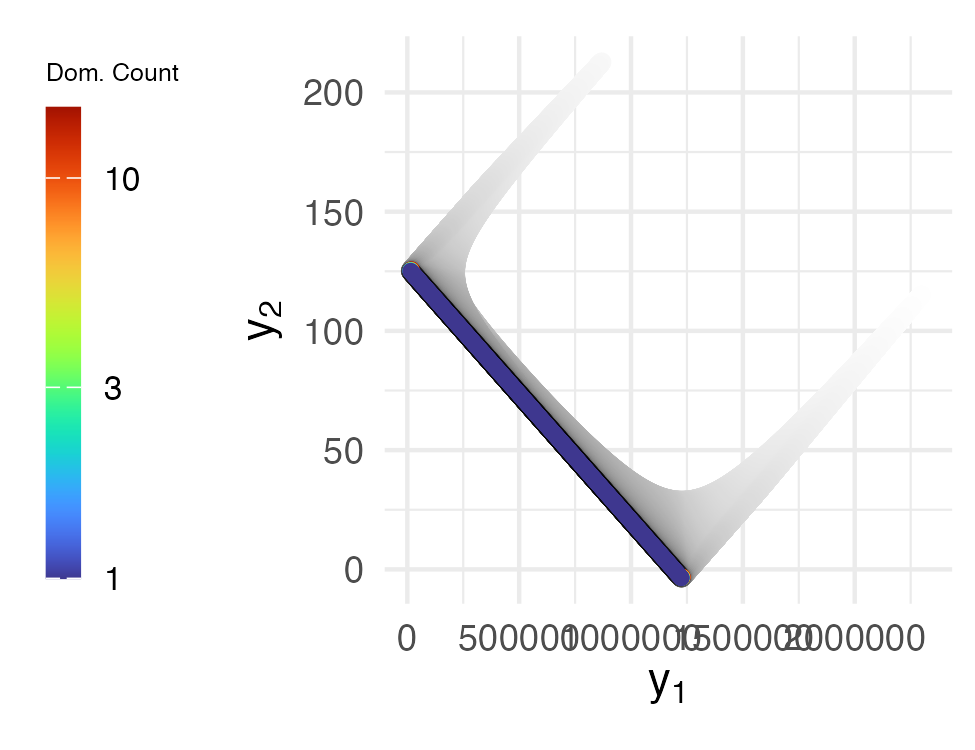}
        \caption{BONO4 instance: Linear ellipsoids.}
    \end{subfigure}
    \begin{subfigure}{\textwidth}
        \includegraphics[width=0.245\linewidth]{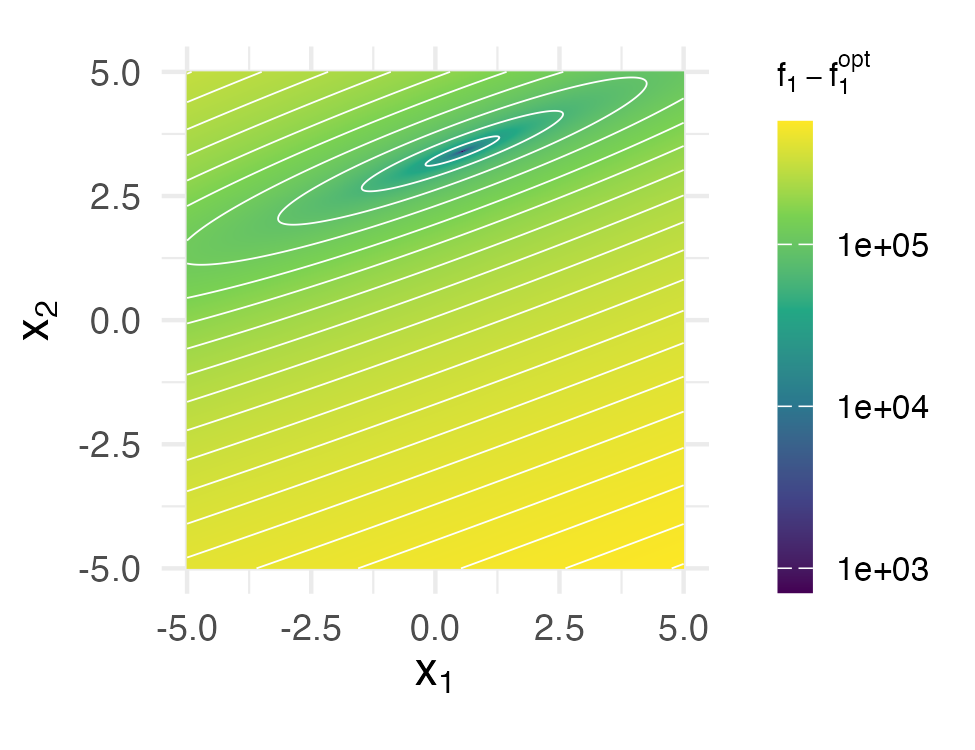}
        \includegraphics[width=0.245\linewidth]{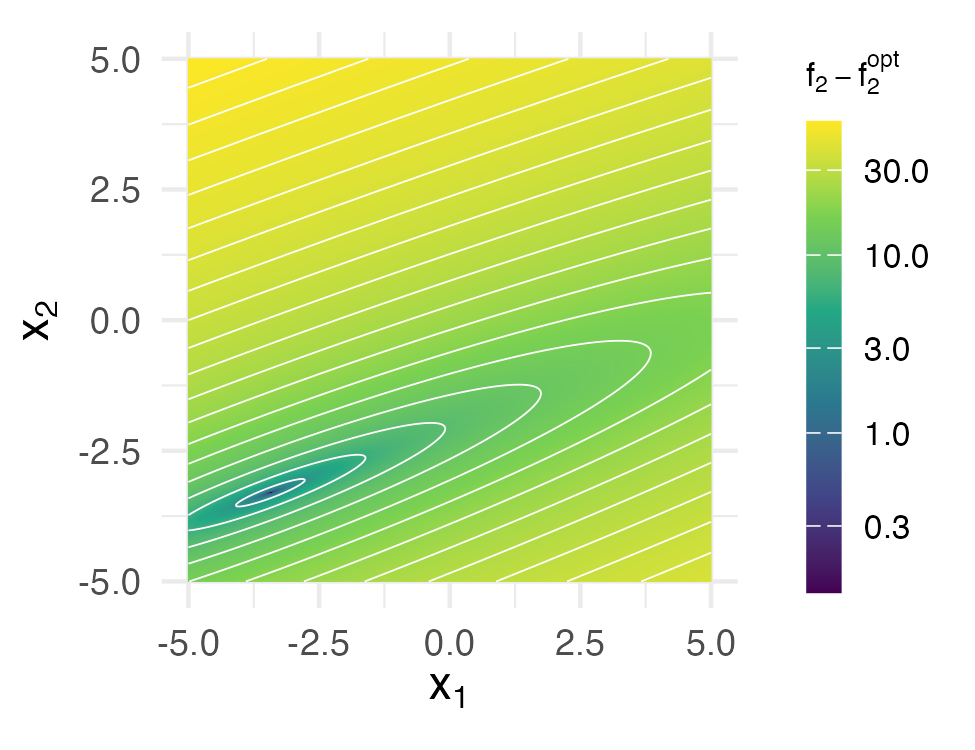}
        \includegraphics[width=0.245\linewidth]{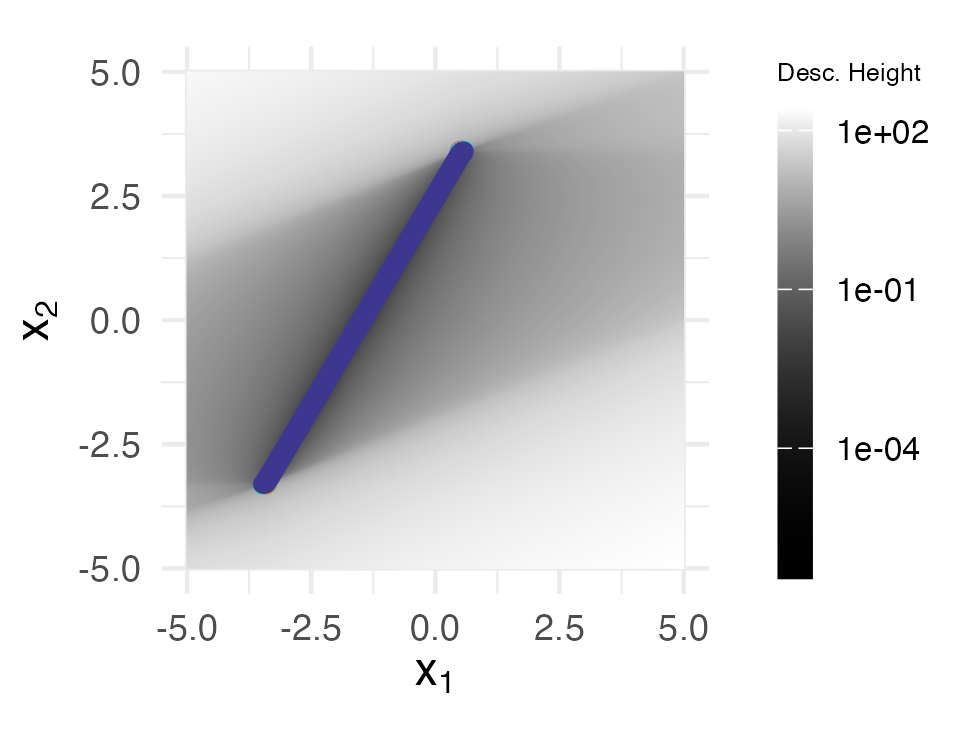}
        \includegraphics[width=0.245\linewidth]{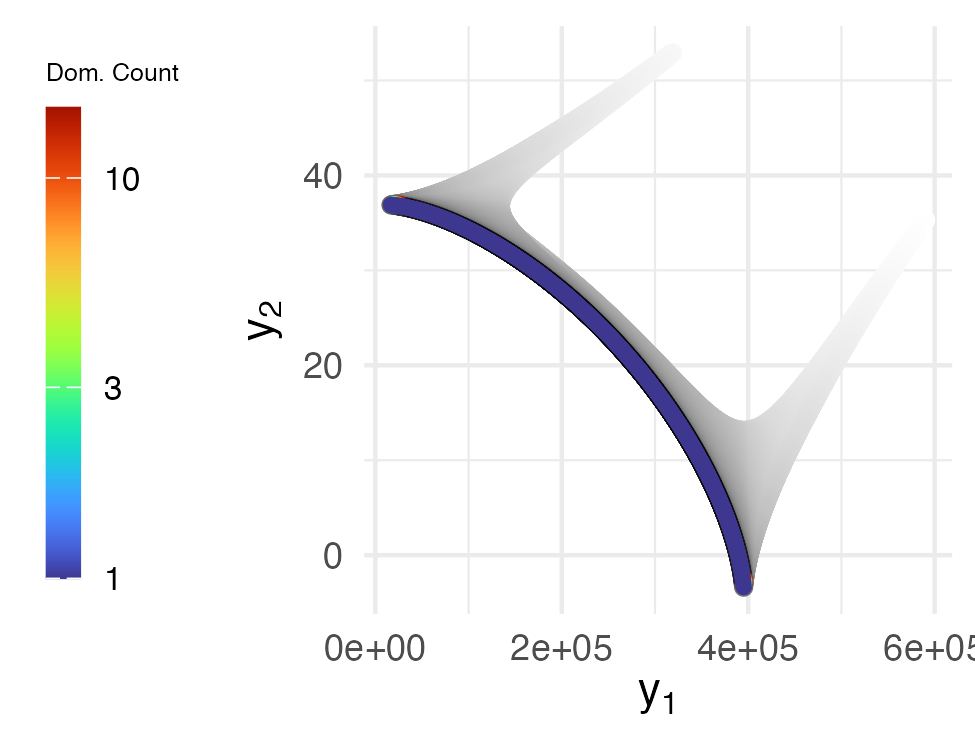}
        \caption{BONO5 instance: Concave ellipsoids.}
    \end{subfigure}
    \caption{Ellipsoids with linear Pareto set, but different Pareto front shapes. The particular distance parameters sampled are $p \approx 2.89$, $p = 1$, and $p \approx 0.64$ for the convex, linear and concave fronts, respectively.}
    \label{fig:linear-ps}
\end{figure}

\medskip

Following this, we have three problem categories with moderately conditioned ellipsoids, i.e., $\kappa \sim \LU(50,200)$.
These still feature linear Pareto sets by sharing identical Hessian matrices ($H_1 = H_2$), but the objective functions are not separable anymore.
All ellipsoid generators with linear Pareto sets are illustrated in \Cref{fig:linear-ps}.

\subsubsection{BONO3: Convex-front Ellipsoids}

This problem class features convex-front problems, where we sample $p \in \LU (1.5,3)$.

\subsubsection{BONO4: Linear-front Ellipsoids}

Problems in this category are identical to the convex-front ellipsoids, but use $p = 1$ to enforce the creation of a linear Pareto front.

\subsubsection{BONO5: Concave-front Ellipsoids}

This is the final category of problems with linear Pareto sets.
Here, we sample $p \sim \LU(\frac 1 3, \frac 2 3)$ to create a concave Pareto front.

\begin{figure}[t]
    \centering
    \begin{subfigure}{\textwidth}
        \includegraphics[width=0.245\linewidth]{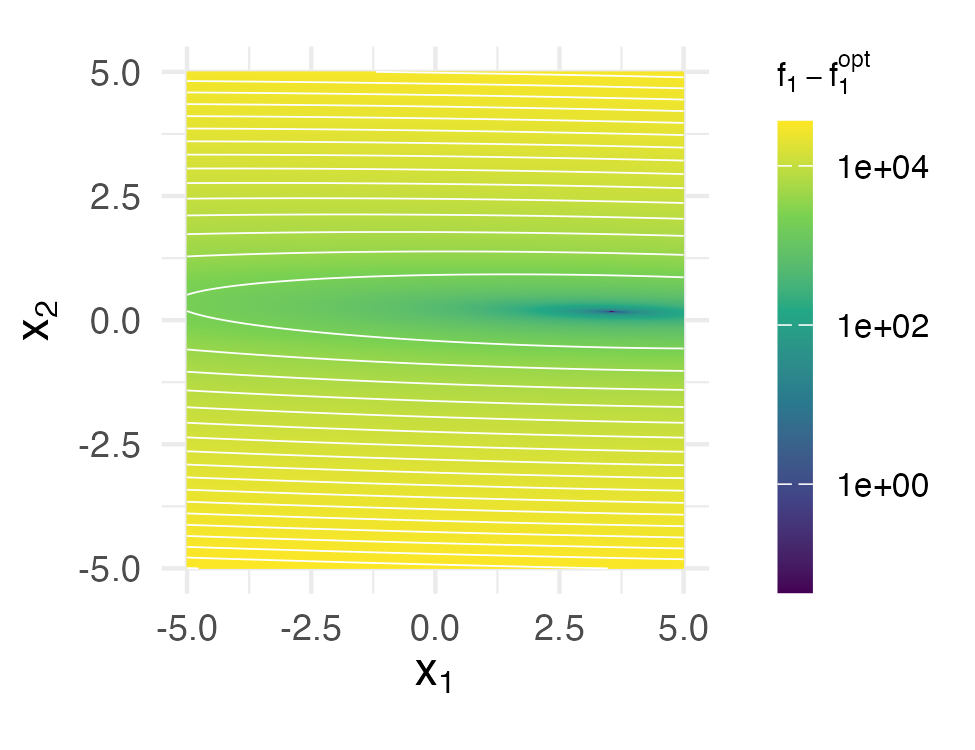}
        \includegraphics[width=0.245\linewidth]{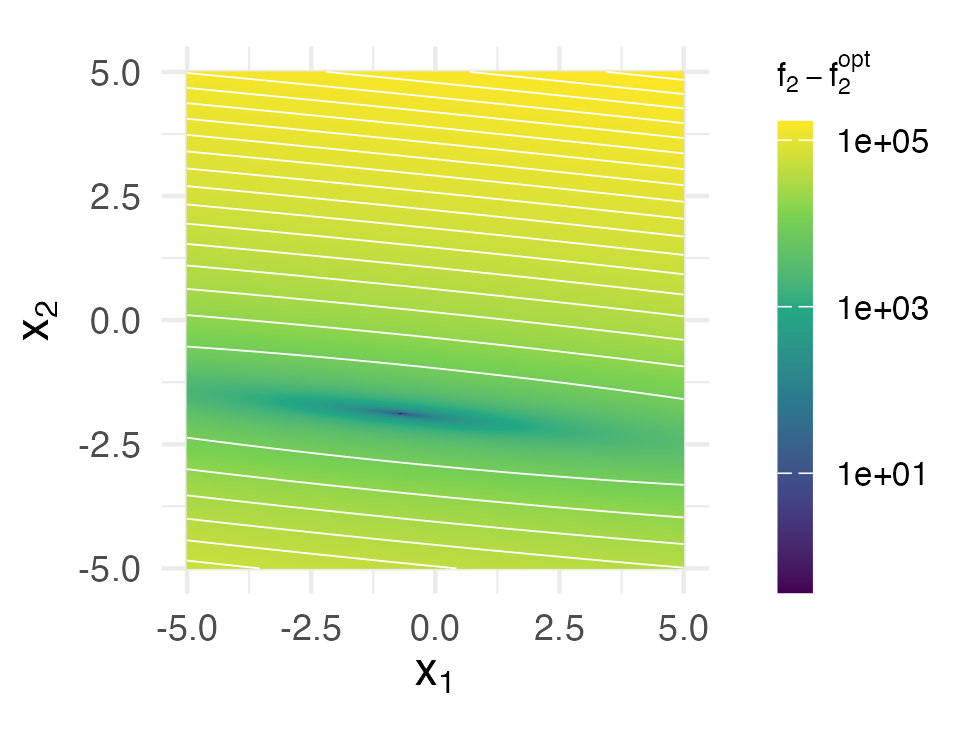}
        \includegraphics[width=0.245\linewidth]{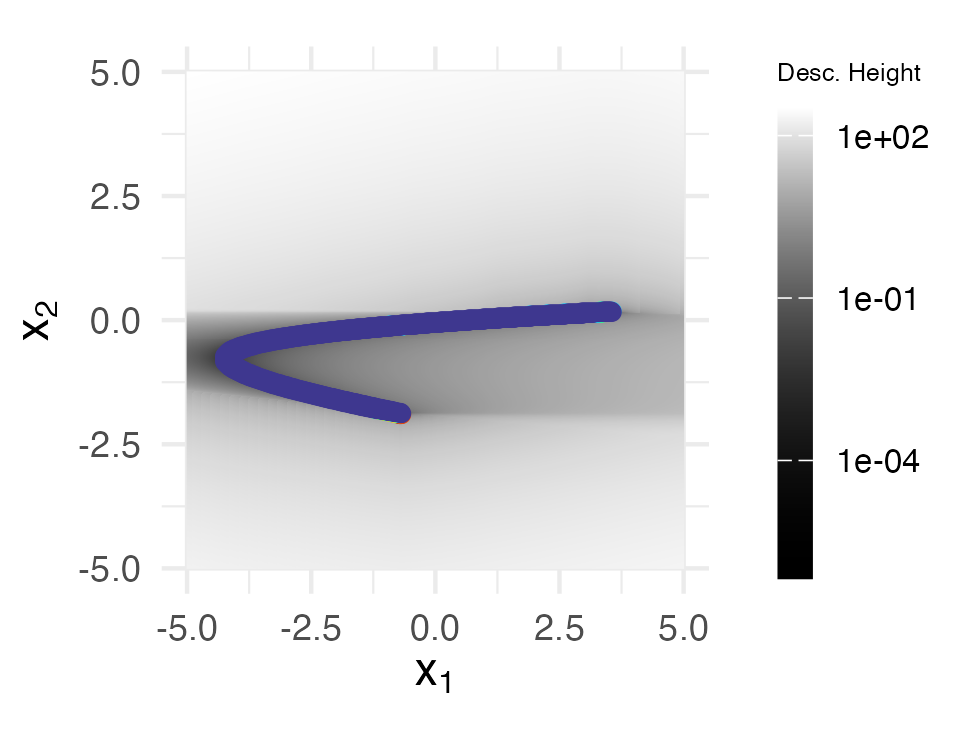}
        \includegraphics[width=0.245\linewidth]{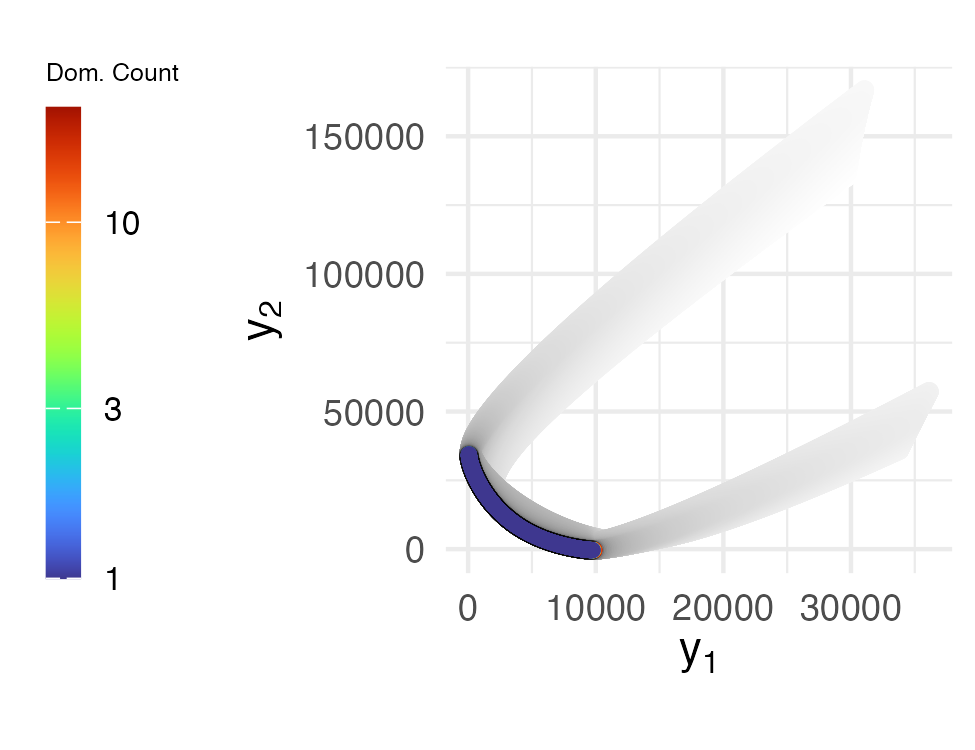}
        \caption{BONO6 instance: Free ellipsoids.}
        \label{fig:free-ellipsoids:a}
    \end{subfigure}
    \begin{subfigure}{\textwidth}
        \includegraphics[width=0.245\linewidth]{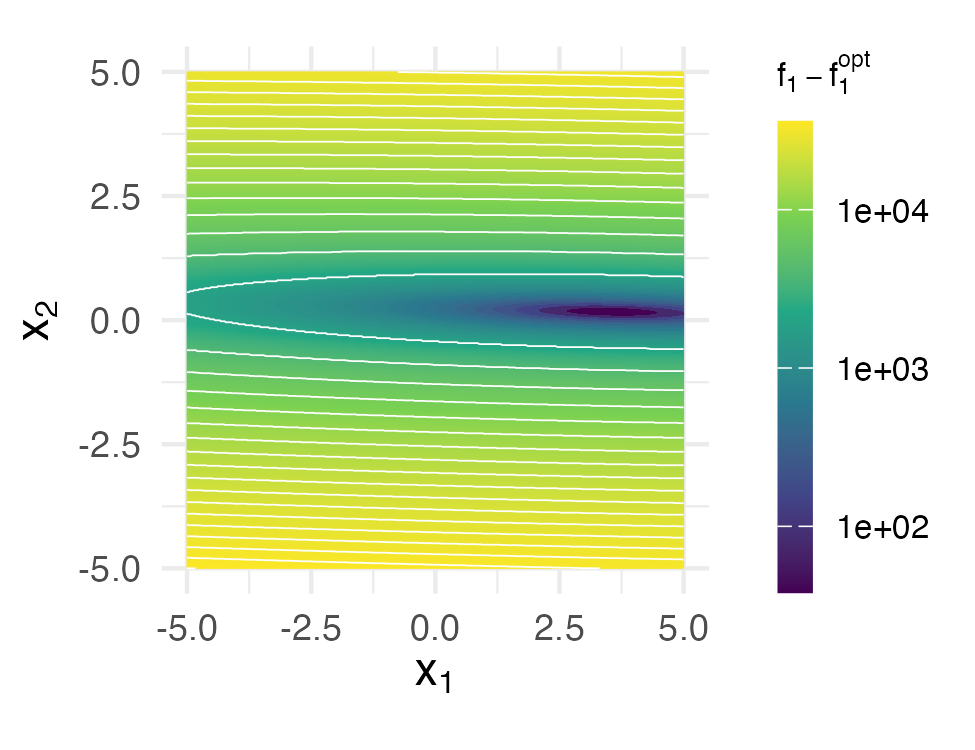}
        \includegraphics[width=0.245\linewidth]{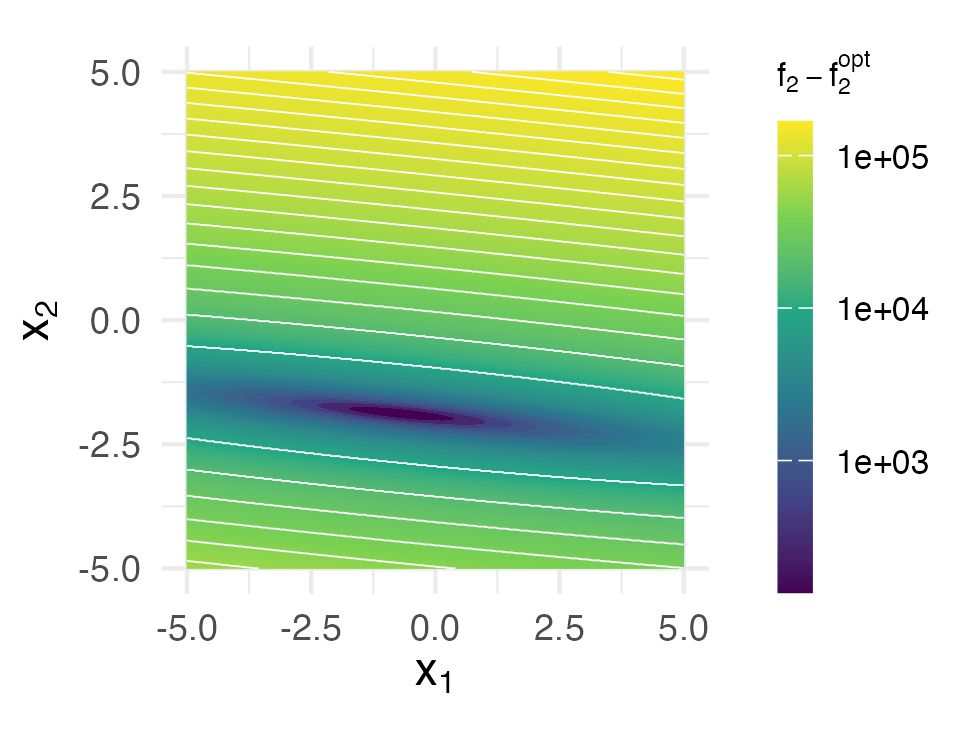}
        \includegraphics[width=0.245\linewidth]{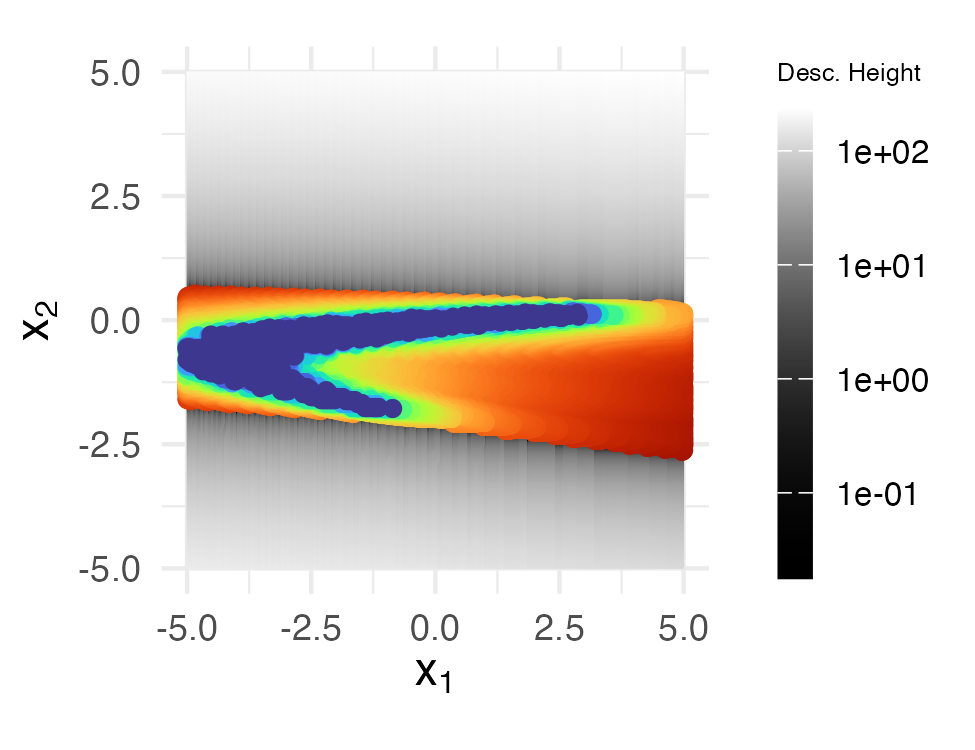}
        \includegraphics[width=0.245\linewidth]{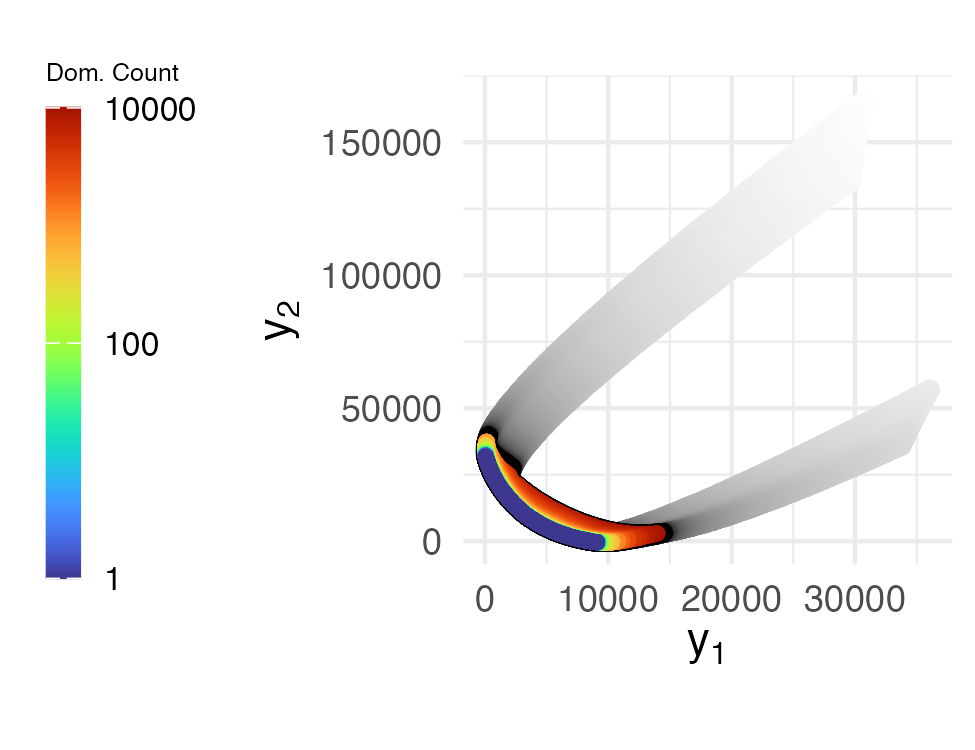}
        \caption{BONO7 instance: Stepped ellipsoids.}
        \label{fig:free-ellipsoids:b}
    \end{subfigure}
    
    \caption{Illustration of two optimization problems with independently sampled Hessians. While (a) is smooth, (b) features objective discretization, introducing many plateaus in both objectives and the resulting multi-objective problem.}
\end{figure}

\begin{figure}[t]
    \label{fig:free-ellipsoids}
    \medskip
    \begin{subfigure}{\textwidth}
        \includegraphics[width=0.245\linewidth]{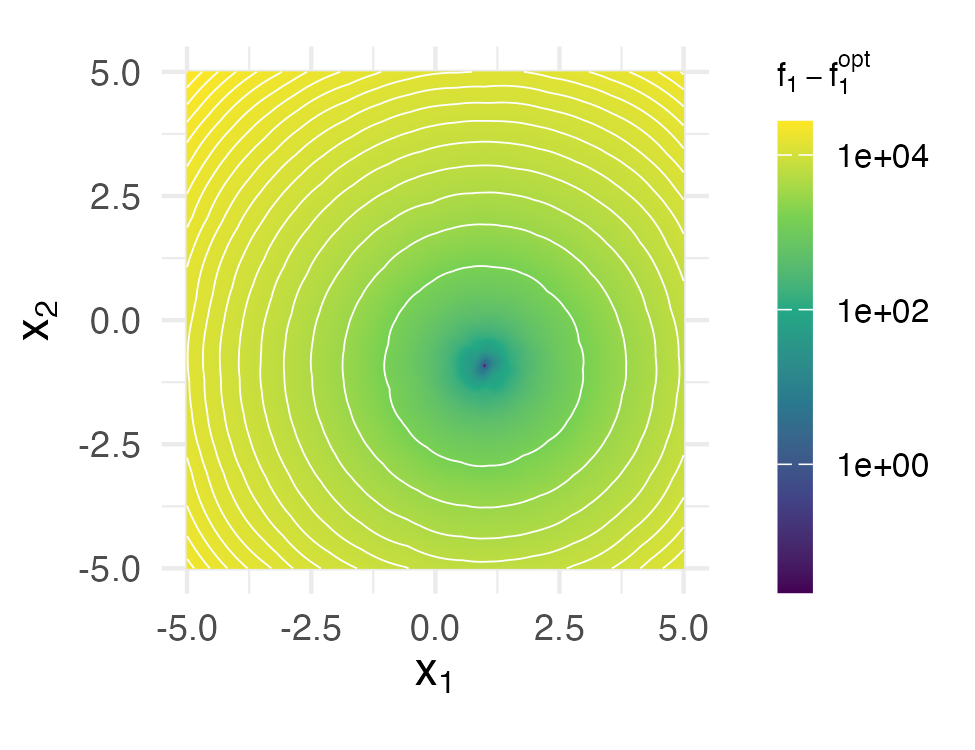}
        \includegraphics[width=0.245\linewidth]{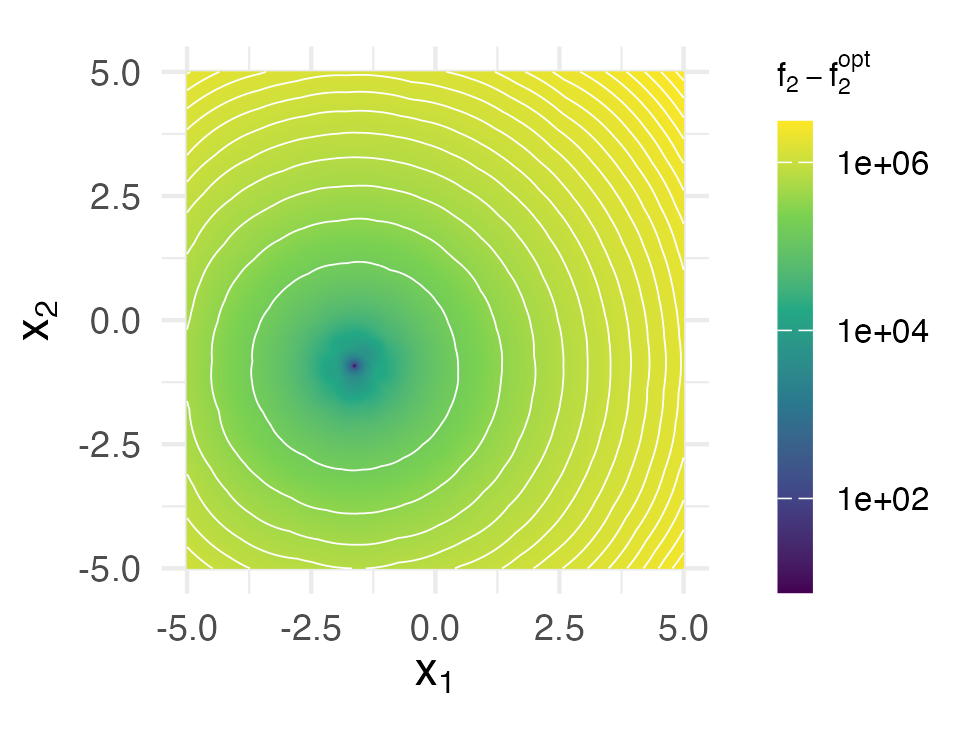}
        \includegraphics[width=0.245\linewidth]{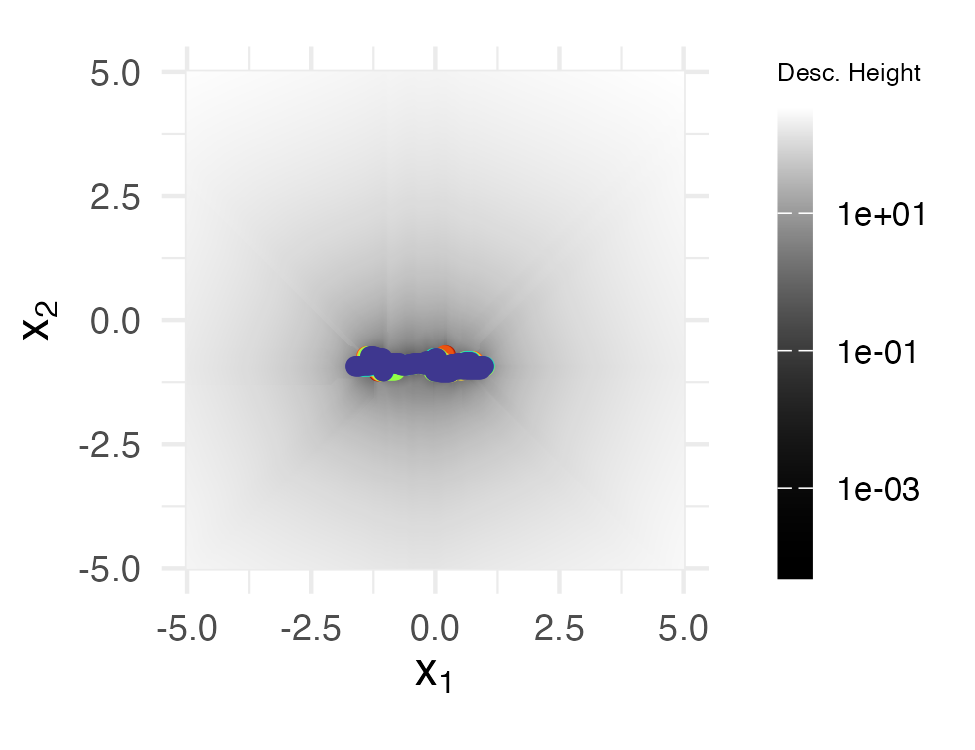}
        \includegraphics[width=0.245\linewidth]{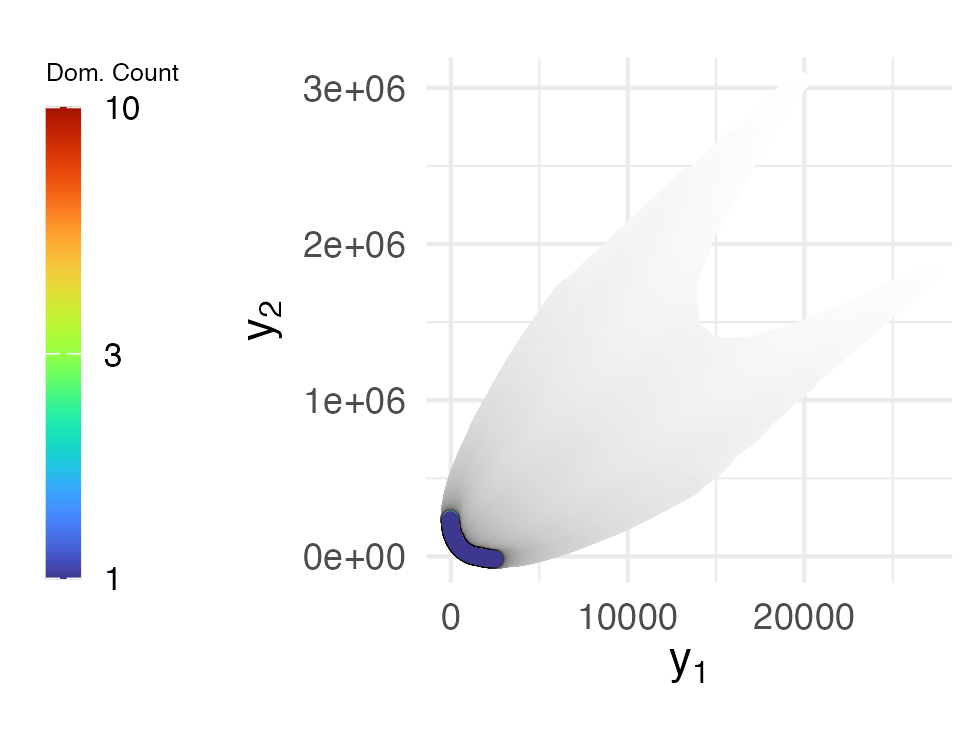}
        \caption{BONO8 instance: Multimodal axis-aligned spheres.}
    \end{subfigure}
    \begin{subfigure}{\textwidth}
        \includegraphics[width=0.245\linewidth]{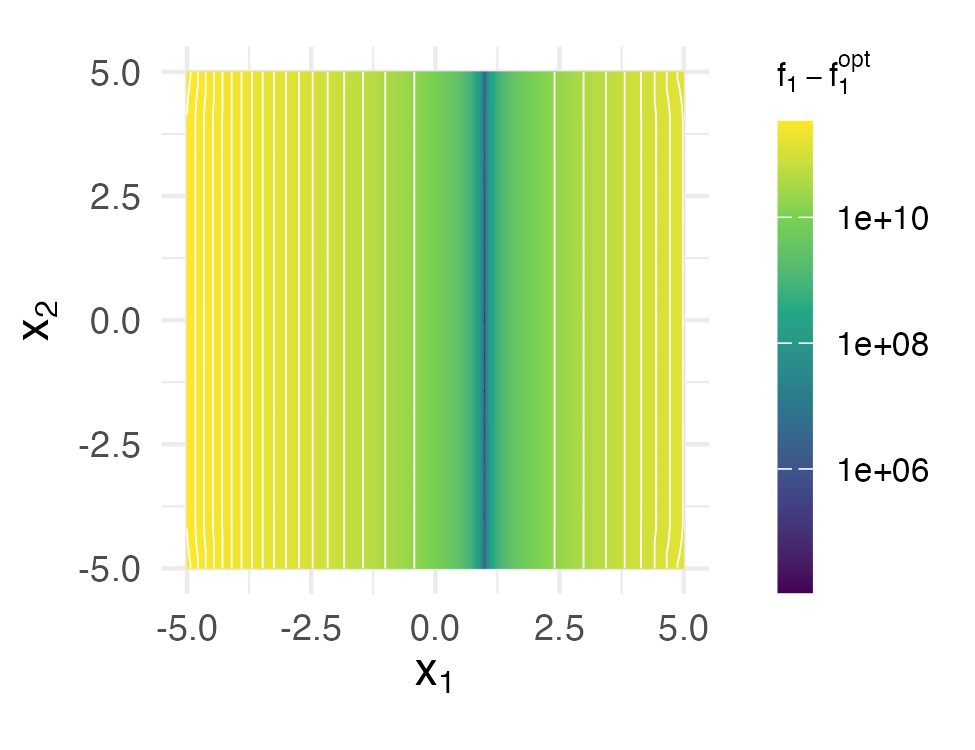}
        \includegraphics[width=0.245\linewidth]{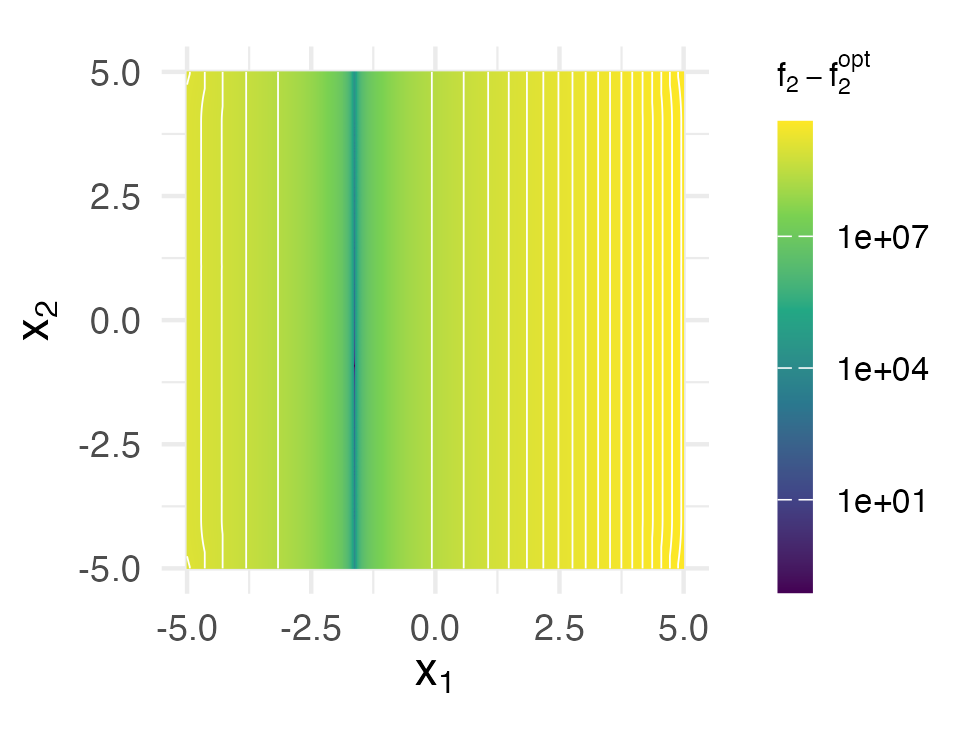}
        \includegraphics[width=0.245\linewidth]{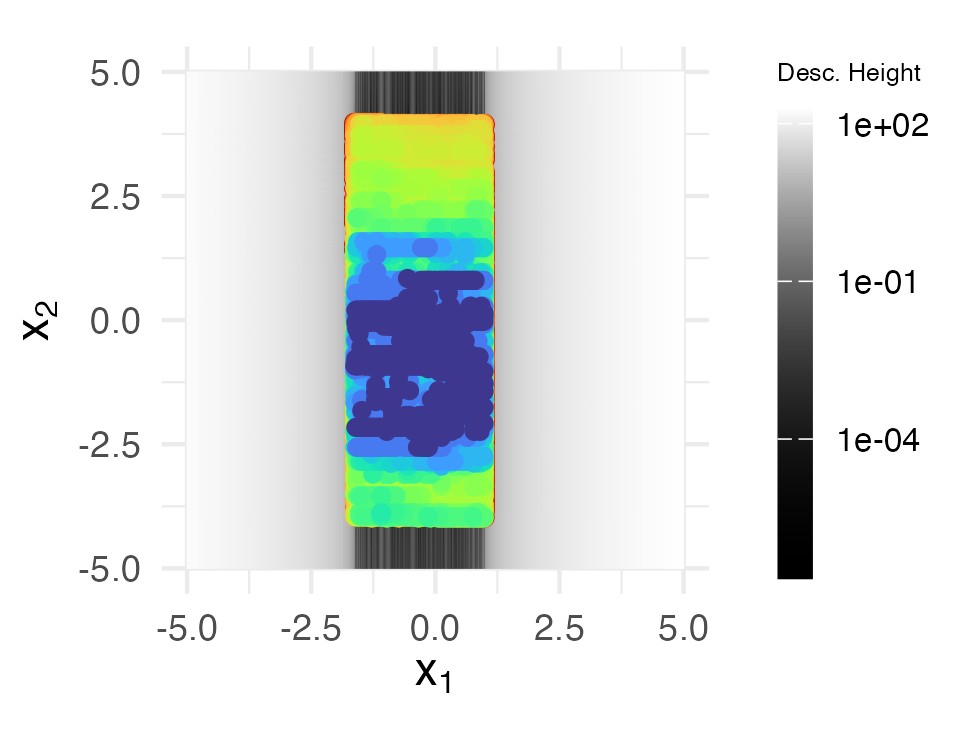}
        \includegraphics[width=0.245\linewidth]{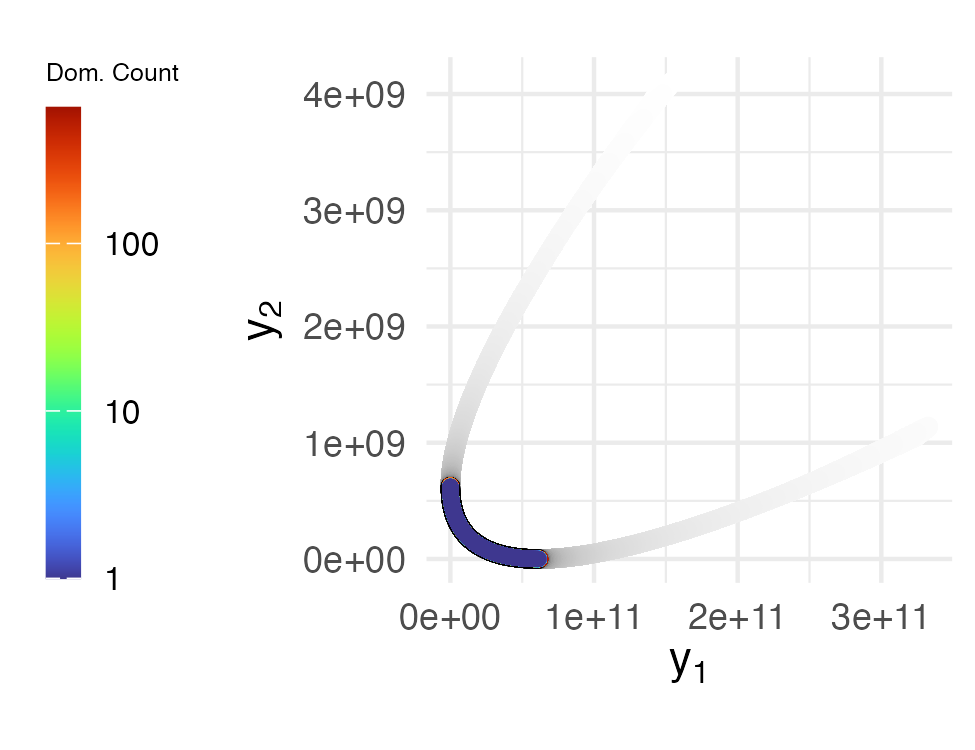}
        \caption{BONO9 instance: Multimodal axis-aligned ellipsoids.}
    \end{subfigure}
    
    \caption{Visualization of instances of the multimodal, axis-aligned problems. While their global shape follows BONO1-2, the Pareto sets and front become more disconnected.}
    \label{fig:mm-axis-aligned}
\end{figure}

\medskip

At last, we have two classes of unimodal problems without a linear Pareto set.
Representatives of these two groups are shown in \Cref{fig:free-ellipsoids}.

\subsubsection{BONO6: Free Ellipsoids}
For the group of free ellipsoid problems, we sample two quadratic problems independently.
We use the same moderate conditioning as in the previous problems, i.e., $\kappa \sim \LU(50, 200)$, but the independence of the Hessians creates more complex, curved Pareto sets, cf.~\Cref{fig:free-ellipsoids}.
This also is the first category for which we cannot guarantee from the outset that the unbounded Pareto-optimal solutions lie within the designated decision space $[-5,5]^d$ due to the curved shape of the Pareto set.
For this category, we densely sample along the Pareto set (cf.~\Cref{sec:approximation}) to verify the admissibility of all solutions and, if the check fails, repeat the problem creation until the Pareto set is confirmed to be within bounds.

\subsubsection{BONO7: Stepped Ellipsoids}

Problems in this final unimodal category are identical to free ellipsoids, but use the discretized formulation from \Cref{eq:stepped_peak}.
The number of discretization steps is sampled from $N_h \sim \lfloor\LU(50, 201)\rfloor$ (resulting in values between $50$ and $200$) and used to discretize the area between ideal and nadir points into the given number of steps using step size $h_i = (Y_{N,i} - Y_{I,i}) / N_h$.
This creates a search space with many plateaus, likely increasing optimization difficulty.

\subsection{Multimodal Problems with Global Structure} \label{sec:bono_multi_with}

Next, we use the unimodal BONO problems BONO1-7 to create a set of multimodal problems \emph{with global structure}.
To enable comparisons to the unimodal functions, we use each unimodal generator as the basis for the global shapes of the problems, and apply $J=500$ perturbation functions each, resulting in seven problem categories: BONO8-14.

\subsubsection{BONO8: Multimodal Axis-aligned Spheres}

Both multimodal axis-aligned problem generators are illustrated in \Cref{fig:mm-axis-aligned}.
In BONO8, the perturbation functions consist of sphere functions.

\subsubsection{BONO9: Multimodal Axis-aligned Ellipsoids}

In BONO9, perturbation functions are separable ellipsoids with conditioning equal to the global function to be approximated.

\begin{figure}
    \centering
    \begin{subfigure}{\textwidth}
        \includegraphics[width=0.245\linewidth]{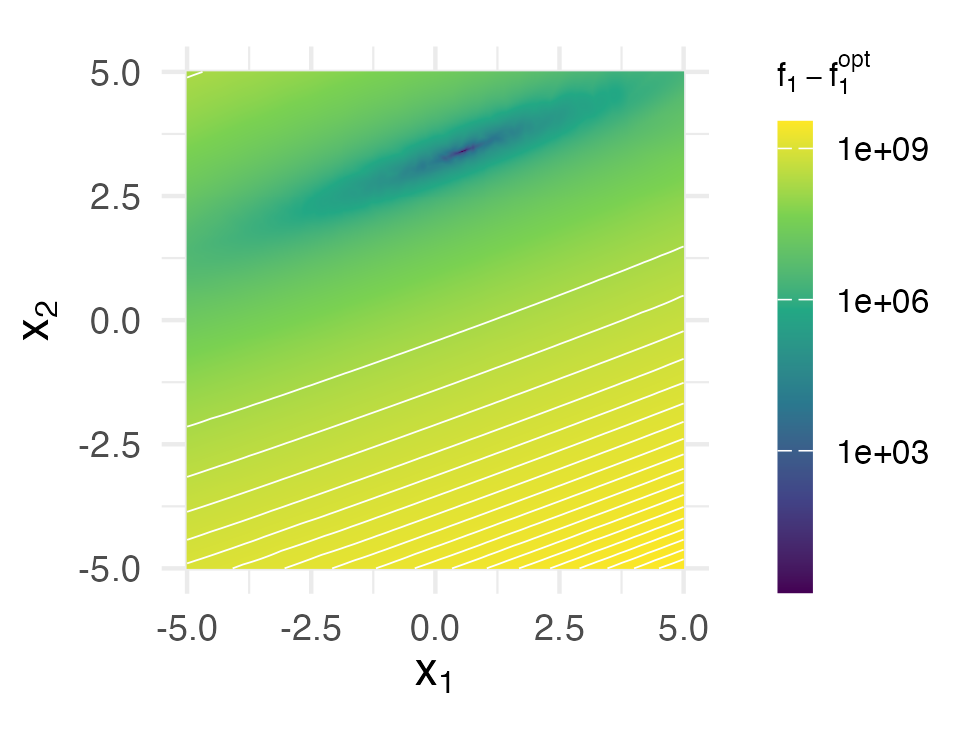}
        \includegraphics[width=0.245\linewidth]{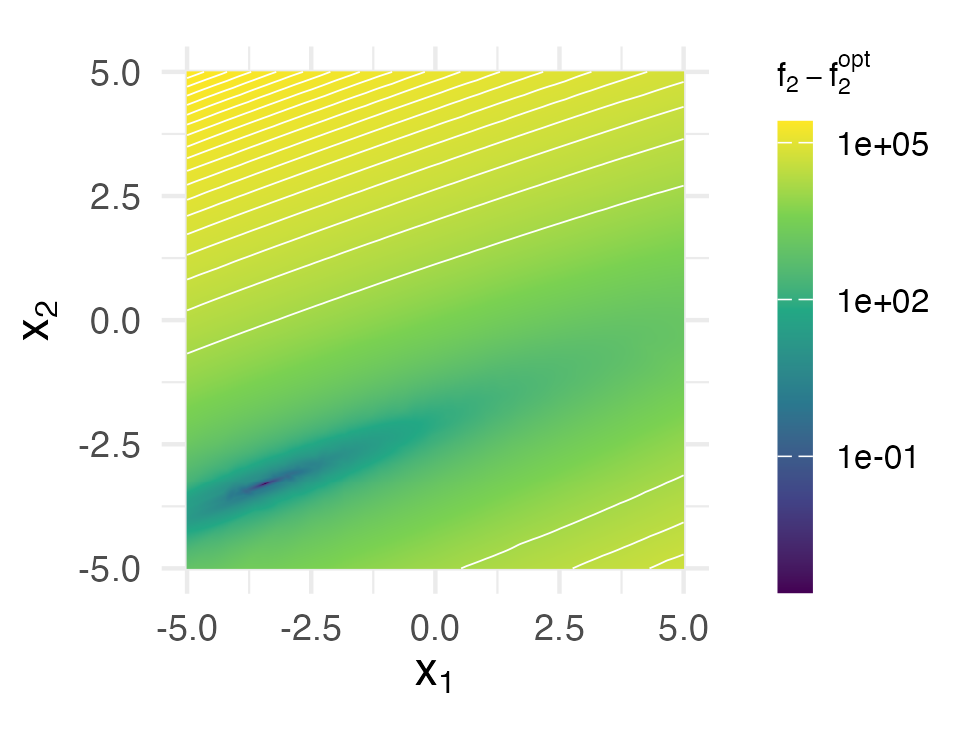}
        \includegraphics[width=0.245\linewidth]{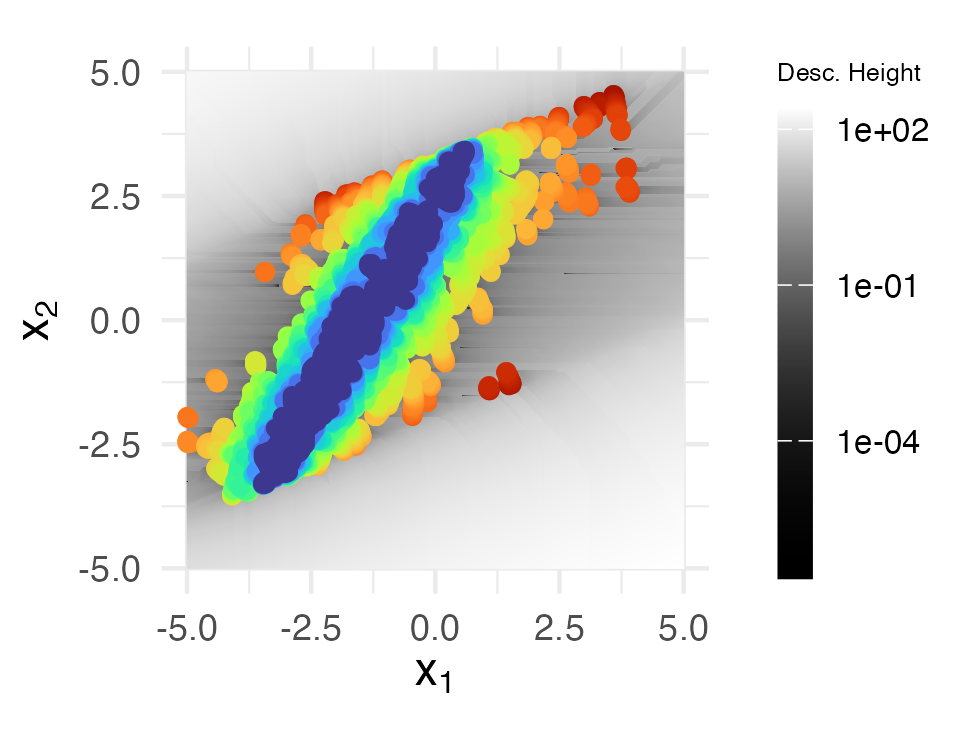}
        \includegraphics[width=0.245\linewidth]{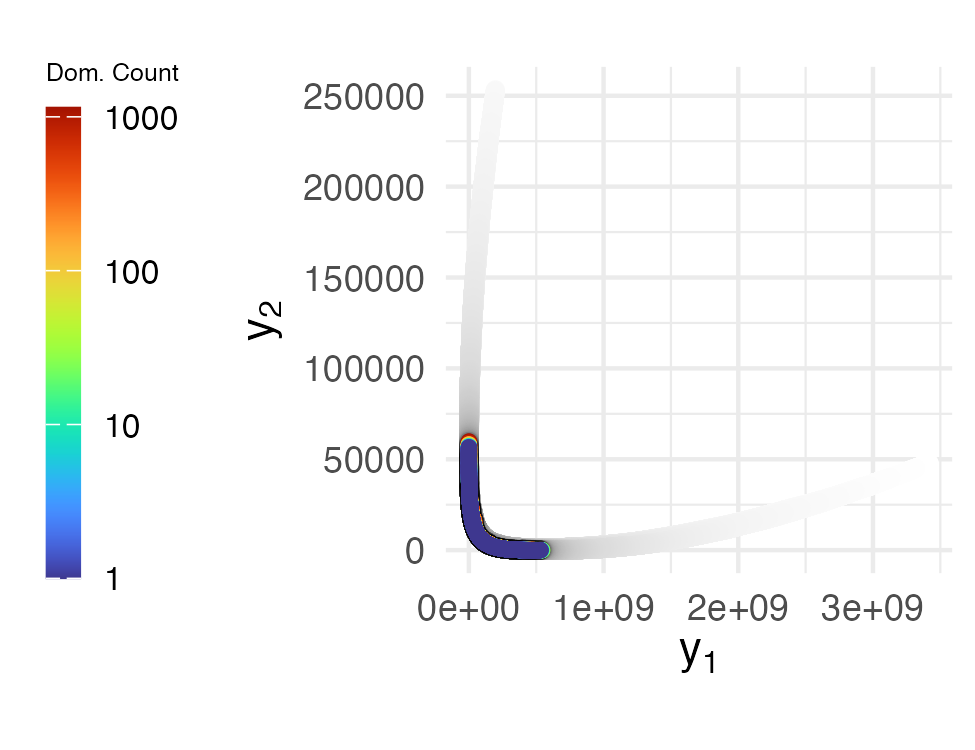}
        \caption{BONO10 instance: Convex ellipsoid global structure.}
    \end{subfigure}
    \begin{subfigure}{\textwidth}
        \includegraphics[width=0.245\linewidth]{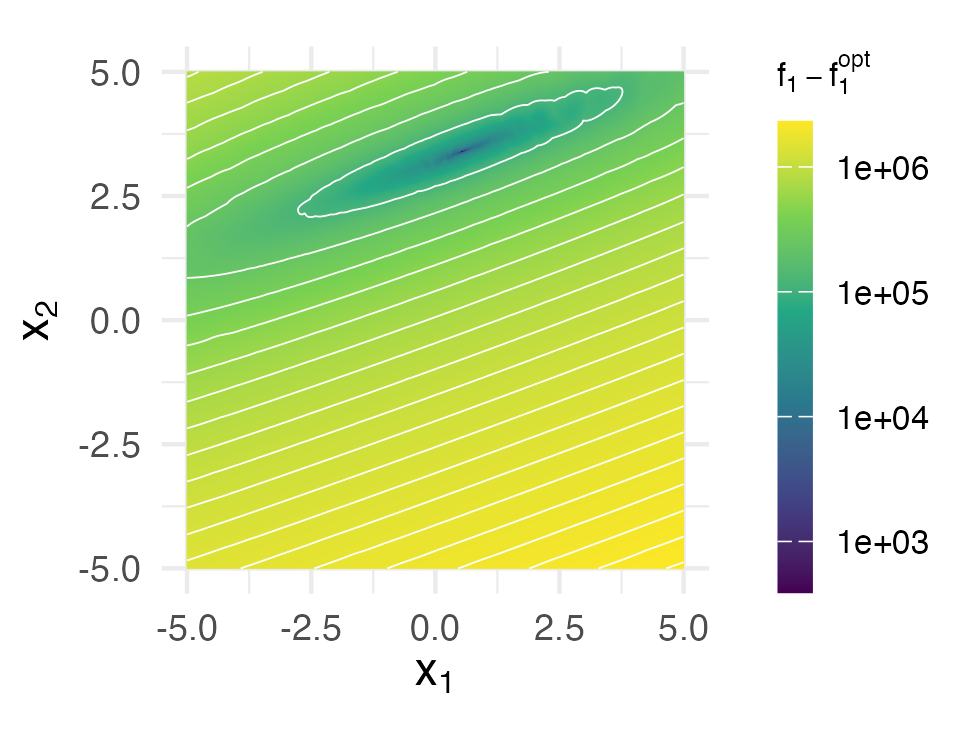}
        \includegraphics[width=0.245\linewidth]{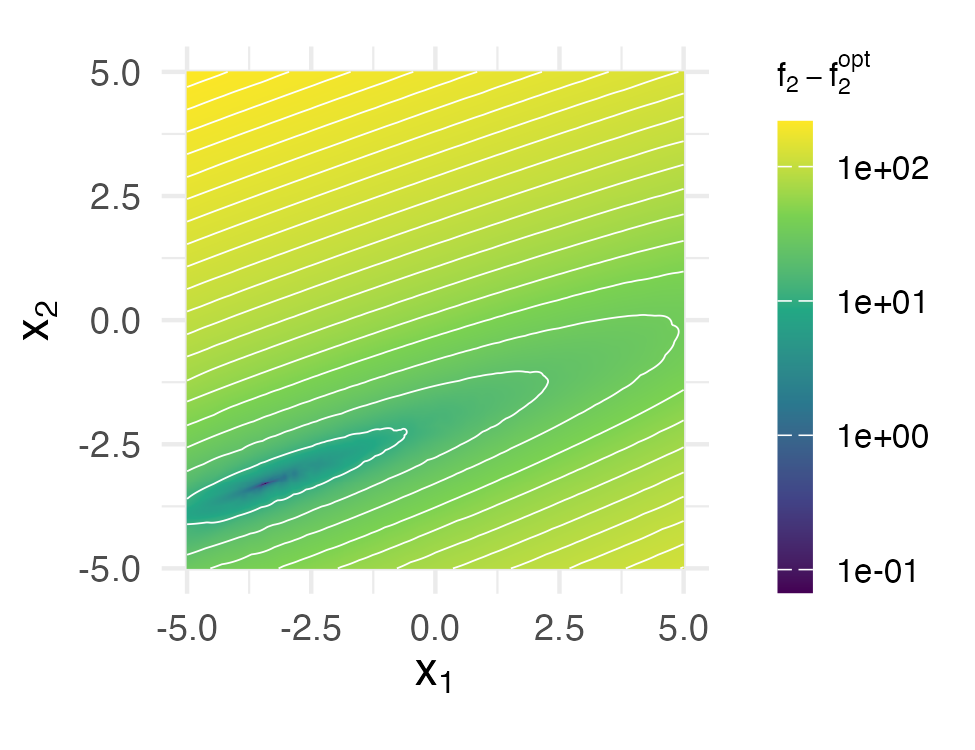}
        \includegraphics[width=0.245\linewidth]{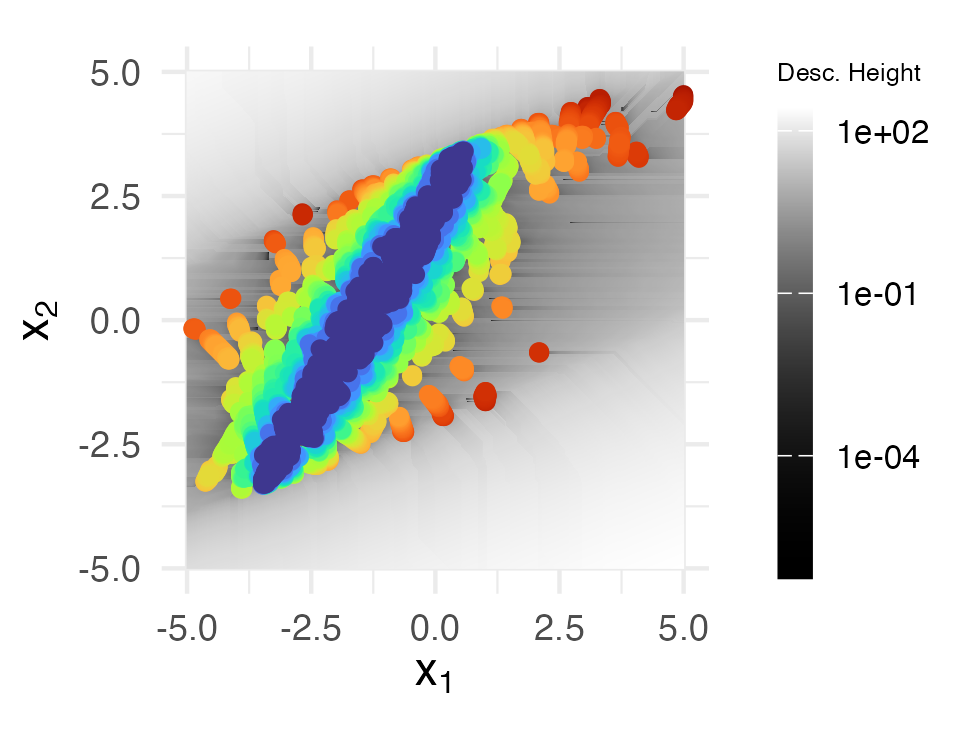}
        \includegraphics[width=0.245\linewidth]{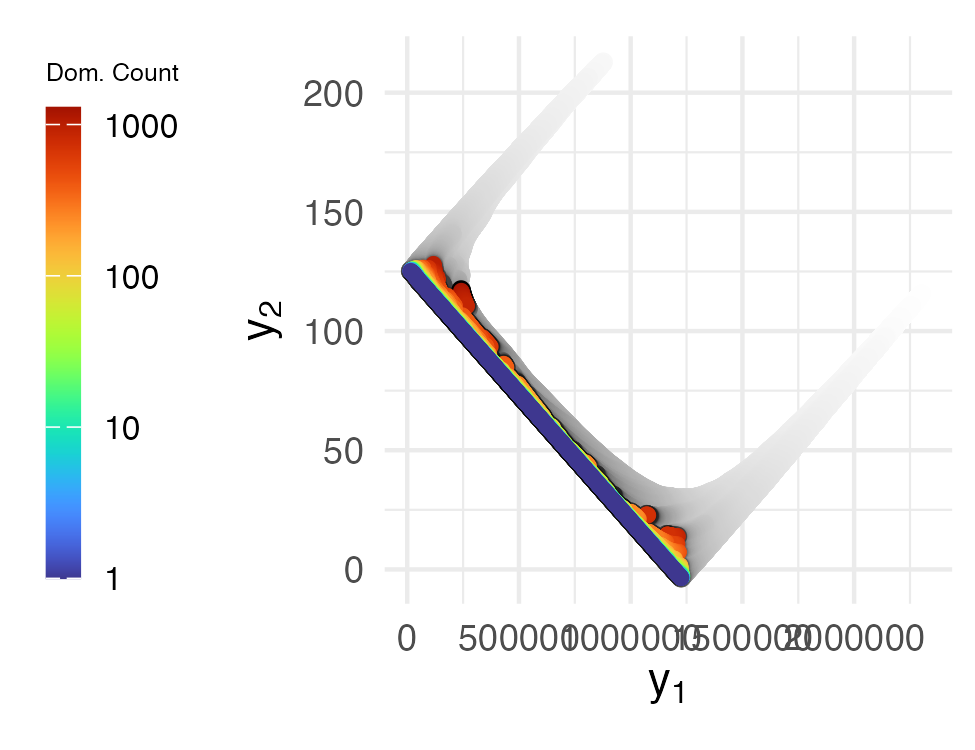}
        \caption{BONO11 instance: Linear ellipsoid global structure.}
    \end{subfigure}
    \begin{subfigure}{\textwidth}
        \includegraphics[width=0.245\linewidth]{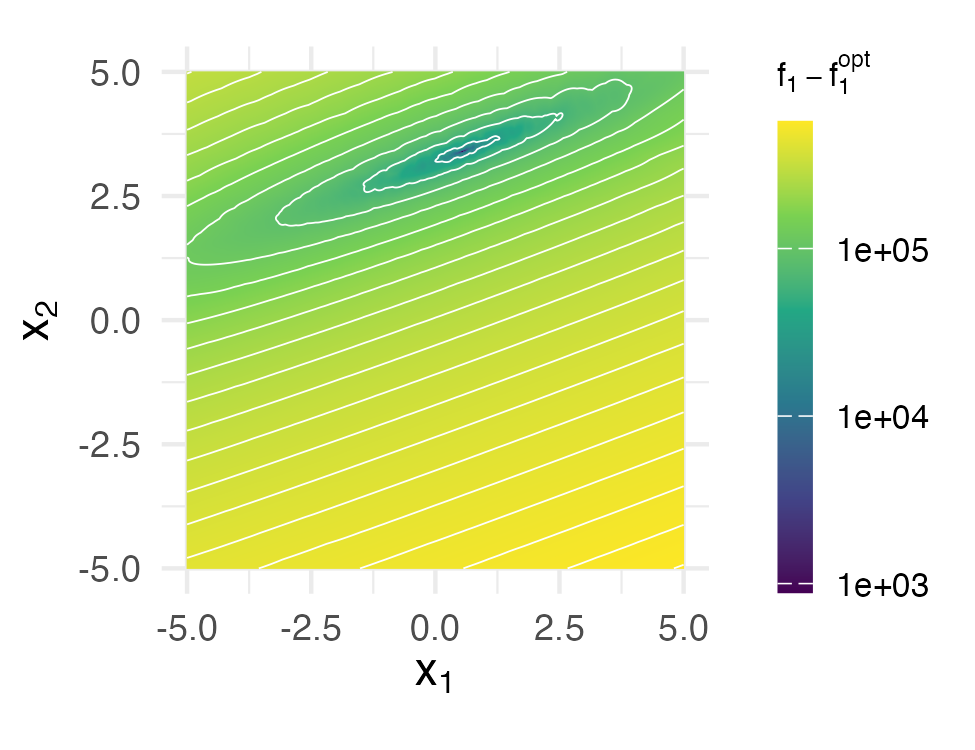}
        \includegraphics[width=0.245\linewidth]{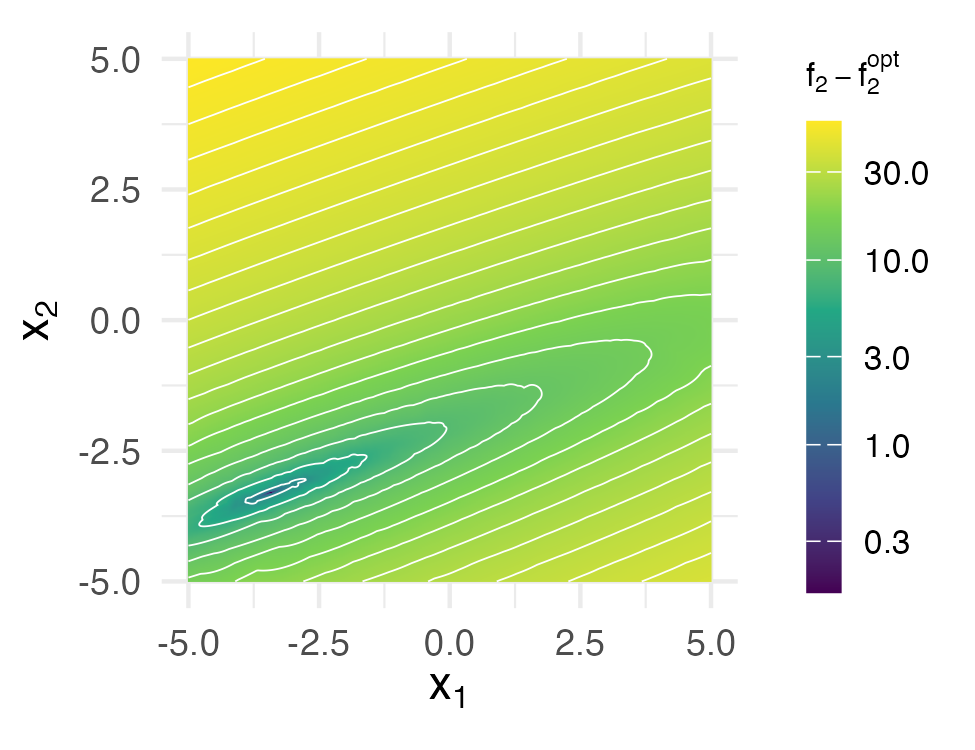}
        \includegraphics[width=0.245\linewidth]{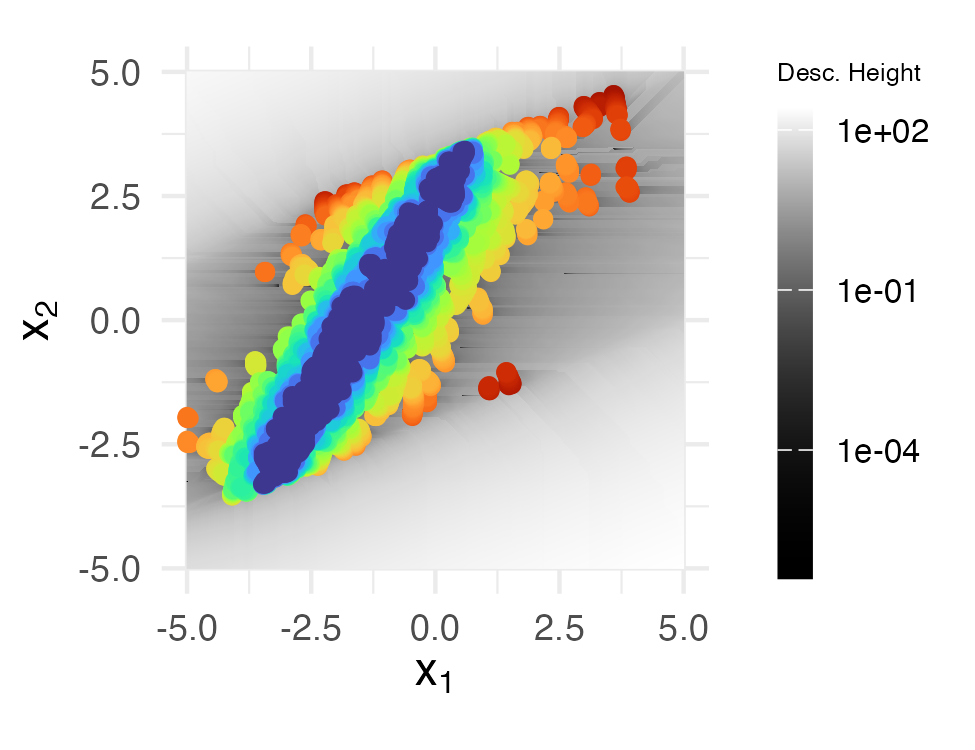}
        \includegraphics[width=0.245\linewidth]{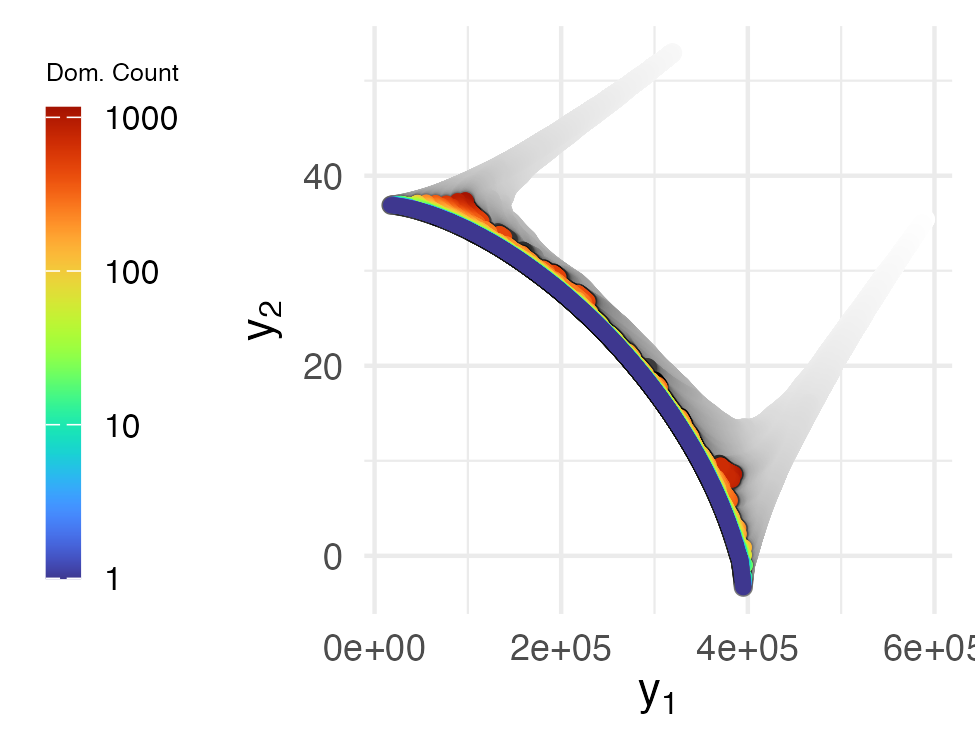}
        \caption{BONO12 instance: Concave ellipsoid global structure.}
    \end{subfigure}
    \caption{Multimodal ellipsoids with linear global Pareto set shape. The Pareto set consists of many local optima, while many additional local optima make the exact approximation more difficult.}
    \label{fig:mm-linear-ps}
\end{figure}

\subsubsection{BONO10: Multimodal Convex-front Ellipsoids}

Exemplary problems for multimodal ellipsoids with linear Pareto sets are illustrated in \Cref{fig:mm-linear-ps}.
As the conditioning of the original problems is $\kappa \sim \LU(50,200)$, Hessian matrices of the perturbation functions will likewise have lower conditioning here as in, e.g., BONO9.
It can be observed that many further, locally efficient points are created in the vicinity of the Pareto set, while the individual objective problems are well-approximated.

\subsubsection{BONO11: Multimodal Linear-front Ellipsoids}

The problems created by the multimodal convex, linear and concave front generators only differ in the approximate Pareto front shape in the objective space.
The multimodal linear front generator demonstrates best that the approximation of the underlying unimodal problem's Pareto front is very close.

\subsubsection{BONO12: Multimodal Concave-front Ellipsoids}

The concave front problem is again the last in this category.
Note that due to the multimodality, neither front shape is guaranteed to be perfect anymore, e.g., a convex front may contain concave parts and vice-versa.

\begin{figure}
    \centering
    \begin{subfigure}{\textwidth}
        \includegraphics[width=0.245\linewidth]{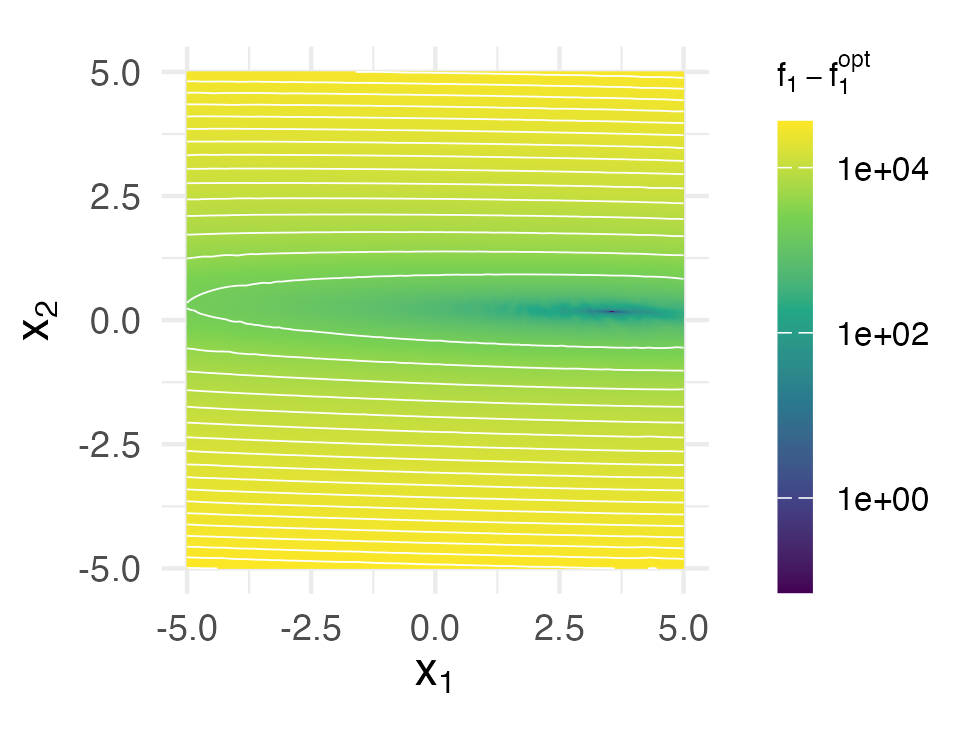}
        \includegraphics[width=0.245\linewidth]{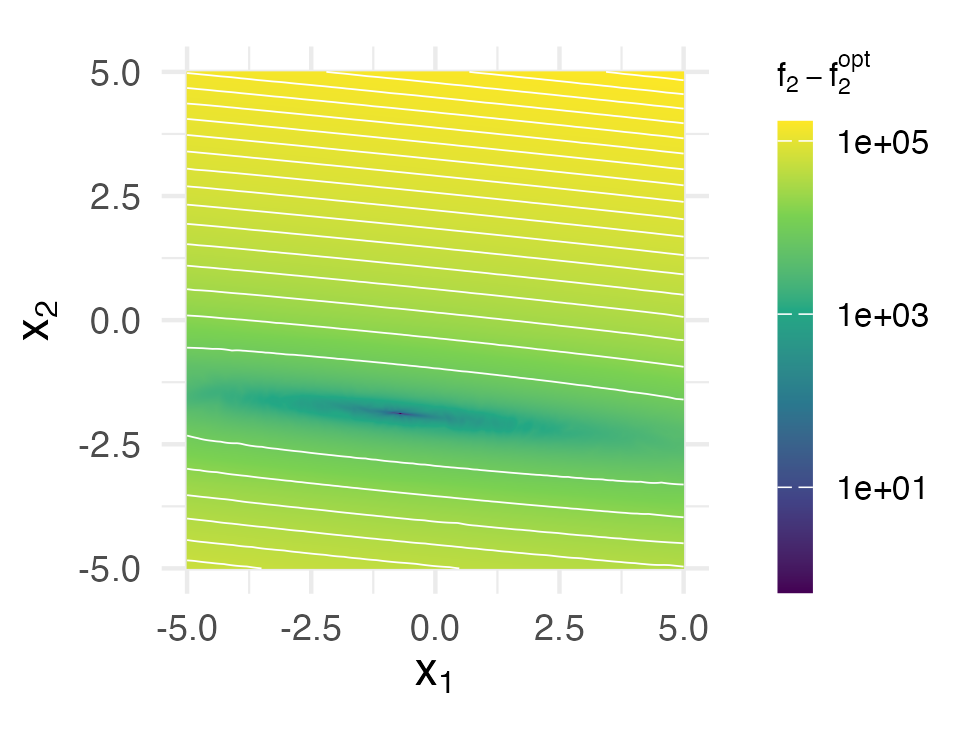}
        \includegraphics[width=0.245\linewidth]{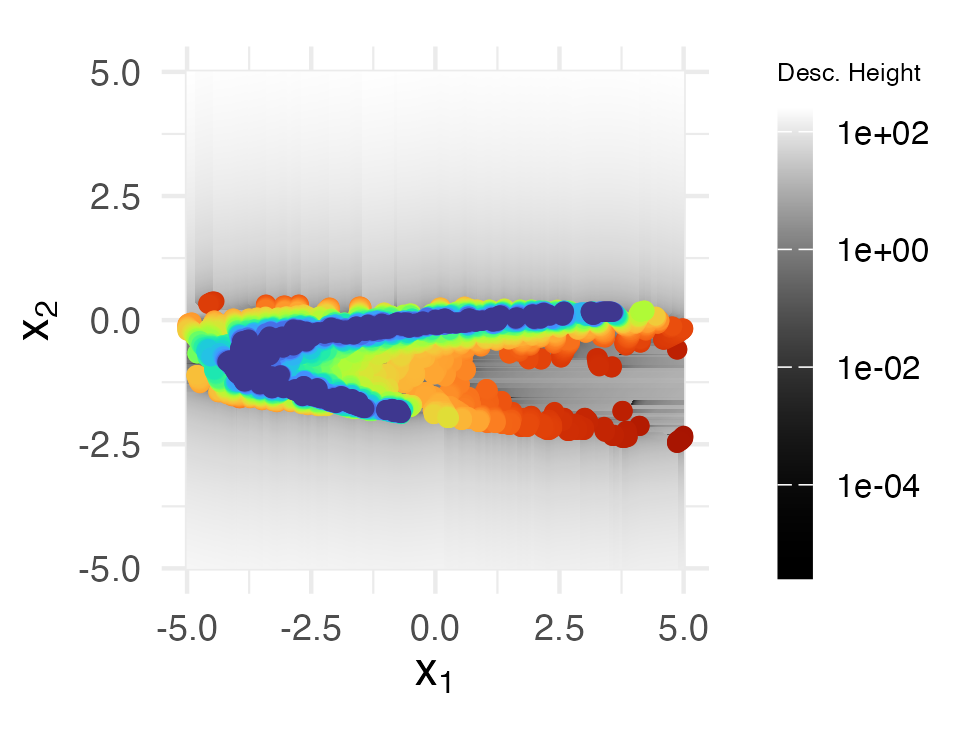}
        \includegraphics[width=0.245\linewidth]{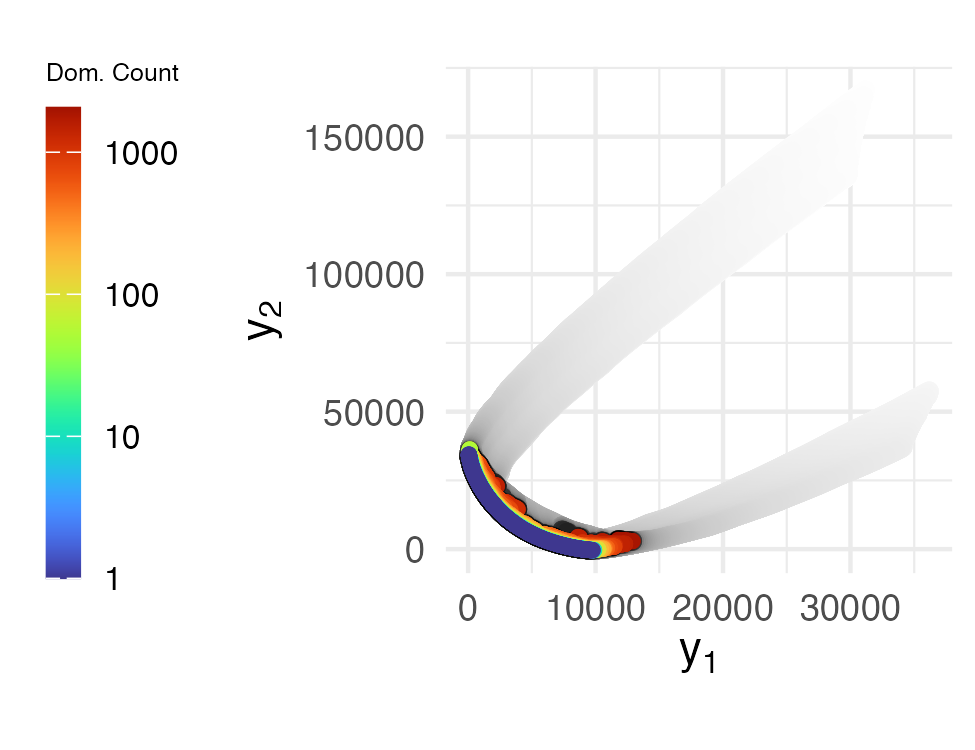}
        \caption{BONO13 instance: Multimodal free ellipsoids.}
    \end{subfigure}
    \begin{subfigure}{\textwidth}
        \includegraphics[width=0.245\linewidth]{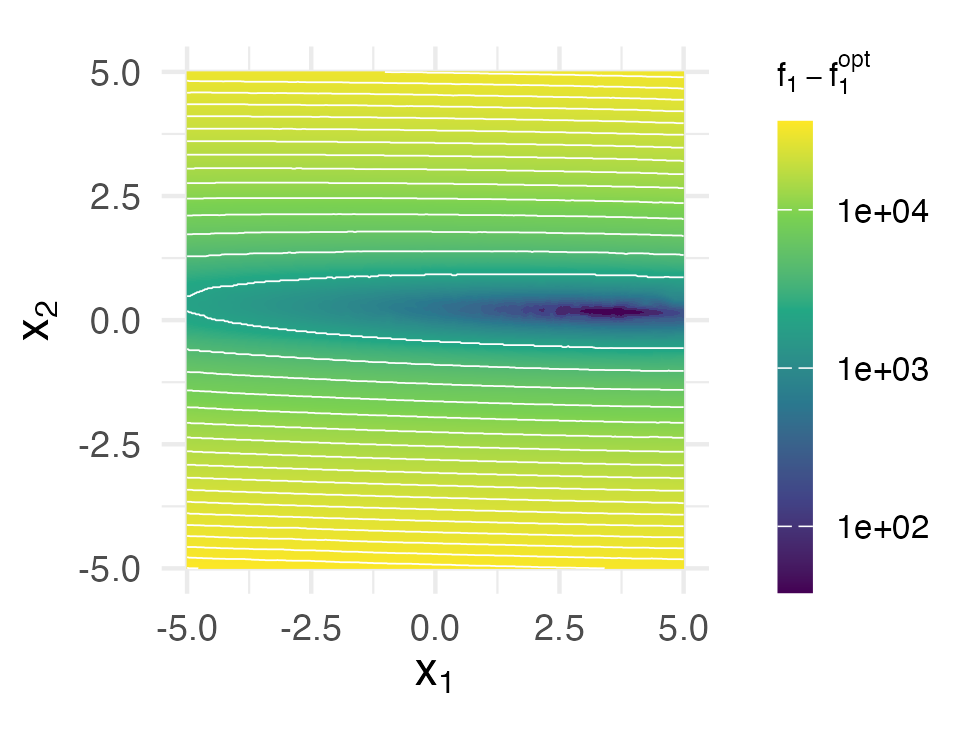}
        \includegraphics[width=0.245\linewidth]{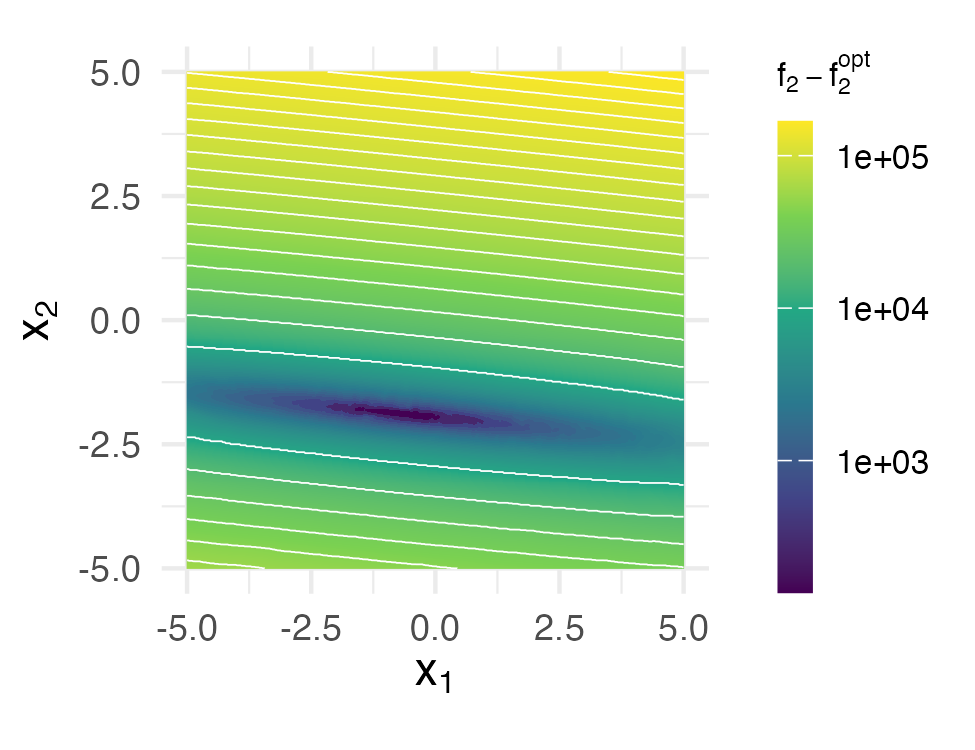}
        \includegraphics[width=0.245\linewidth]{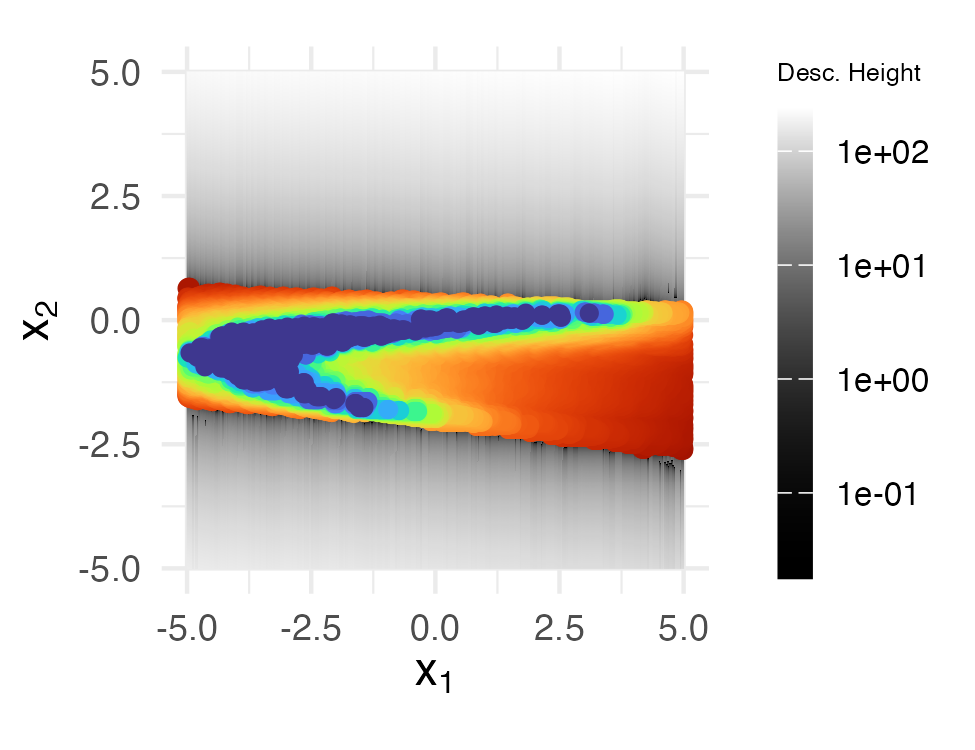}
        \includegraphics[width=0.245\linewidth]{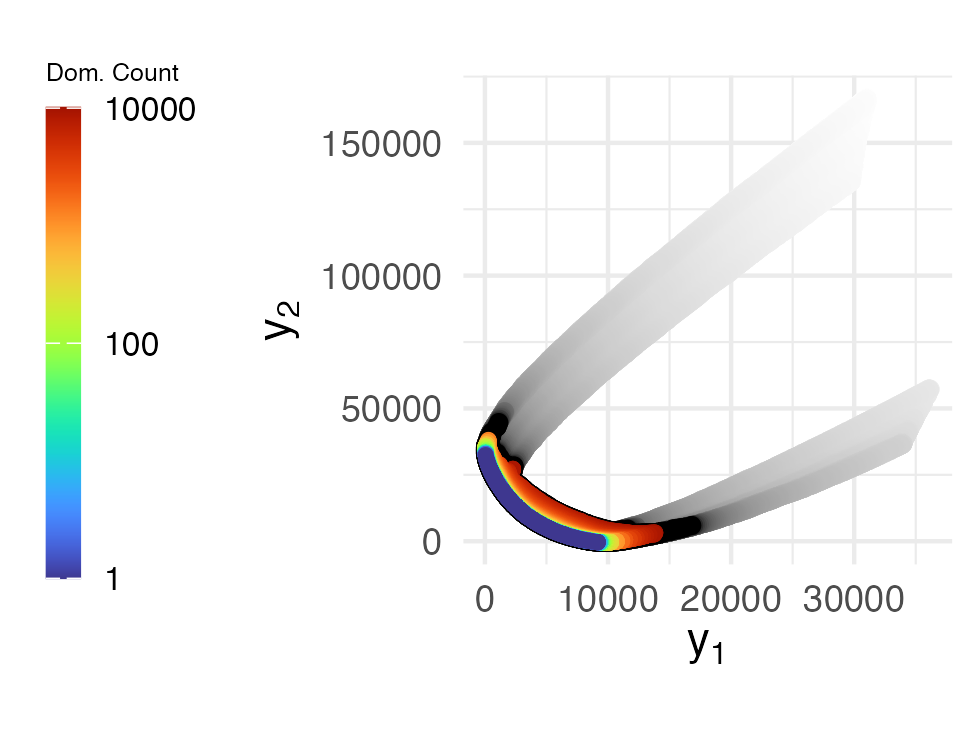}
        \caption{BONO14 instance: Multimodal stepped ellipsoids.}
    \end{subfigure}

    \caption{Multimodal functions with global free ellipsoid shapes. Many local optima are created around the Pareto set, which can form difficult plateaus in the stepped version.}
    \label{fig:mm-free-ellipsoids}
\end{figure}

\subsubsection{BONO13: Multimodal Free Ellipsoids}

The multimodal free ellipsoid functions are presented in \Cref{fig:mm-free-ellipsoids}.
Similar to BONO10-12, many locally efficient points are created in the vicinity of the Pareto set.

\subsubsection{BONO14: Multimodal Stepped Ellipsoids}

In contrast to BONO13, the stepped ellipsoids problems include discretization, which lets the number of locally efficient points increase even further.
Note that the sampling density of the visualization may influence how many locally efficient points can be detected, as a decreasing step size will eventually always find neighboring points rounded to the same level sets.

\subsection{Multimodal Problems without Global Structure} \label{sec:bono_multi_without}

\begin{table}[t]
    \centering
    \caption{Overview of the most important parameter variations within the multimodal BONO-Bench problems without global structure: Conditioning ($\kappa$), distance parameter ($p$), number of discretization steps ($N_h$), and number of constituent peak functions ($J$).
    }
    \label{tab:bono_multi_without}
    \begin{tabular}{ccccc}
        \toprule \bfseries Problem ID & $\boldsymbol{\kappa}$ & $\boldsymbol{p}$ & $\boldsymbol{N_h}$ & $\boldsymbol{J}$\\
        \midrule
        BONO15  & $1$ & $2$ & - & $10$ \\
        BONO16  & $1$ & $\LU(\frac 1 3, 3)$ & - & $100$ \\
        BONO17  & $1$ & $\LU(\frac 1 3, 3)$ & $\lfloor \LU(50,201) \rfloor$ & $100$ \\ \midrule
        BONO18  & $\LU(50,200)$ & $\LU(\frac 1 3, 3)$ & - & $10$ \\
        BONO19  & $\LU(50,200)$ & $\LU(\frac 1 3, 3)$ & - & $100$ \\
        BONO20  & $\LU(50,200)$ & $\LU(\frac 1 3, 3)$ & $\lfloor \LU(50,201) \rfloor$ & $100$ \\ \midrule
    \end{tabular}
\end{table}

\begin{figure}
    \centering
    \begin{subfigure}{\textwidth}
        \includegraphics[width=0.245\linewidth]{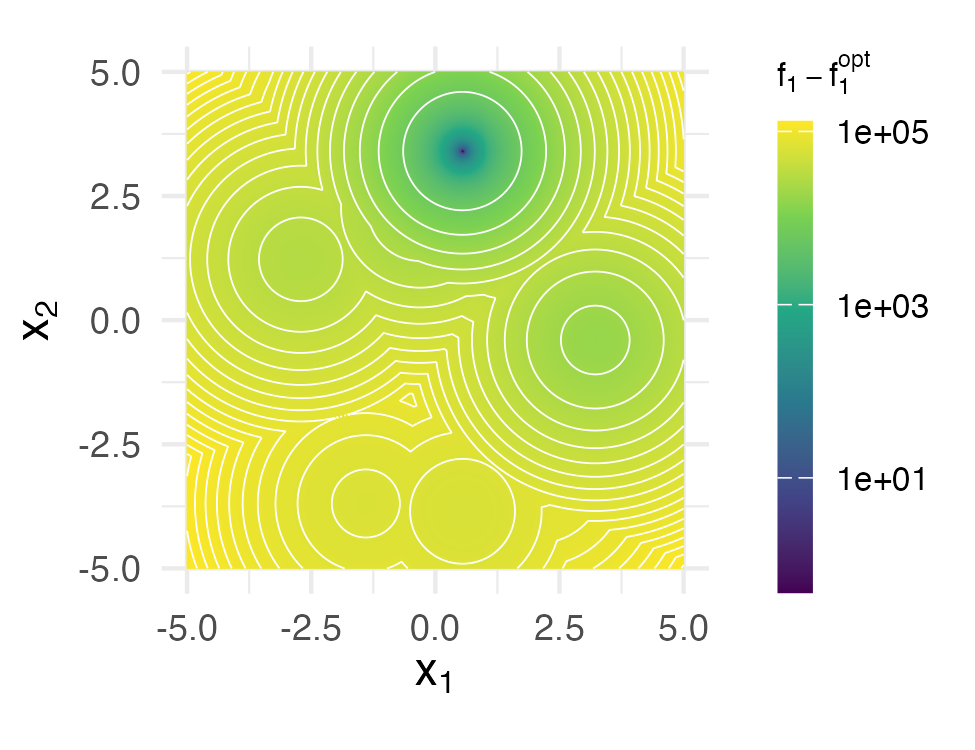}
        \includegraphics[width=0.245\linewidth]{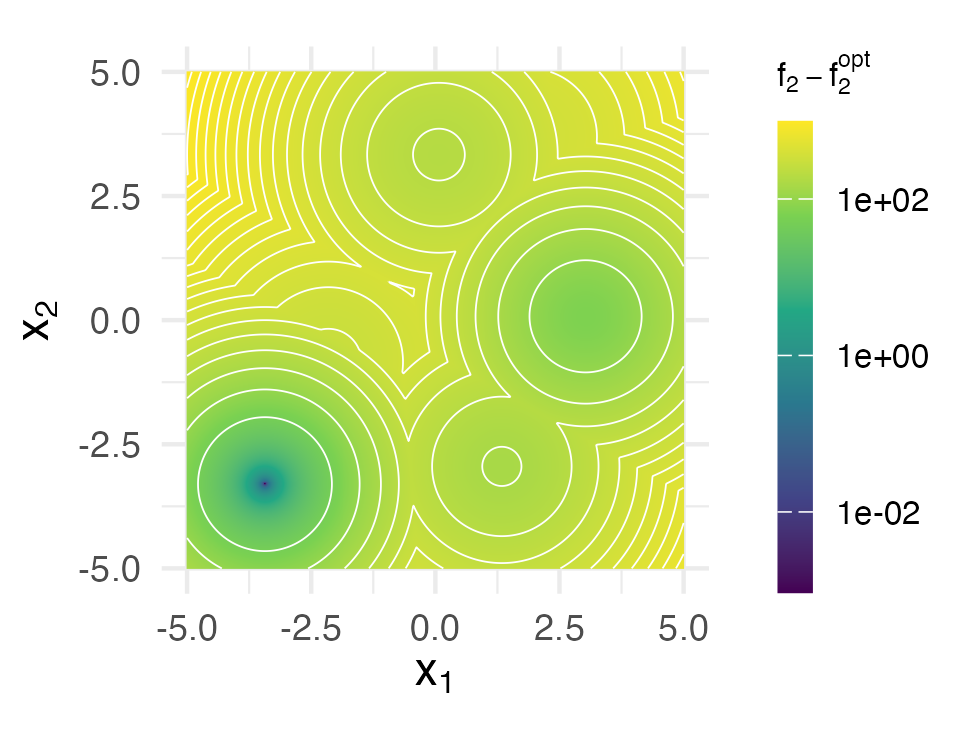}
        \includegraphics[width=0.245\linewidth]{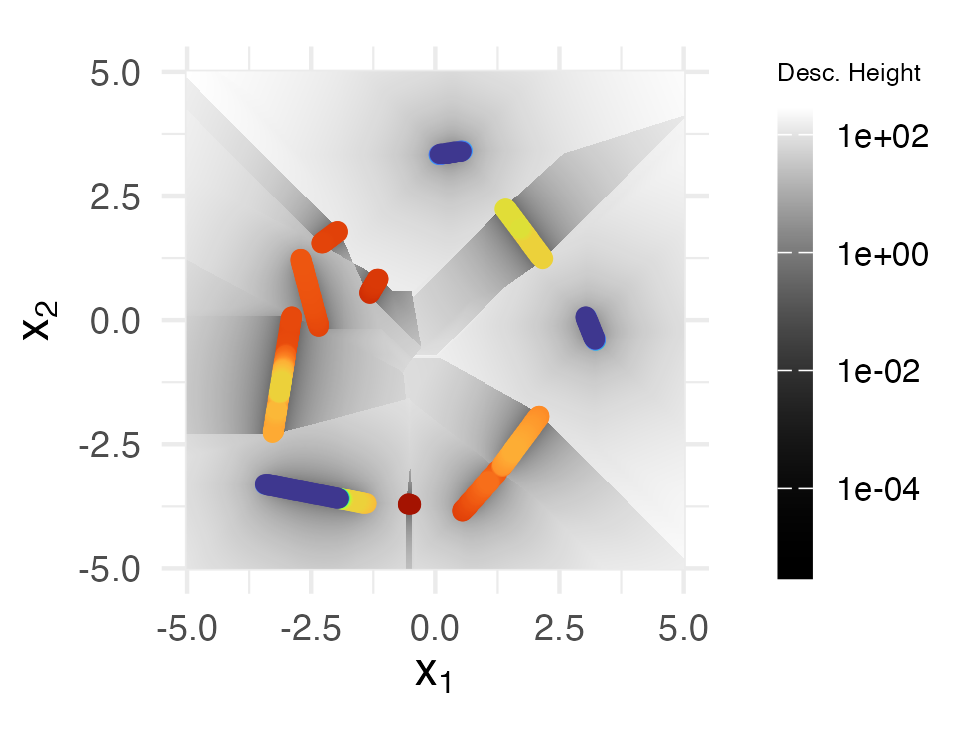}
        \includegraphics[width=0.245\linewidth]{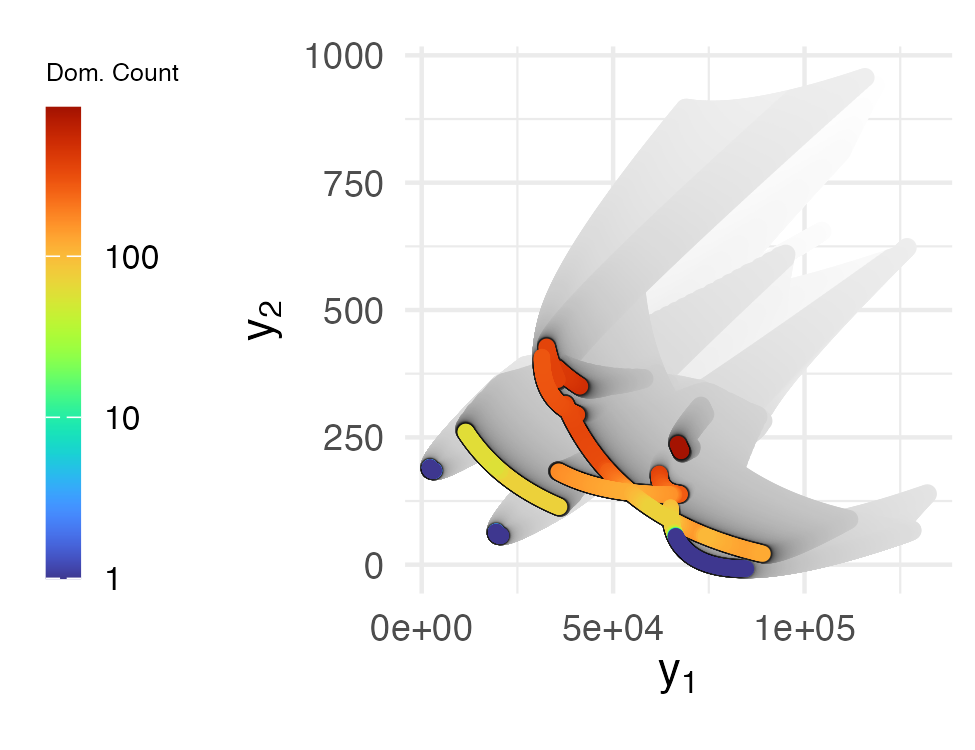}
        \caption{BONO15 instance: Few spheres.}
    \end{subfigure}
    \begin{subfigure}{\textwidth}
        \includegraphics[width=0.245\linewidth]{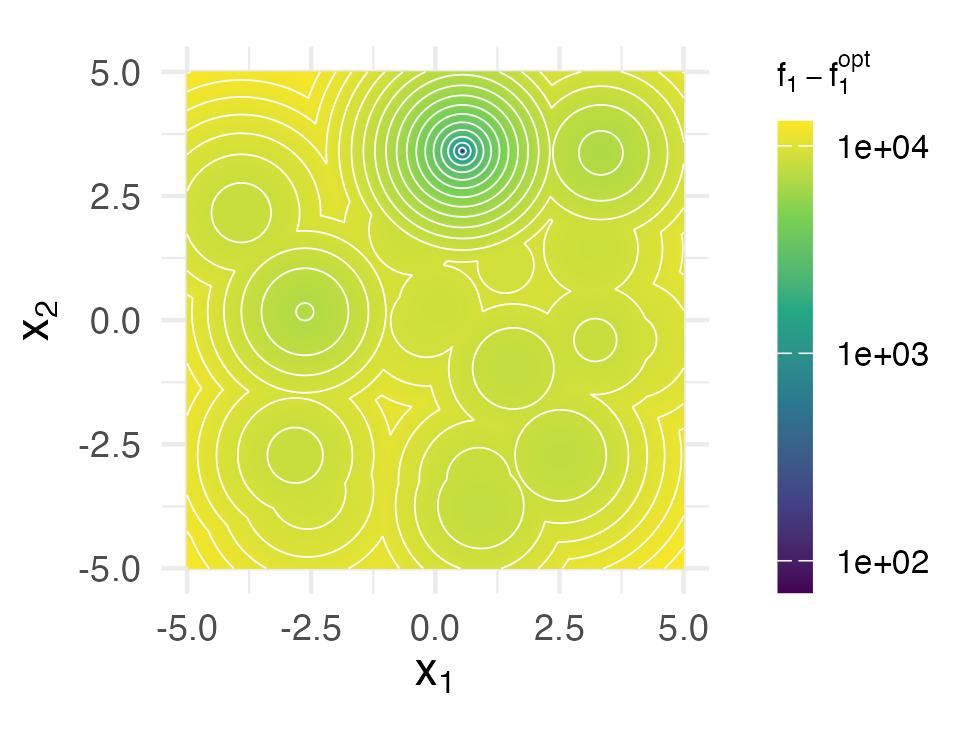}
        \includegraphics[width=0.245\linewidth]{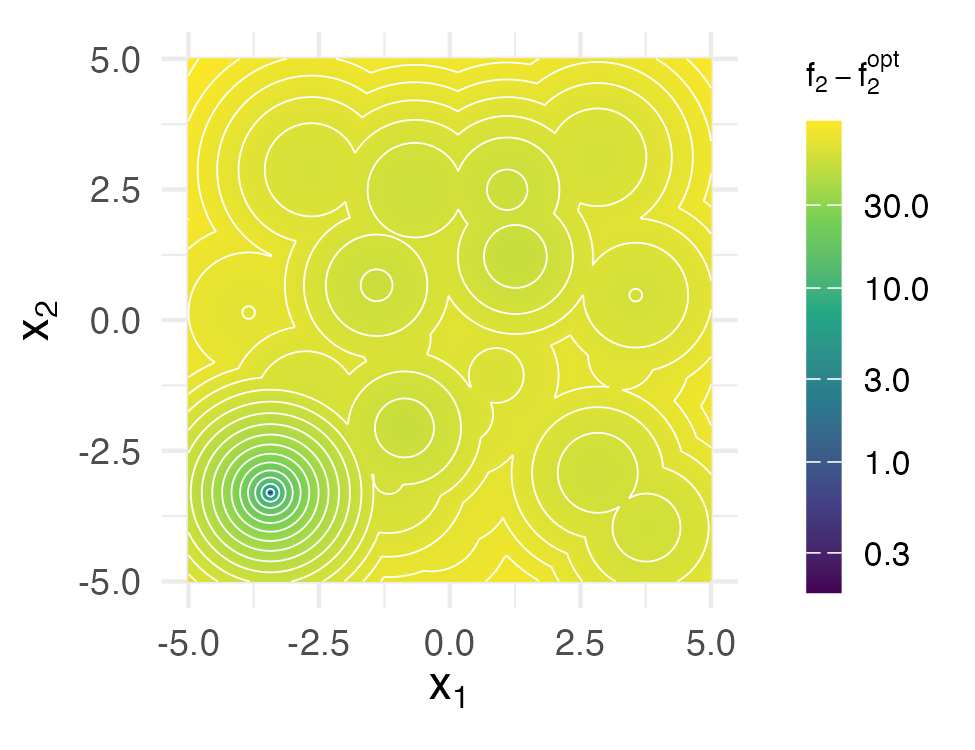}
        \includegraphics[width=0.245\linewidth]{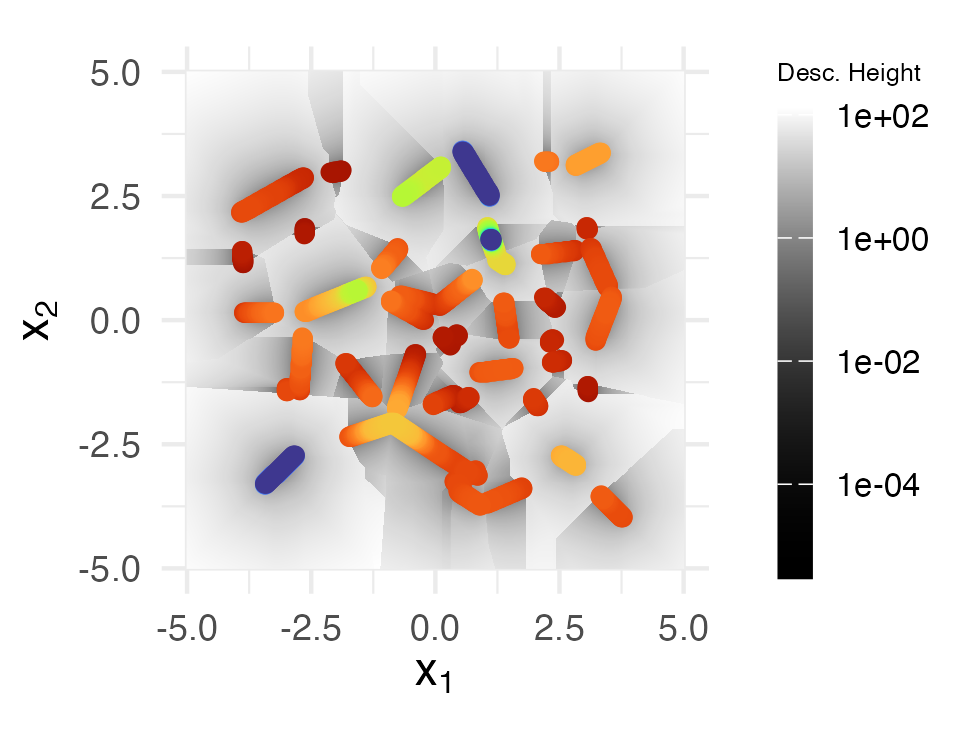}
        \includegraphics[width=0.245\linewidth]{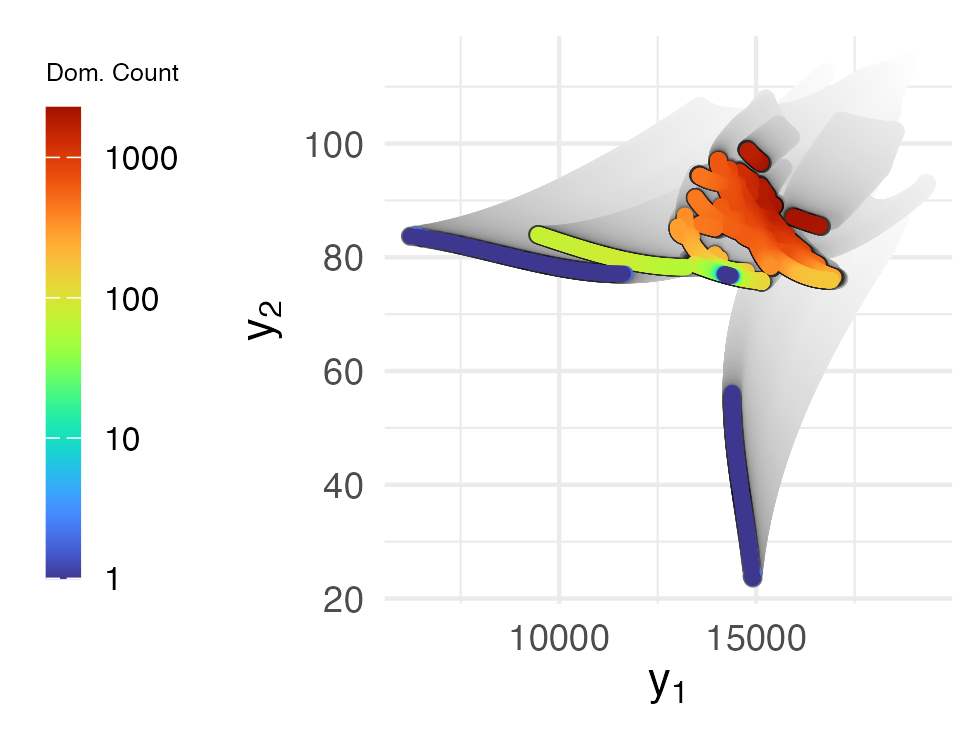}
        \caption{BONO16 instance: Many spheres.}
    \end{subfigure}
    \begin{subfigure}{\textwidth}
        \includegraphics[width=0.245\linewidth]{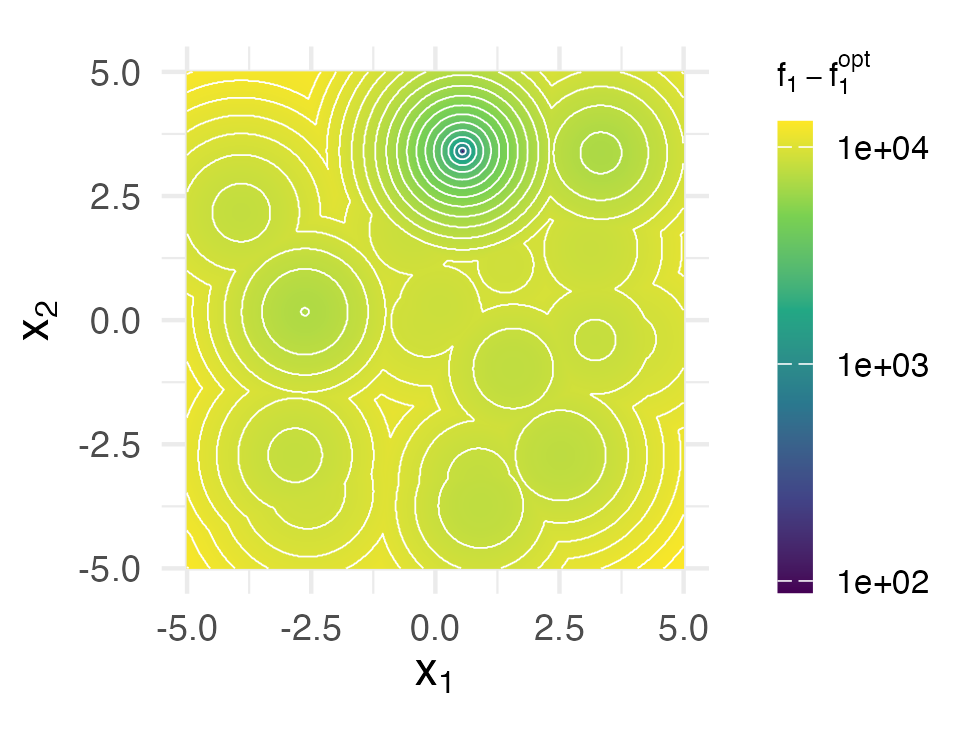}
        \includegraphics[width=0.245\linewidth]{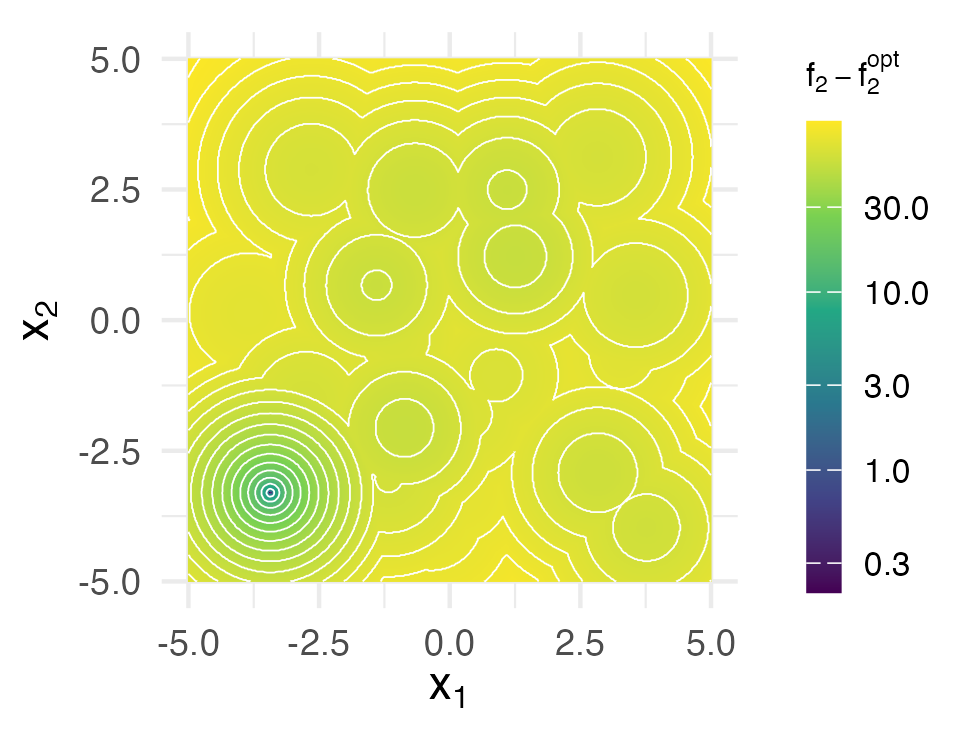}
        \includegraphics[width=0.245\linewidth]{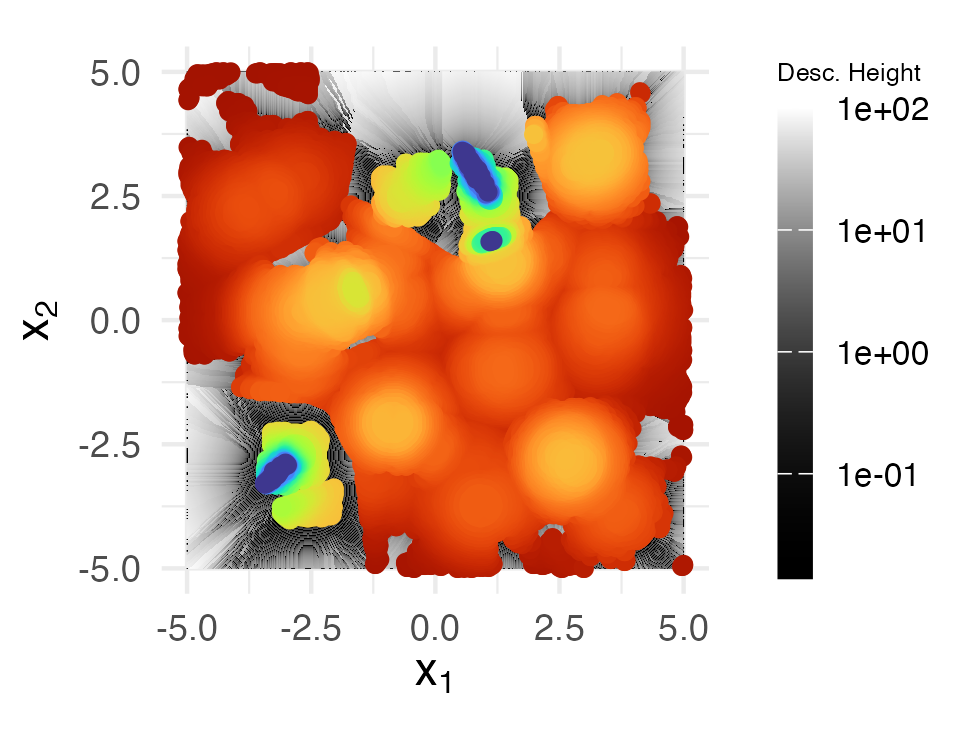}
        \includegraphics[width=0.245\linewidth]{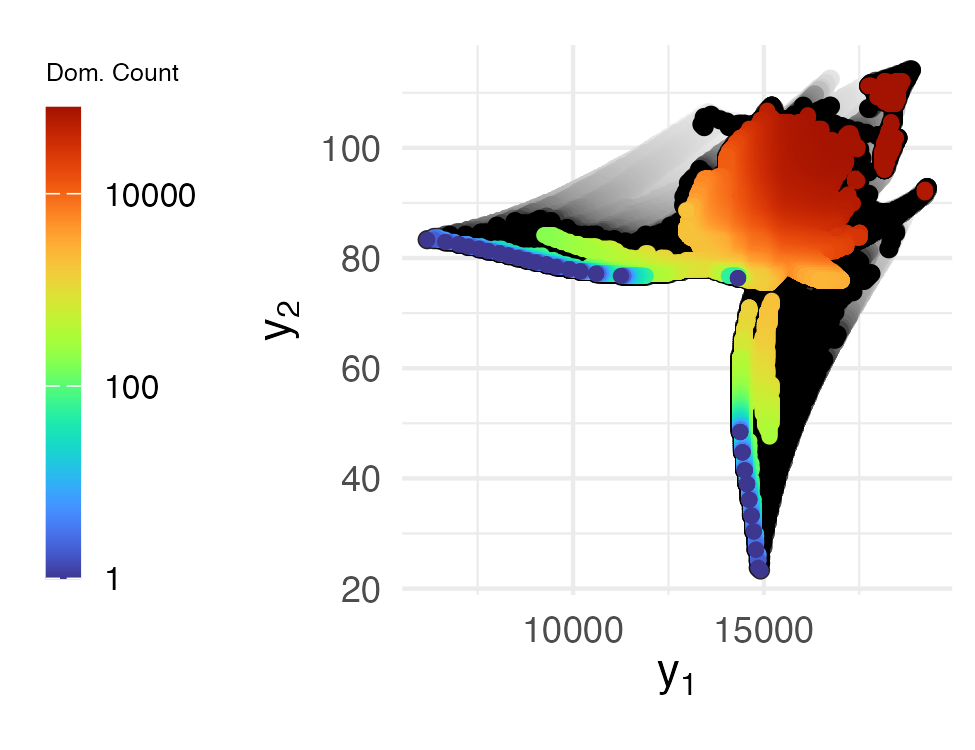}
        \caption{BONO17 instance: Stepped many spheres.}
    \end{subfigure}

    \caption{Sphere-based multimodal problems without global structure. The amount of multi-objective local optima increases quickly with the number of constituent sphere functions in the search landscape.}
    \label{fig:sphere-mm}
\end{figure}

\begin{figure}
    \centering
    \begin{subfigure}{\textwidth}
        \includegraphics[width=0.245\linewidth]{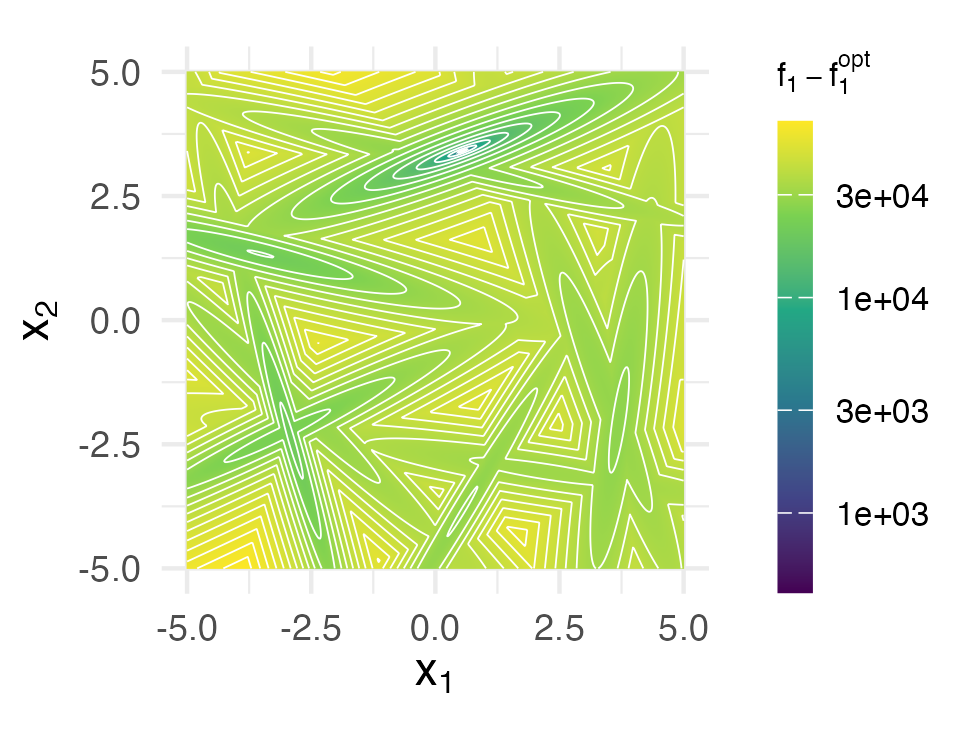}
        \includegraphics[width=0.245\linewidth]{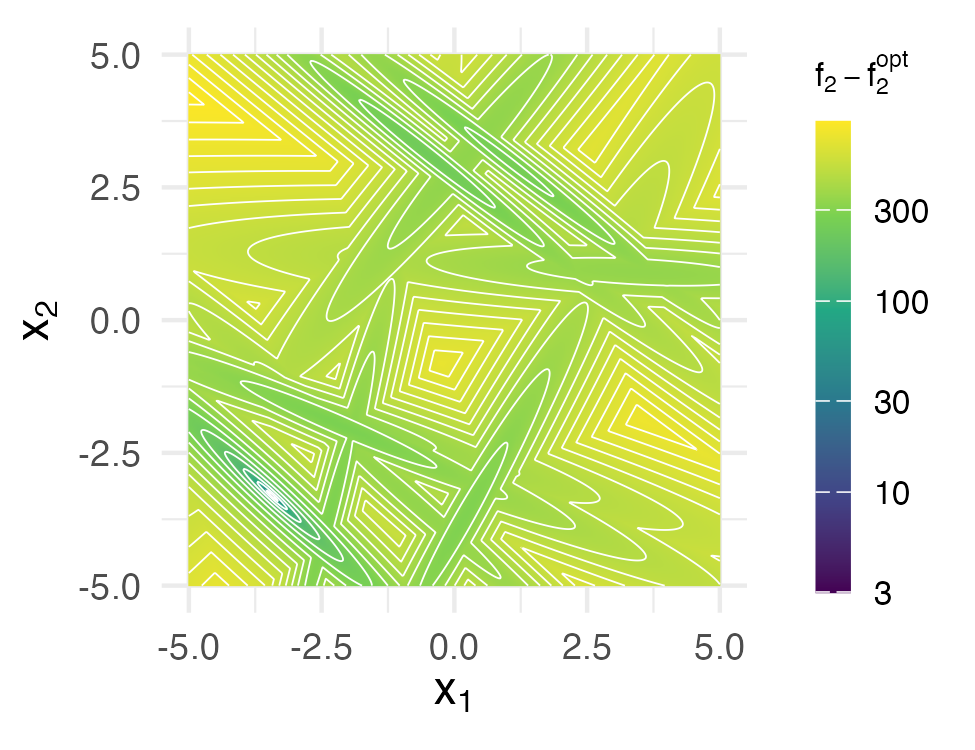}
        \includegraphics[width=0.245\linewidth]{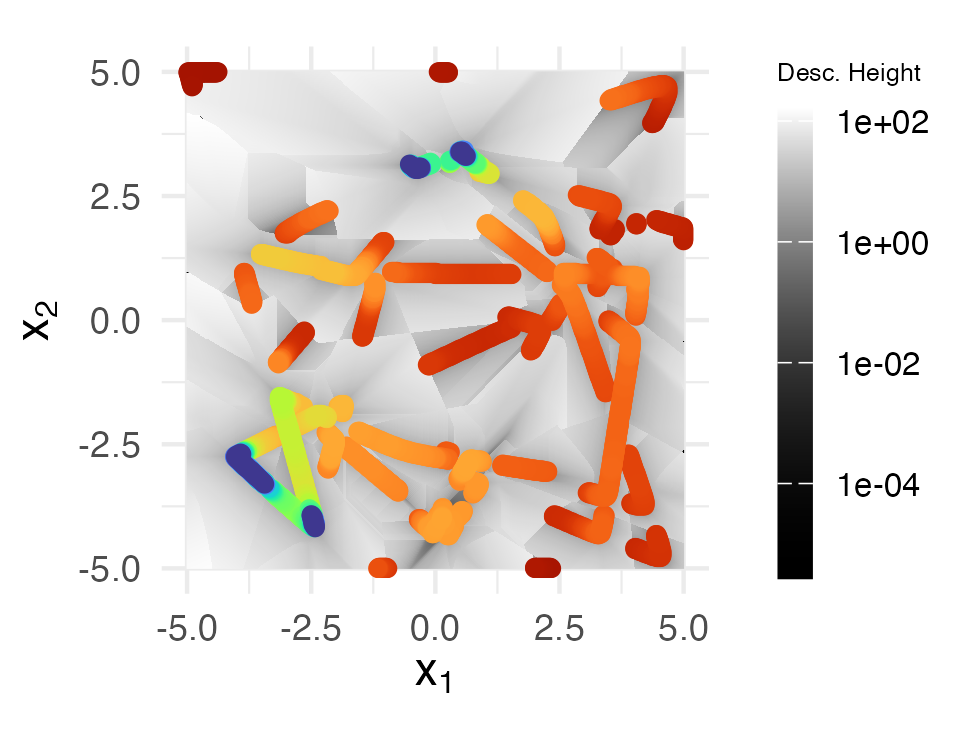}
        \includegraphics[width=0.245\linewidth]{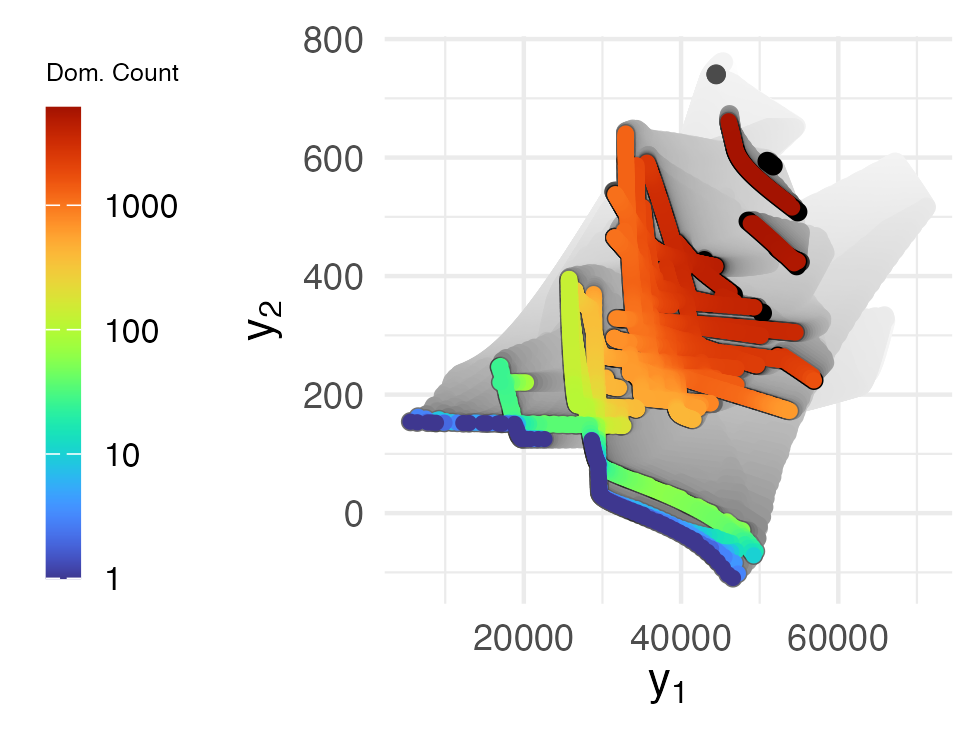}
        \caption{BONO18 instance: Few ellipsoids.}
    \end{subfigure}
    \begin{subfigure}{\textwidth}
        \includegraphics[width=0.245\linewidth]{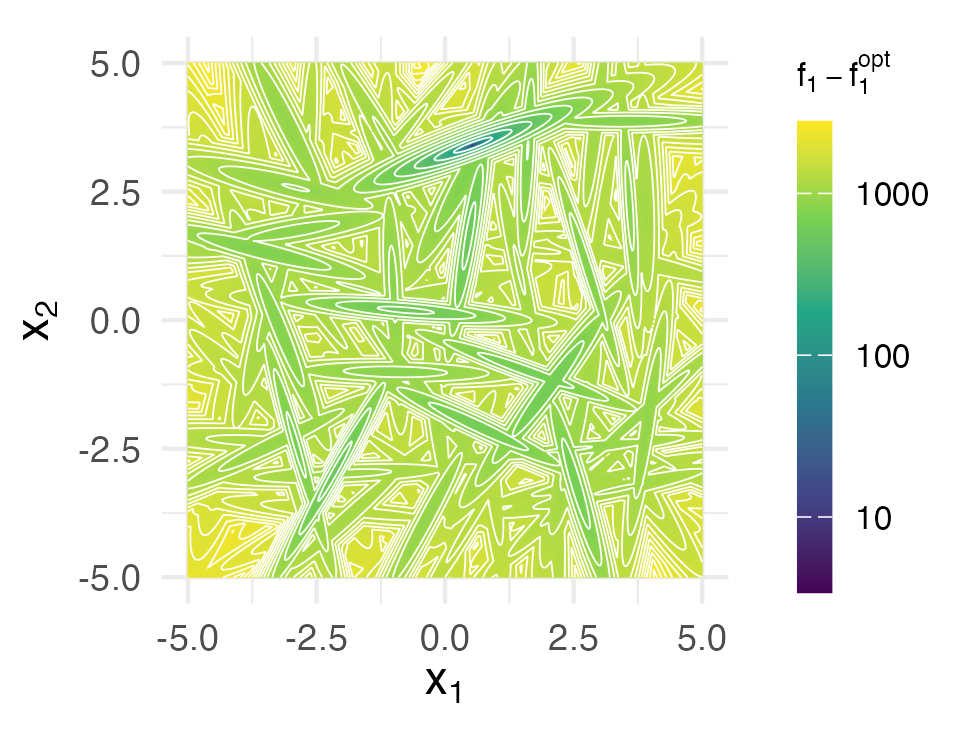}
        \includegraphics[width=0.245\linewidth]{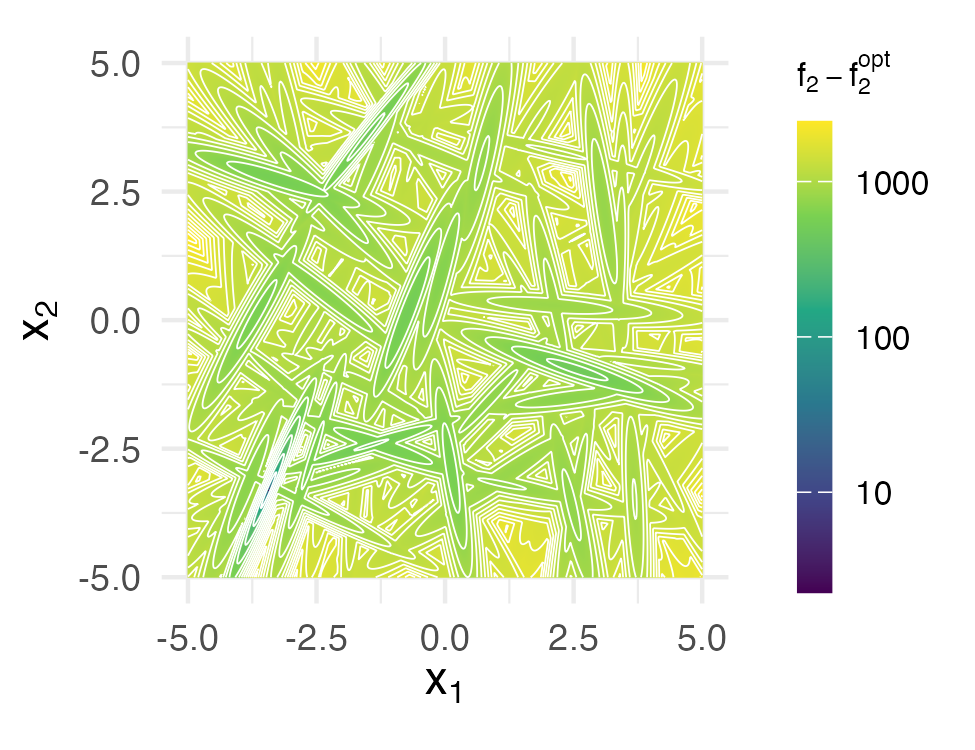}
        \includegraphics[width=0.245\linewidth]{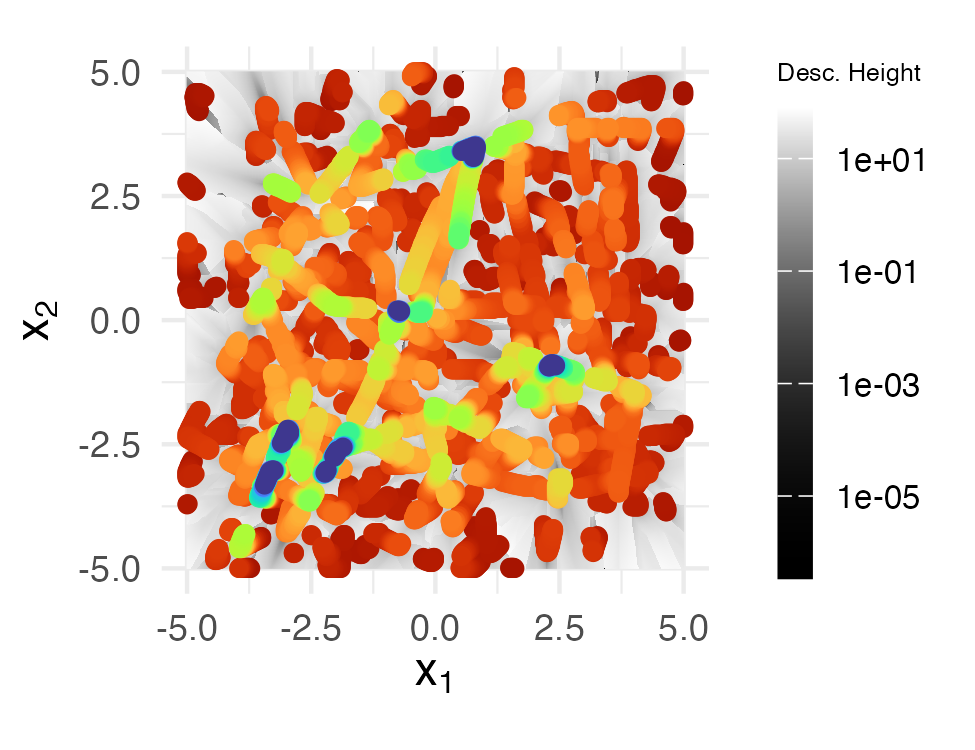}
        \includegraphics[width=0.245\linewidth]{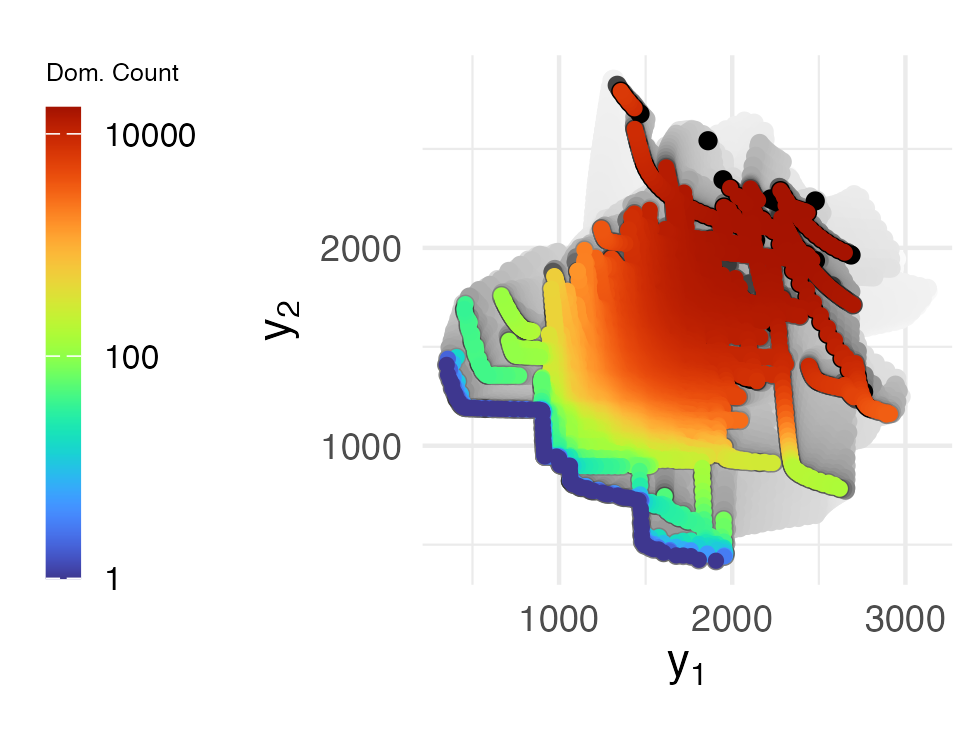}
        \caption{BONO19 instance: Many ellipsoids.}
    \end{subfigure}
    \begin{subfigure}{\textwidth}
        \includegraphics[width=0.245\linewidth]{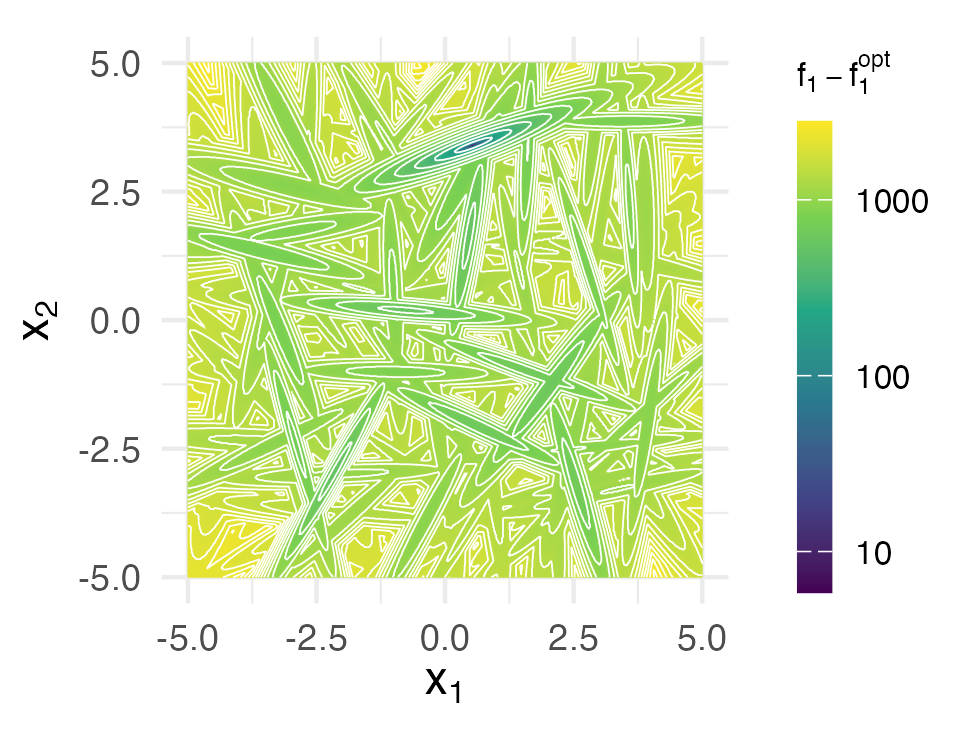}
        \includegraphics[width=0.245\linewidth]{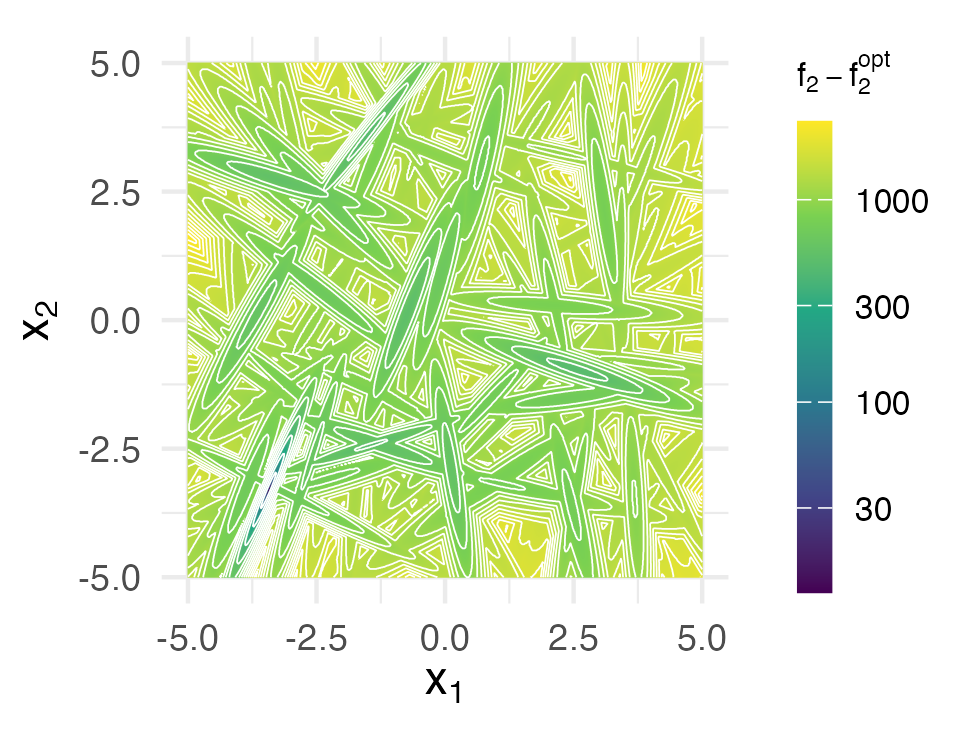}
        \includegraphics[width=0.245\linewidth]{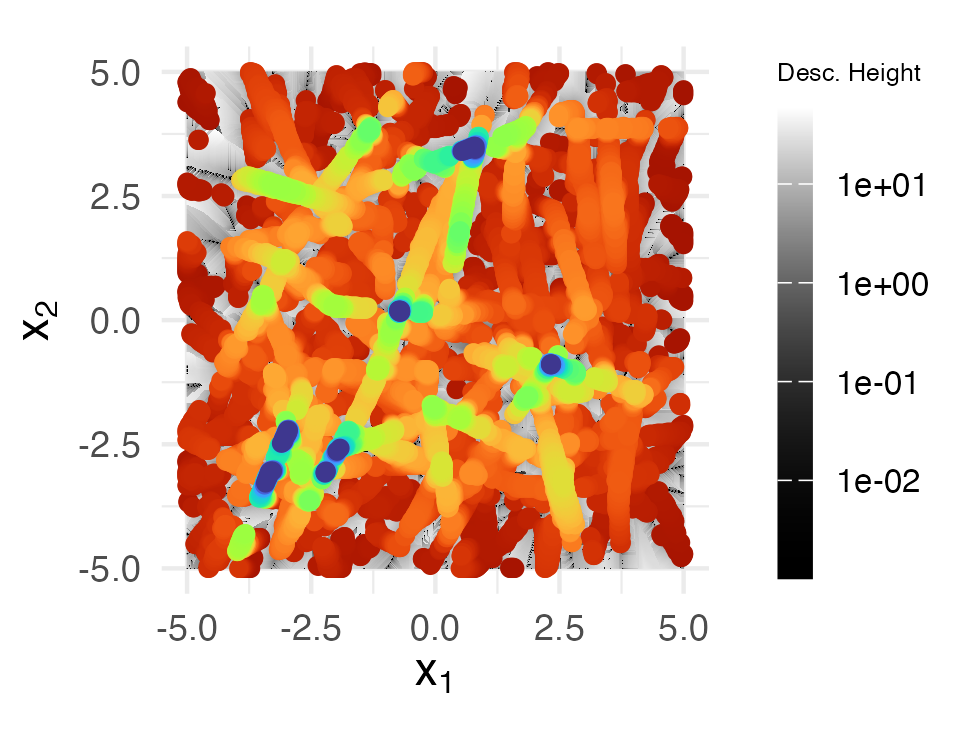}
        \includegraphics[width=0.245\linewidth]{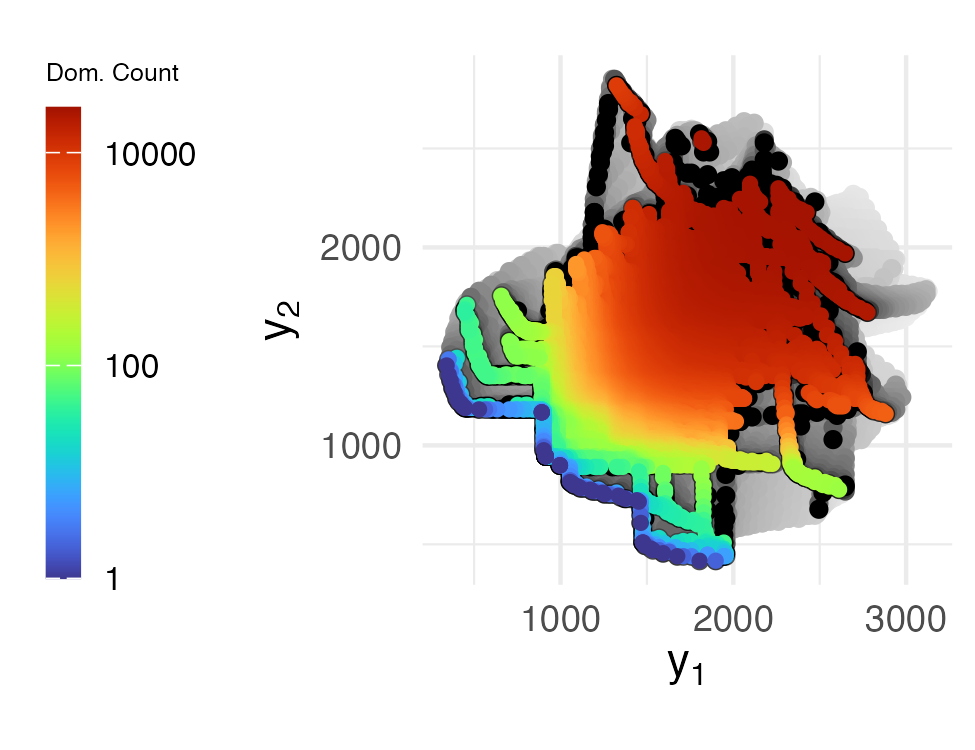}
        \caption{BONO20 instance: Stepped many ellipsoids.}
    \end{subfigure}

    \caption{Ellipsoid-based multimodal problems without global structure. The number of local optima as well as constituent pieces of the Pareto front increases with the number of local optima.}
    \label{fig:ellipsoid-mm}
\end{figure}

Finally, we create a set of multimodal problems \textit{without global structure} using spherical and ellipsoidal base functions, respectively.
While the multi-sphere problems (BONO15-17) are shown in \Cref{fig:sphere-mm}, the ellipsoidal problems without global structure (BONO18-20) are depicted in \Cref{fig:ellipsoid-mm}.
An overview of the most important parameter variations is given in \Cref{tab:bono_multi_without}.
This class of problems gives rise to complex structures, with local and global optima, i.e., (local) Pareto sets, being rather isolated and widely spread across the search space.
Consequently, PLOT visualizations reveal curved sets in the search space and multiple local Pareto fronts, which are located well behind the global Pareto front, in the objective space.

\subsubsection{BONO15: Few Spheres}

In this problem class, $10$ random sphere functions ($\kappa=1, p=2$) are sampled per problem.

\subsubsection{BONO16: Many Spheres}

BONO16 follows the construction of the BONO15 problem, though $100$ sphere functions are used per objective, and the global distance parameter is sampled from $p \sim \LU(\frac 1 3, 3)$.

\subsubsection{BONO17: Stepped Many Spheres}

BONO17 problems are identical to BONO16 problems, though they are discretized into $N_h\sim\lfloor\LU(50,201)\rfloor$ steps per objective within the area of interest.

\medskip
In contrast to the sphere-based problems, the ellipsoidal problems use $\kappa \sim \LU(50,200)$ for the sampling of the Hessians.
Also, all ellipsoidal problems use $p \sim \LU(\frac 1 3, 3)$.

\subsubsection{BONO18: Few Ellipsoids}
BONO18 uses $10$ ellipsoids per objective.

\subsubsection{BONO19: Many Ellipsoids}

BONO19 uses $100$ ellipsoids rather than the $10$ used in BONO18.

\subsubsection{BONO20: Stepped Many Ellipsoids}

BONO20 problems are equal to BONO19 problems, but again are discretized into $N_h\sim\lfloor\LU(50,201)\rfloor$ steps per objective.

\section{Experimental Study} \label{sec:study}

To demonstrate the properties and usability of the BONO-Bench test suite, we perform an experimental study involving a diverse set of optimizers.
This study does not claim to represent the complete history and state of the art in bi-objective optimization, but rather tries to demonstrate the applicability of our approach.
For this reason, we do not perform extensive algorithm configuration, but rather focus on applying well-known optimizers out of the box with (mostly) default settings.

\subsection{Experimental Setup} \label{sec:setup}

Below, we outline the experimental setup, detailing the problem set, the optimizers used, the evaluation procedure, and considerations regarding implementation and reproducibility.

\subsubsection{Problem Set}

Our problem set is based on the $20$ problem classes introduced in \Cref{sec:bono}.
To study the scaling effect of the number of variables, we consider decision space dimensionalities $d \in \{2,3,5,10,20\}$.
Moreover, by varying the random seed, we create $15$ unique instances for each combination of BONO problem and decision space dimensionality, resulting in a total of $20 \times 5 \times 15 = 1\,500$ problem instances.

\subsubsection{Optimizers}

For this study, we select several well-known optimizers that operate without problem-specific configurations and can thus be executed in a true black-box manner.
Specifically, we use the SPEA2 \cite{zitzler2001spea2}, NSGA2 \cite{deb2002fast} and SMS-EMOA \cite{beume2007sms} implementations from the \texttt{pymoo} Python package \cite{blank2020pymoo}.
All \texttt{pymoo} optimizers are run with their default configurations, including a population size of $100$, binary tournament selection, simulated binary crossover, and polynomial mutation.
In addition, we use the default GDE3 implementation from the \texttt{pymoode} package \cite{leite2023multi}, which employs a DE/rand/1/bin strategy, as a stand-in for differential evolution algorithms as well as the implementation of the MO-CMA-ES \cite{igel2007covariance} from the \texttt{deap} package \cite{DEAP_JMLR2012}.
For NSGA2 and SMS-EMOA, we also perform an ablation regarding their population size, for which we additionally use population sizes of $50$ and $200$ individuals.
All algorithms are initialized using uniform random sampling within the decision space $\mathcal X = [-5,5]^d$.
As a baseline, we include random search, which continues sampling uniformly at random from the decision space until the budget is exhausted.

\subsubsection{Evaluation Procedure}

All optimizers are executed with a budget of $10^5 d$ function evaluations, where $d$ is the search space dimensionality.
To record performance, we first approximate the Pareto front using \Cref{alg:approximation} for the R2 and HV indicators, using a target precision of $\delta_{R2}=10^{-6}$ and $\delta_{HV}=10^{-5}$, respectively, yielding $\hat I_{R2}^*$ and $\hat I_{HV}^*$.
Empirically, $\delta_{R2}=10^{-6}$ and $\delta_{HV}=10^{-5}$ lead to roughly similarly well-approximated Pareto fronts.
Before the indicator computation, the area between the ideal and nadir points is normalized to the unit square $[0,1]^2$.
We use the resulting ideal $(0,0)$ and nadir $(1,1)$ points as reference points for the exact R2 and hypervolume indicator, respectively.

Then, we distribute 101 indicator targets log-uniformly in $[10^{-5},10^0]$ (R2) and $[10^{-4},10^0]$ (HV) w.r.t. the approximated indicator values $\hat I_{R2}^*$ and $\hat I_{HV}^*$, which gives a factor of 10 as buffer for the approximation error.
Whenever an indicator target value is surpassed by the unbounded archive of non-dominated points during a run, the number of function evaluations performed up to this point is recorded.
We found that these targets result in a roughly comparable level of difficulty for both indicators, though the resulting hypervolume targets tend to be solved a little faster (cf.~\Cref{sec:results}).
The performance of each solver is then summarized using runtime profiles \cite{hansen2022anytime}, which display the number of targets solved over time, i.e., consumed budget.
We also compute a virtual best solver (VBS) per indicator that combines the best runtimes across all algorithms for each problem instance.
In addition to runtime profiles depicting the performance across individual problems or the whole set of generators, we can also analyze performance differentiated by generator type (unimodal, multimodal with global structure, multimodal without global structure), shape of the global Pareto front (concave, linear, convex), objective space discretization, and other emergent properties (e.g., approximate number of peaks contributing to the Pareto front).

\subsubsection{Implementation and Reproducibility}

The BONO-Bench test suite is implemented in Python 3.11, utilizing the packages \texttt{numpy}, \texttt{pandas}, \texttt{matplotlib}, \texttt{moocore}, and \texttt{sortedcontainers}.
The code to reproduce all experiments, including the raw experimental results, is published at Zenodo \cite{zenodo2025bonobench}.
We make the accompanying Python package \texttt{bonobench} available as a standalone tool at \url{https://github.com/schaepermeier/bonobench}, providing easy access to the optimization problems and evaluation procedures.

\subsection{Verifying Benchmark Properties} \label{sec:properties}

\begin{figure}
    \centering
    \includegraphics[width=0.495\linewidth]{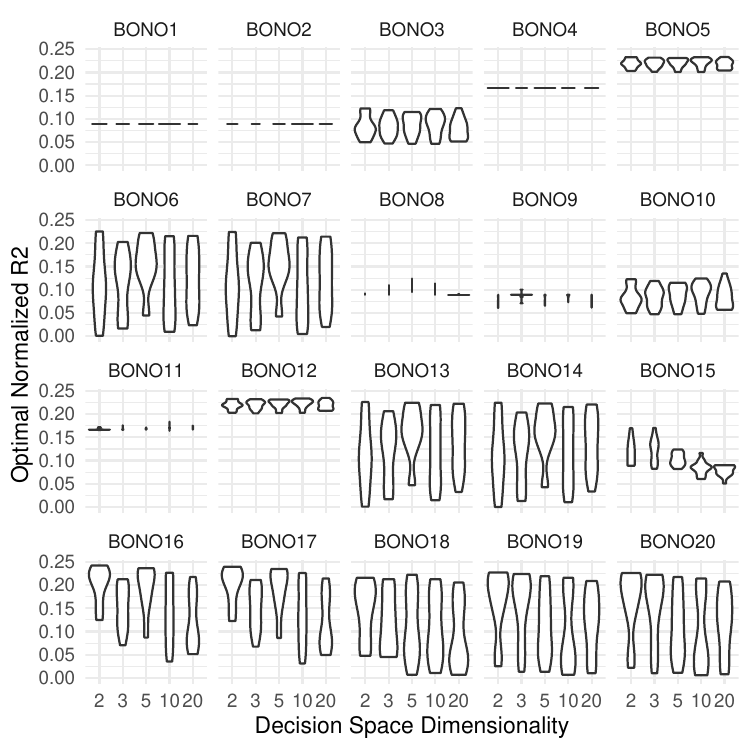}
    \hfill
    \includegraphics[width=0.495\linewidth]{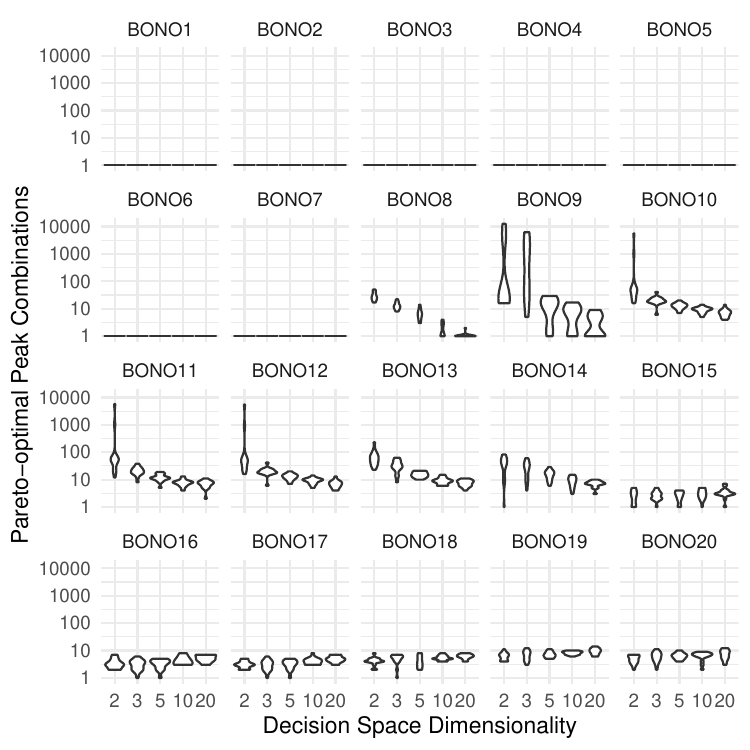}
    \caption{Distribution of BONO-Bench properties by generator and dimensionality.
    \emph{Left:} The optimal normalized R2 value verifies the optimal value in cases of fixed $p$ (BONO1,2,4).
    Other problems also show the expected distributions in the unimodal problems and multimodal problems with global structure.
    For multimodal problems without global structure, a slight drift to lower value ranges with increasing dimensionality can be observed.
    \emph{Right:}
    We can identify a drift to lower number of relevant peak combinations in BONO8-14, though no such drift is visible in the problems without global structure.}
    \label{fig:front-analysis}
\end{figure}

In \Cref{fig:front-analysis}, we visualize the distribution of the approximated normalized R2 value and the number of peak combinations that contribute to the Pareto front, which are derived during the indicator approximation.
It can be seen that there are only few drifts in the indicator value with increasing dimensionality, mostly relating to the last group of multimodal functions without global structure (BONO15-20).
Regarding the Pareto-optimal peak combinations, there is a drift to reduced multimodality of the globally optimal solutions within the middle function group with global structure (BONO8-14).
These issues are likely due to the curse of dimensionality, which affects distance calculations in high-dimensional spaces such that it is unlikely that any two points are close to each other.

The plots also highlight another issue: Some of the BONO9-12 problems have an excessive number of peak combinations contributing to the Pareto front for fewer decision variables ($d \in \{2,3\}$).
For some instances of the BONO9 problem ($d=2$), this led to convergence issues during the indicator approximation, where we early-stopped the approximation after $10^7$ iterations.
Further, we want to note that for in total $18$ (out of $1\,500$ of the) problem instances, we noticed that points outside the designated search area of $[-5,5]^d$ are evaluated in the Pareto front approximation.
These only occurred in multimodal problems (where we do not perform an ad hoc check at problem creation), and focused on higher dimensionalities.

We do not believe that any of these reported issues had a severe impact on the performance evaluation, but we want to transparently report them anyway.

\subsection{Experimental Results} \label{sec:results}

Here, we discuss the experimental results of the algorithm runs on the BONO-Bench suite.

\subsubsection{Runtime Profiles for All Problems}

\begin{figure}
    \centering
    \includegraphics[width=0.495\linewidth]{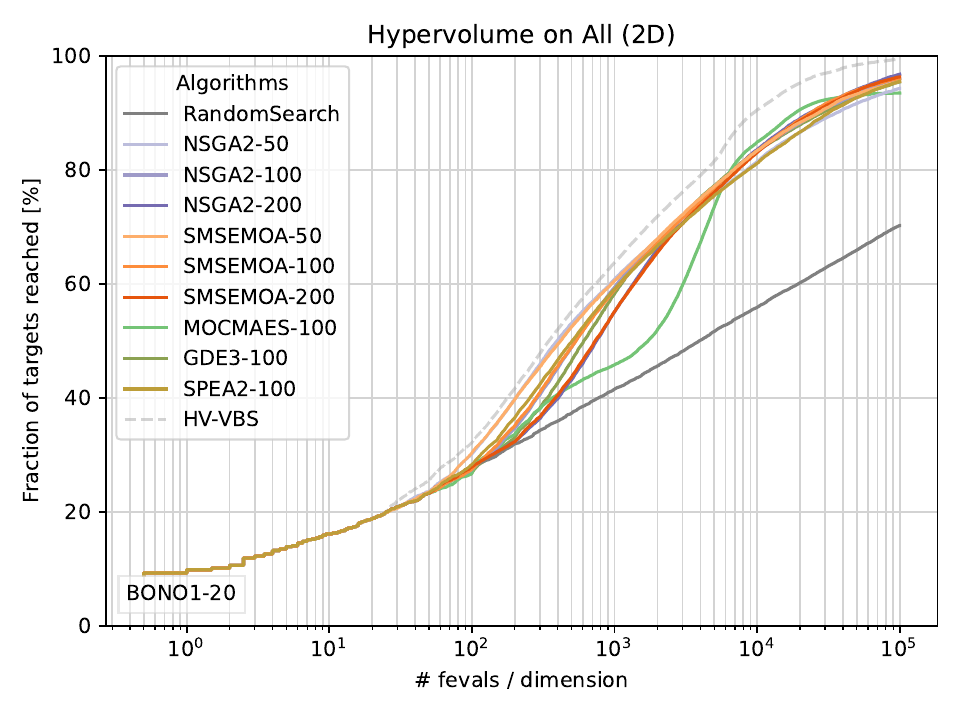}
    \includegraphics[width=0.495\linewidth]{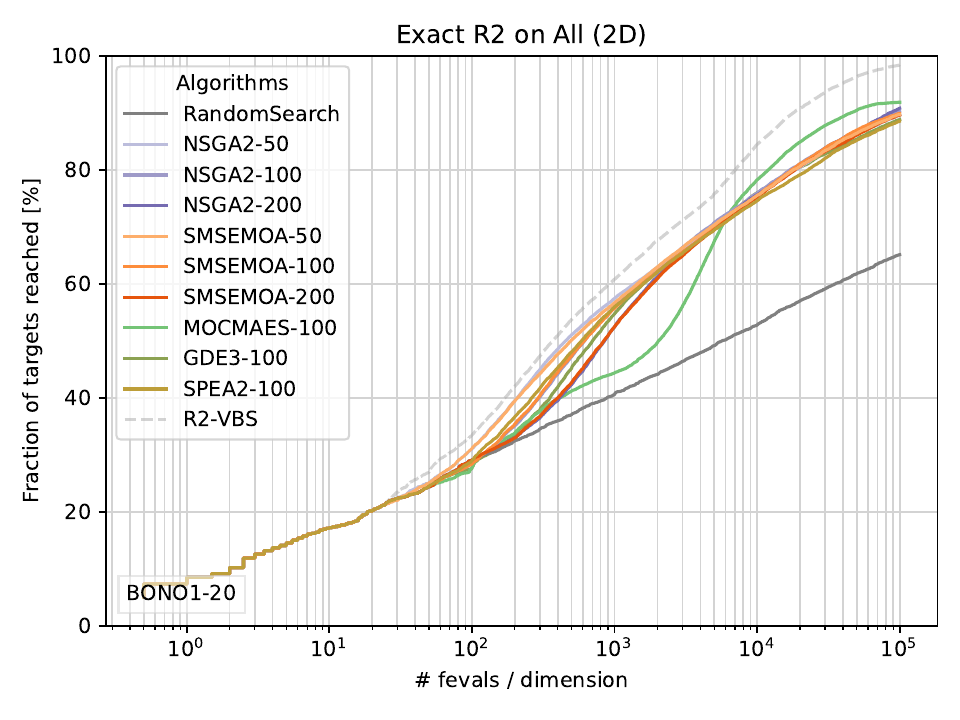}
    \includegraphics[width=0.495\linewidth]{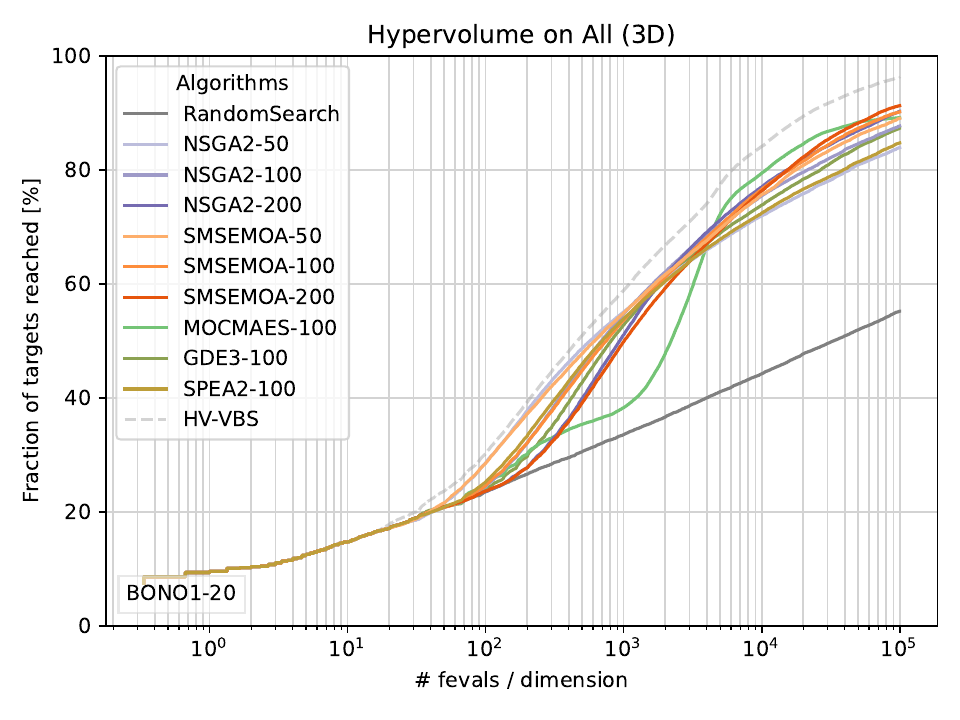}
    \includegraphics[width=0.495\linewidth]{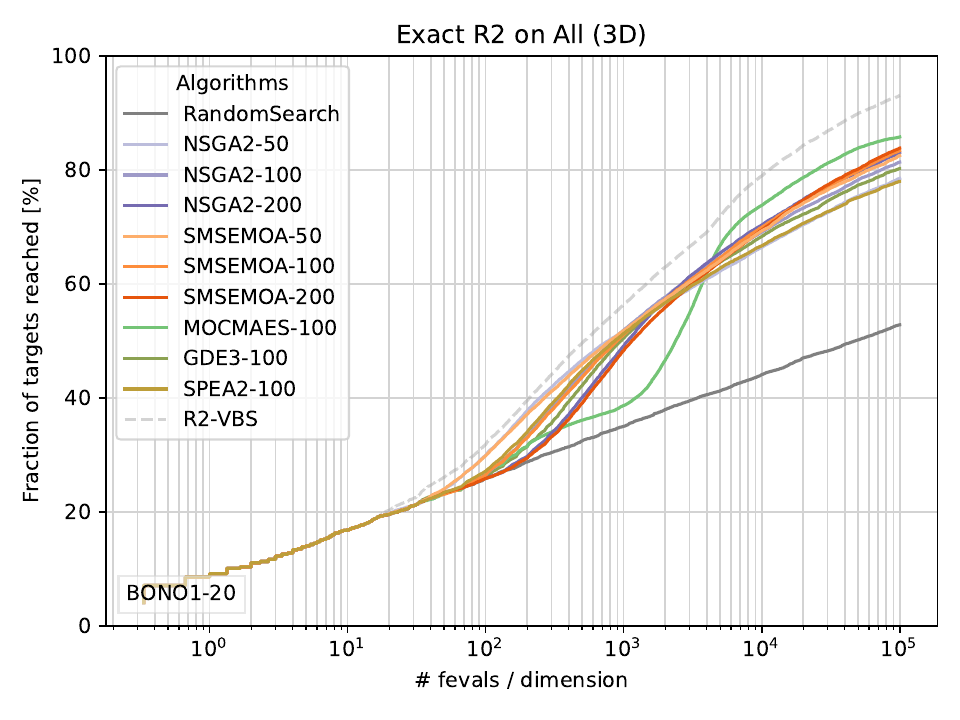}
    \includegraphics[width=0.495\linewidth]{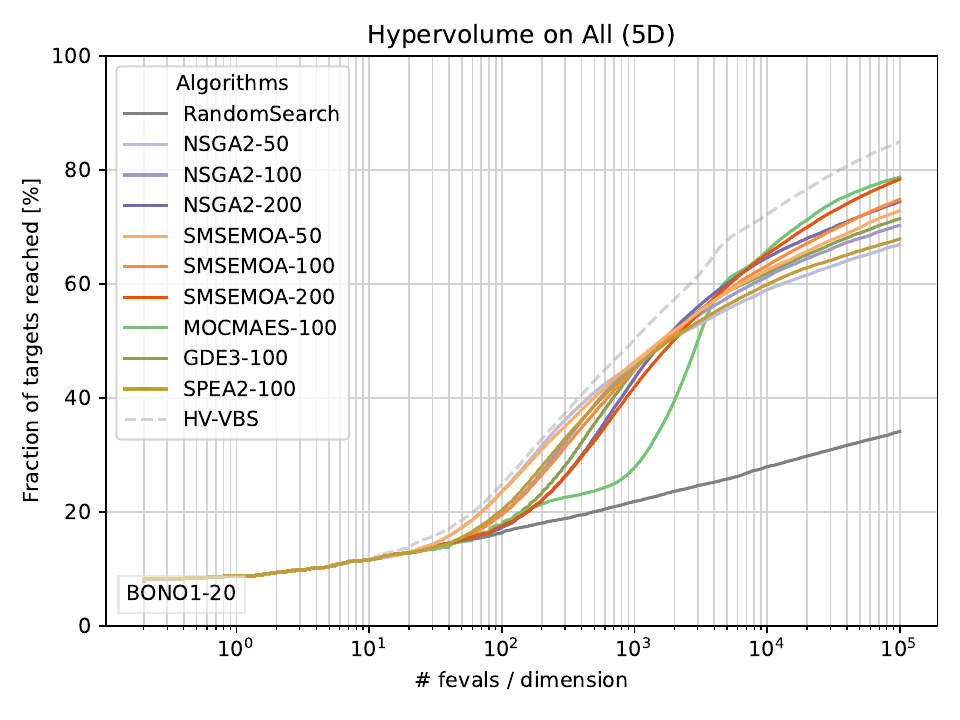}
    \includegraphics[width=0.495\linewidth]{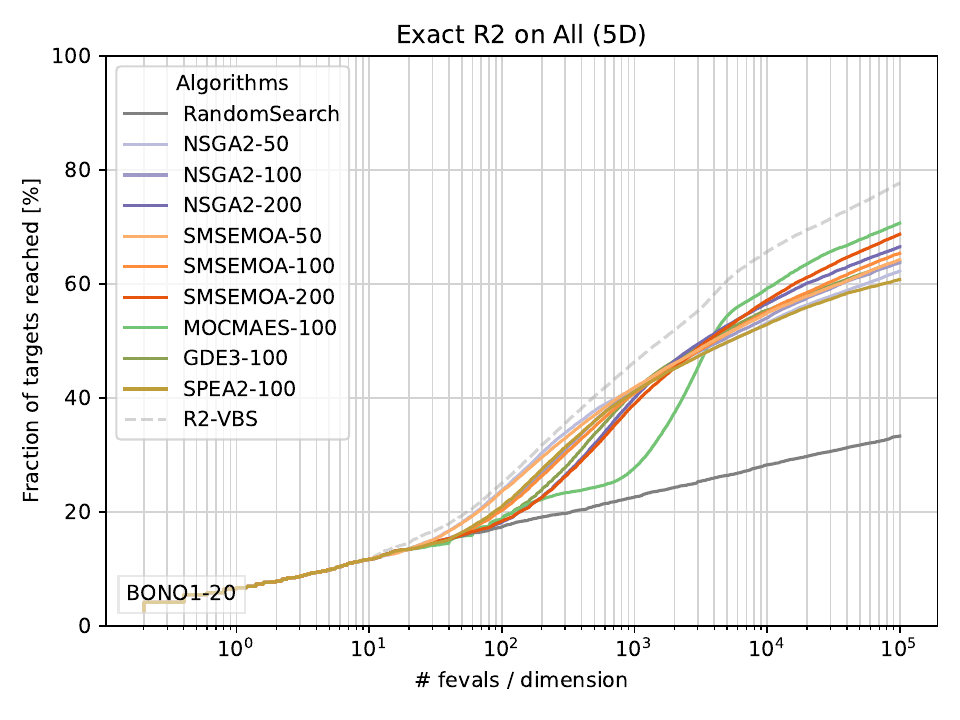}
    \caption{Runtime profiles for all problems with $d \in \{2,3,5\}$. Qualitatively speaking, hypervolume and R2 targets are similarly difficult to solve, though performance gaps between algorithms differ slightly. }
    \label{fig:rtp-all-235D}
\end{figure}

\begin{figure}
    \centering
    \includegraphics[width=0.495\linewidth]{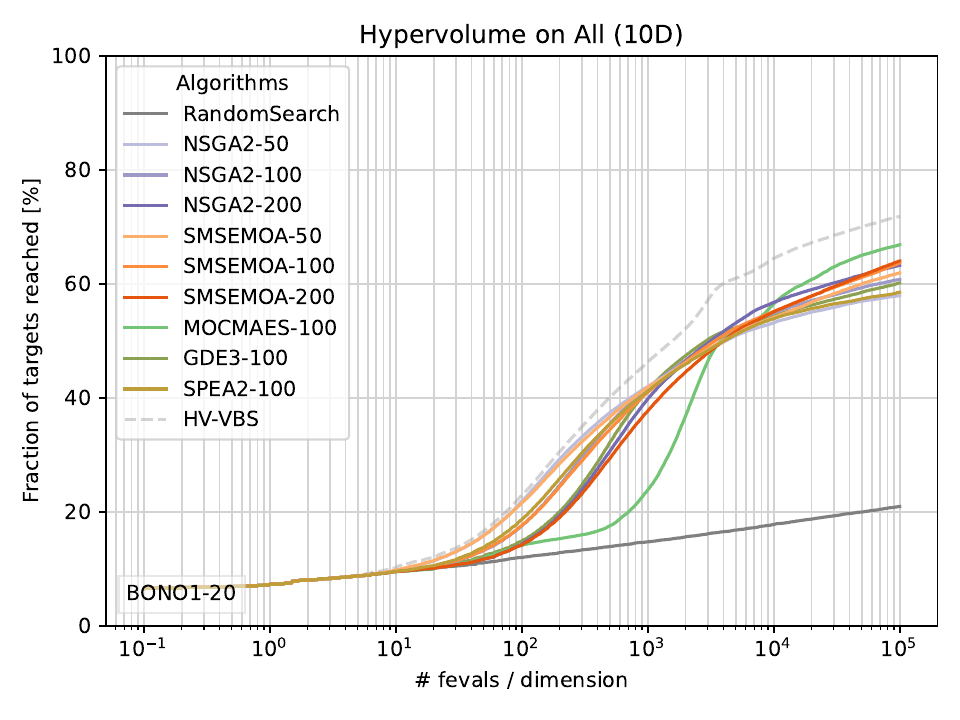}
    \includegraphics[width=0.495\linewidth]{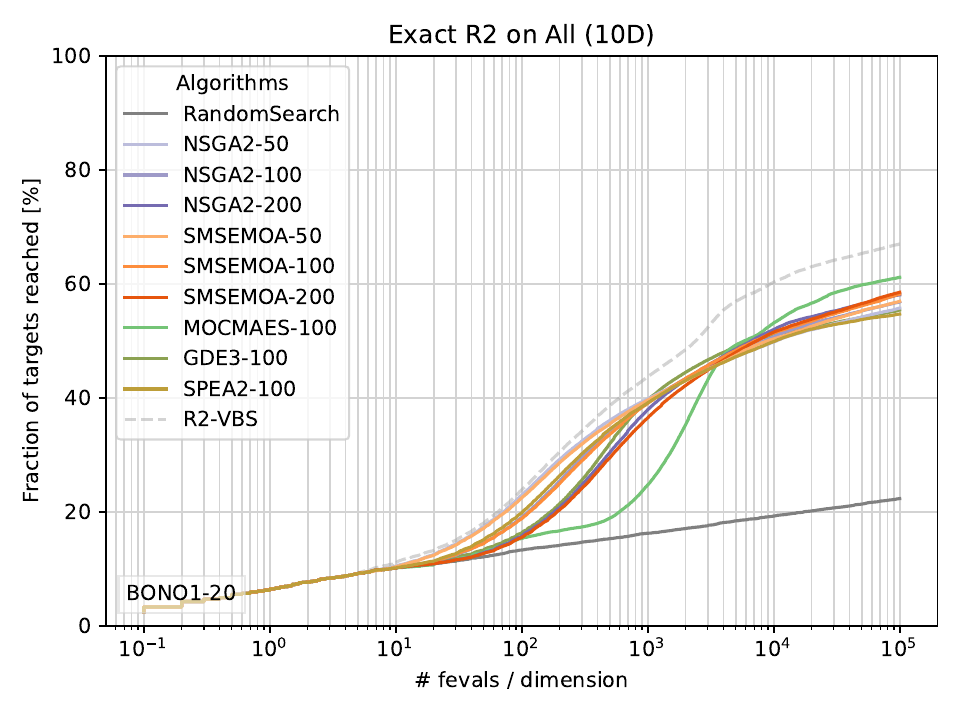}
    \includegraphics[width=0.495\linewidth]{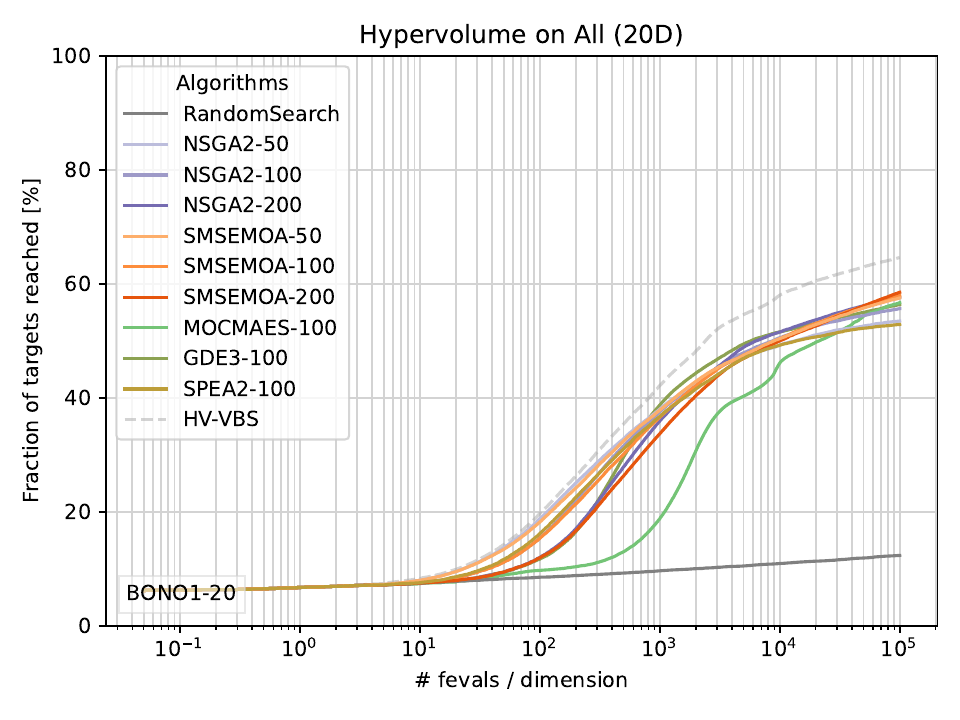}
    \includegraphics[width=0.495\linewidth]{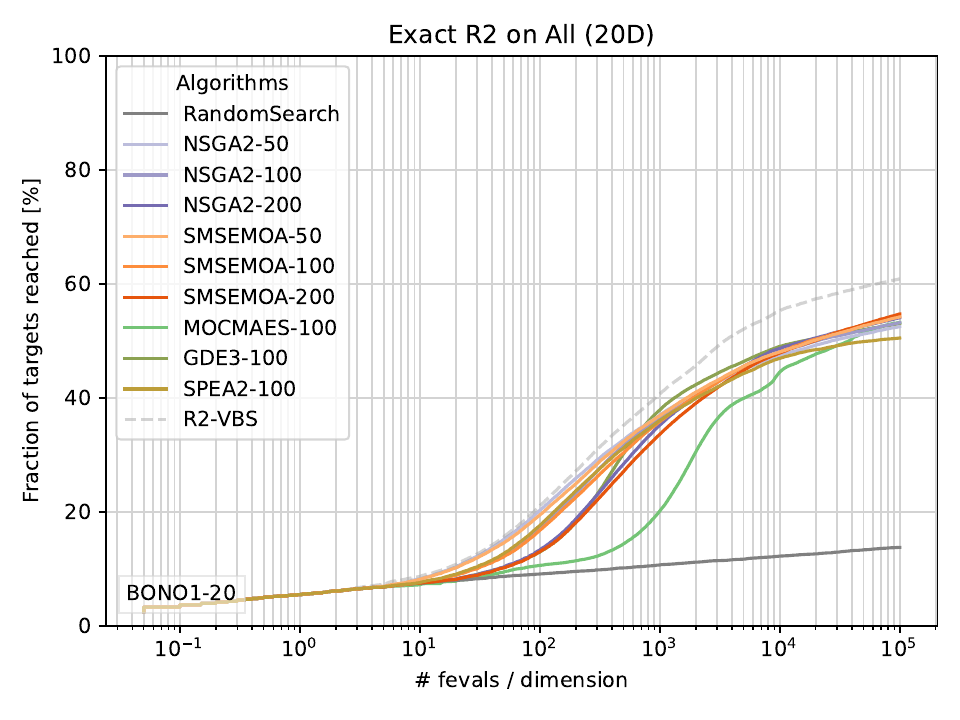}
    \caption{Runtime profiles for all problems with $d \in \{10,20\}$. With increasing decision space dimensionality, the problems become harder to solve. MO-CMA-ES still finally outperforms the competition for $10$-dimensional problems. For $d=20$, while SMS-EMOA still solves the most targets towards the end, GDE3 is the overall best-performing algorithm for intermediate runtimes between $10^3d-10^4d$ evaluations.}
    \label{fig:rtp-all-1020D}
\end{figure}

\Cref{fig:rtp-all-235D,fig:rtp-all-1020D} show the runtime profiles for all problems and algorithms, aggregated by the number of parameters and each indicator separately.
Firstly, we can see that both indicators used result in similar runtime profiles, though some small differences in preference can be distinguished. For example, MO-CMA-ES is evaluated better (relative to the other algorithms) by the exact R2 indicator than the hypervolume.
Furthermore, we observe slightly more solved hypervolume targets, which indicates that additionally, the hypervolume targets tend to focus a little bit more on an earlier part of the convergence to the global front.

Starting with the baseline algorithm, random search, we see it solving comparatively well for the lower-dimensional problems, though it is the worst solver at almost every point in time.
With increasing decision space dimensionality, random search becomes less effective, both absolutely w.r.t.~solved targets as well as relative to the other algorithms; the gap between the optimization algorithms and the baseline widens with an increasing number of decision variables.

Regarding the algorithms, we see that small population sizes lead to faster initial convergence, while NSGA2 and SMS-EMOA variants with larger population sizes perform better beginning with intermediate runtimes (ca. $10^3d$ evaluations), and in the long run (starting with ca. $5 \cdot10^3 d - 10^4 d$ evaluations) MO-CMA-ES dominates in overall performance except for $d=20$.
Despite the small differences between NSGA2, SMS-EMOA, and their population sizes, algorithms from the \texttt{pymoo} package all tend to follow quite similar trajectories.
This can likely be attributed to all other components shared between the different algorithms (cf.~\Cref{sec:setup}), highlighting that just switching out the name-defining selection methodology does not necessarily lead to diverse algorithm behavior.
From our algorithm portfolio, only the MO-CMA-ES provides an essentially different behavior: While its convergence is slower at the beginning, it often ends up solving most targets at the end of the recorded runs.

For the remaining analyses, we focus on the 10-dimensional problems evaluated on the exact R2 indicator.
Runtime profiles for all individual problems, problem groups, and dimensionalities are available at Zenodo \cite{zenodo2025bonobench}.
The repository also includes runtime profiles for individual instances. These reveal that the overall trends shown in our analyses persist throughout most problem instances within a BONO function, though differences can be observed in some cases. Additionally, critical distance plots \cite{demvsar2006statistical} for the performance at the ends of runs are provided.

\subsubsection{Function Groups}

\begin{figure}
    \centering
    \includegraphics[width=0.495\linewidth]{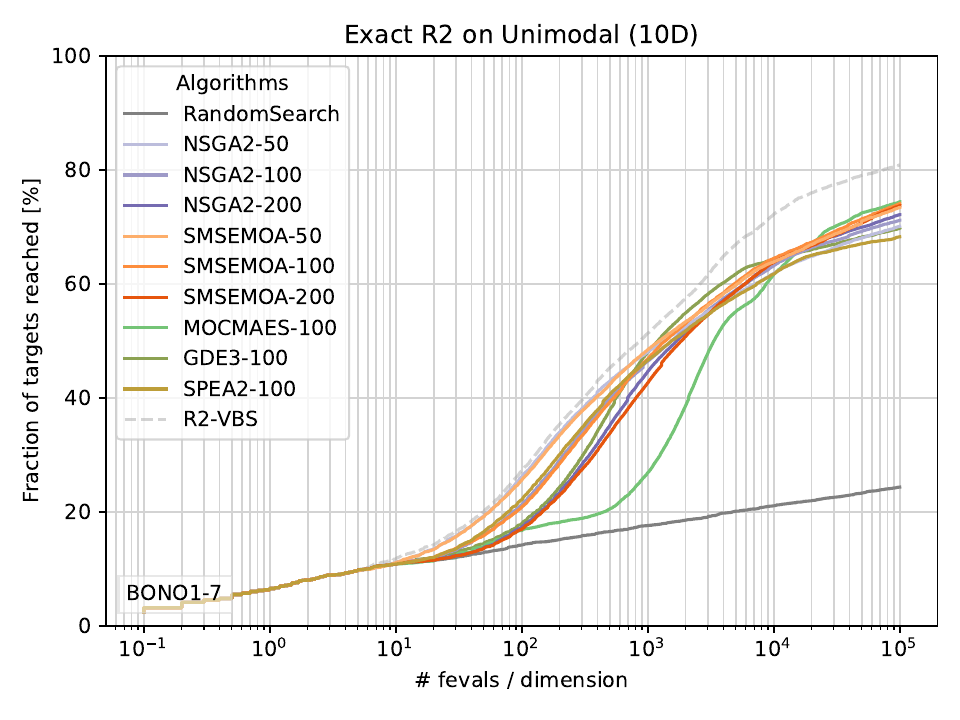}
    \includegraphics[width=0.495\linewidth]{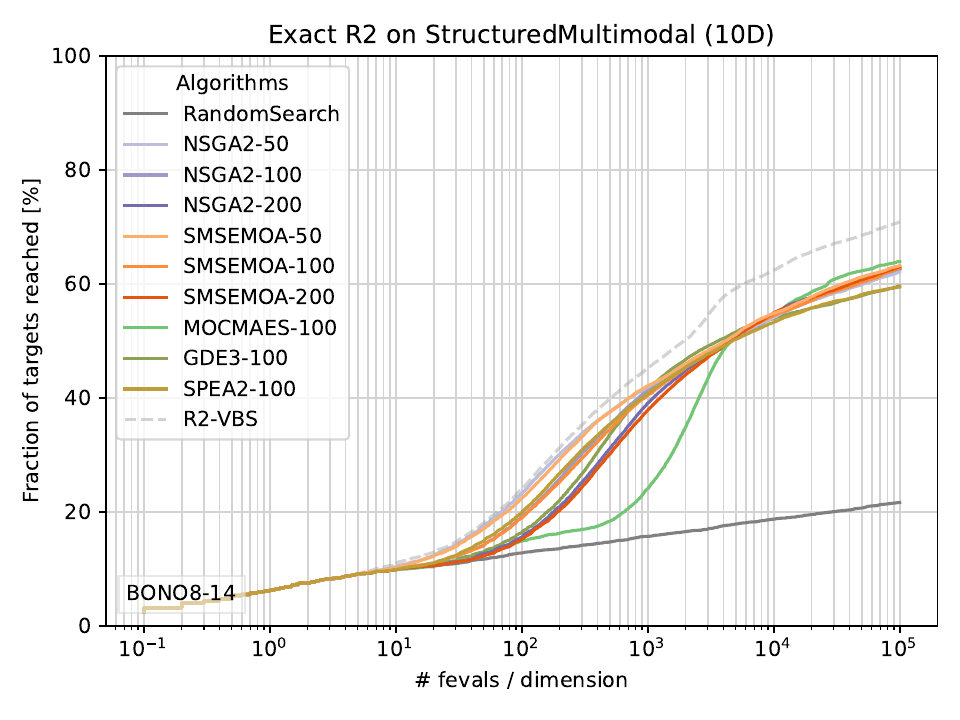}
    \includegraphics[width=0.495\linewidth]{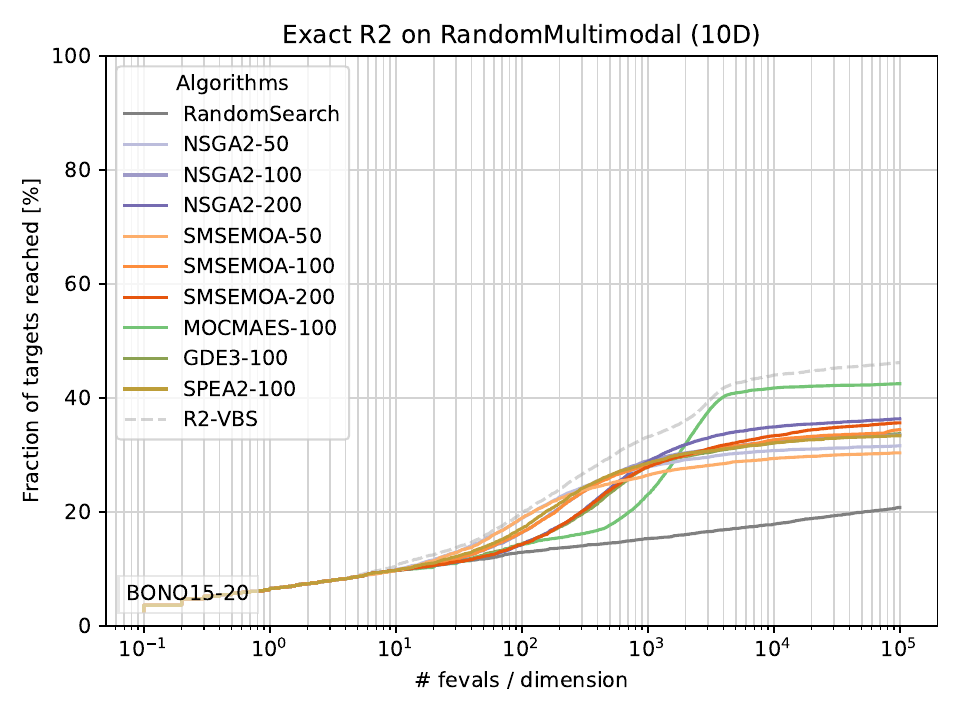}
    \includegraphics[width=0.495\linewidth]{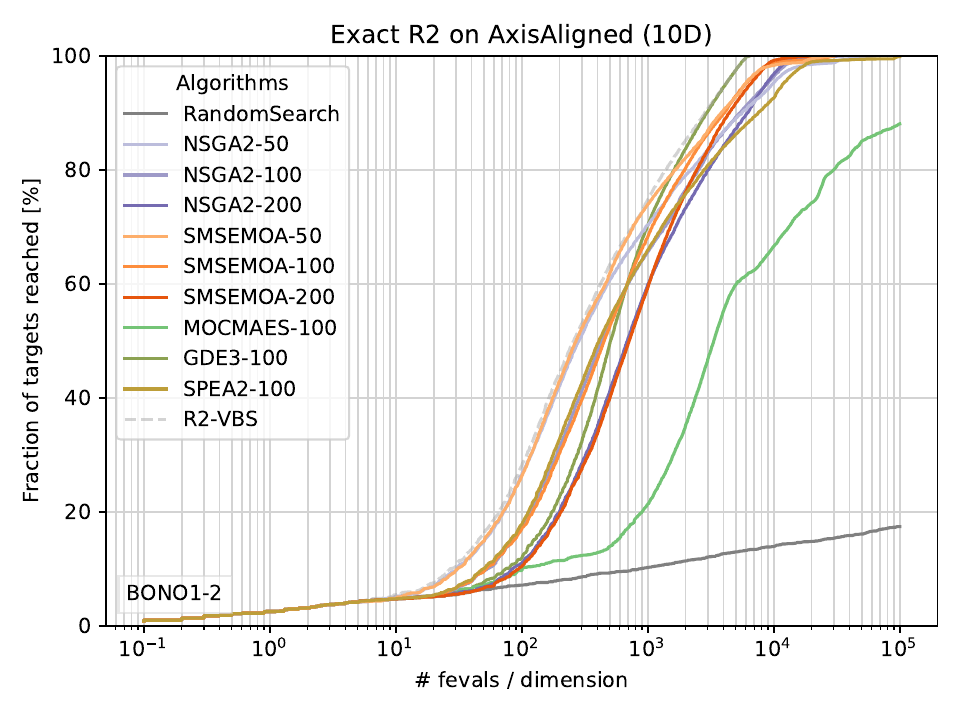}
    \caption{Runtime profiles separated by function group, as well as the combination of just those problems with perfectly axis-aligned Pareto set (BONO1,2).
    The performance on the multimodal problems without global structure is noticeably different from the other two groups. Furthermore, the \texttt{pymoo}-based optimizers can effectively exploit the axis-aligned Pareto sets (lower right plot), which MO-CMA-ES struggles to fully optimize.}
    \label{fig:fgroups}
\end{figure}

In \Cref{fig:fgroups}, we distinguish algorithm performance per function group.
Between the unimodal and the corresponding structured multimodal problems, a performance deterioration can be observed, which also tightens differences between algorithm performances at the end of the runs.
Furthermore, we notice that multimodal problems without global structure ("RandomMultimodal", i.e., BONO15-20) are harder to solve on average, with MO-CMA-ES being clearly the best solver after ca.~$2 \cdot 10^3 d$ evaluations despite plateauing performance around ca.~$4 \cdot 10^3 d$.
The performance plateaus indicate that the algorithms often fail to find all relevant subsets of the Pareto set in one run, and may benefit from being restarted.
In addition to the function groups, \Cref{fig:fgroups} also depicts the performance of the optimizers on the problems with axis-aligned Pareto set: Here, it becomes apparent that the \texttt{pymoo}-based solvers can much more effectively utilize this bias than MO-CMA-ES, which can likely be traced back to the polynomial mutation and simulated binary crossover which exploit separability, while MO-CMA-ES operates in a rotation-invariant manner.

\subsubsection{Front Shape}

\begin{figure}
    \centering
    \includegraphics[width=0.495\linewidth]{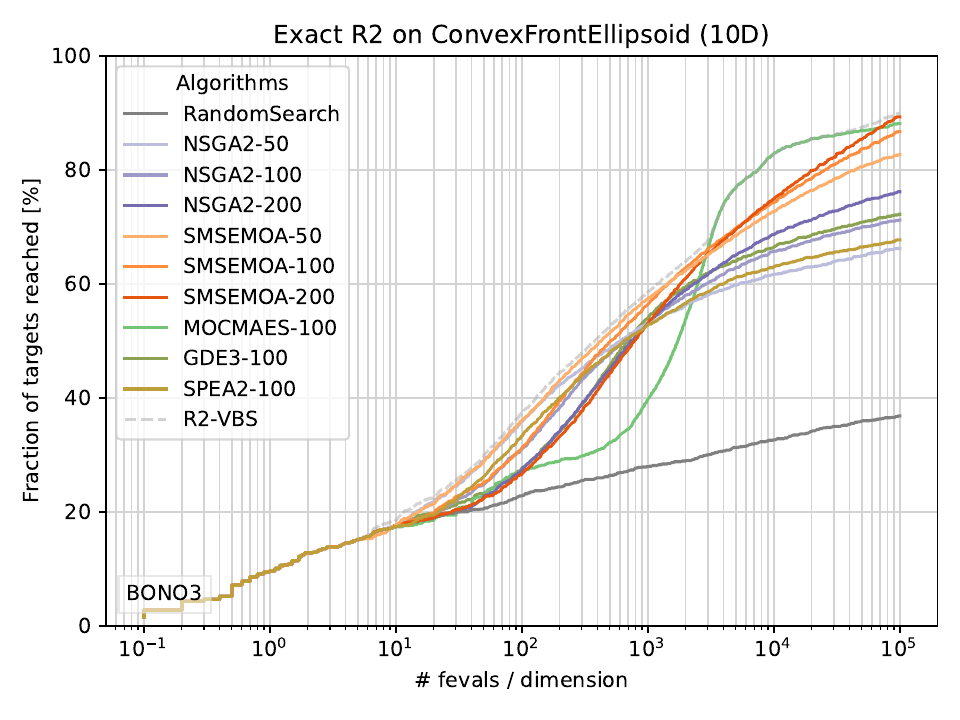}
    \includegraphics[width=0.495\linewidth]{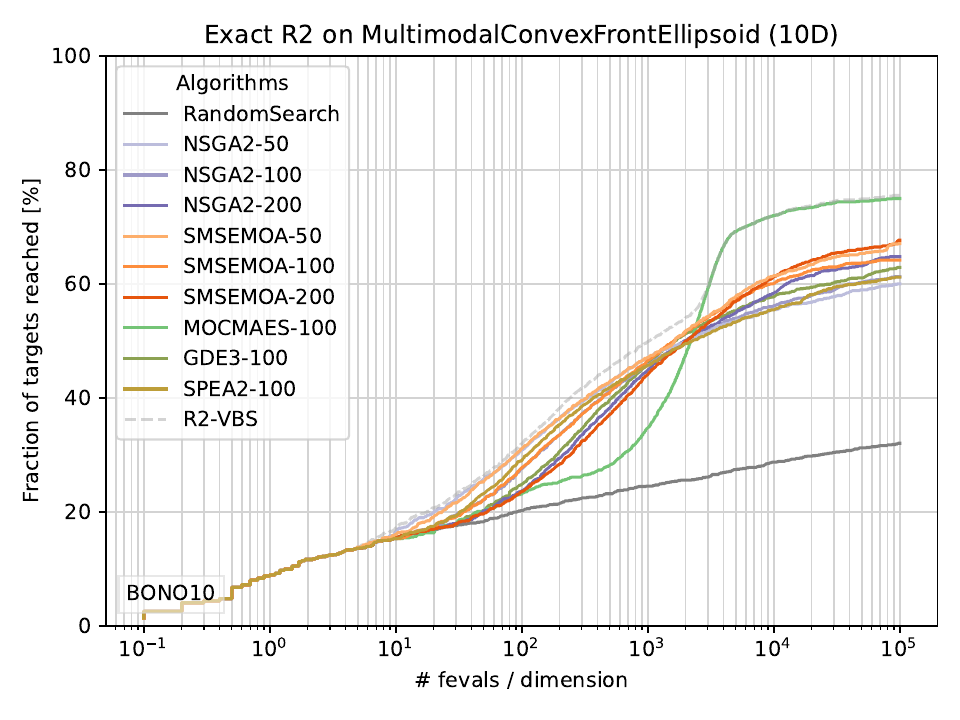}
    \includegraphics[width=0.495\linewidth]{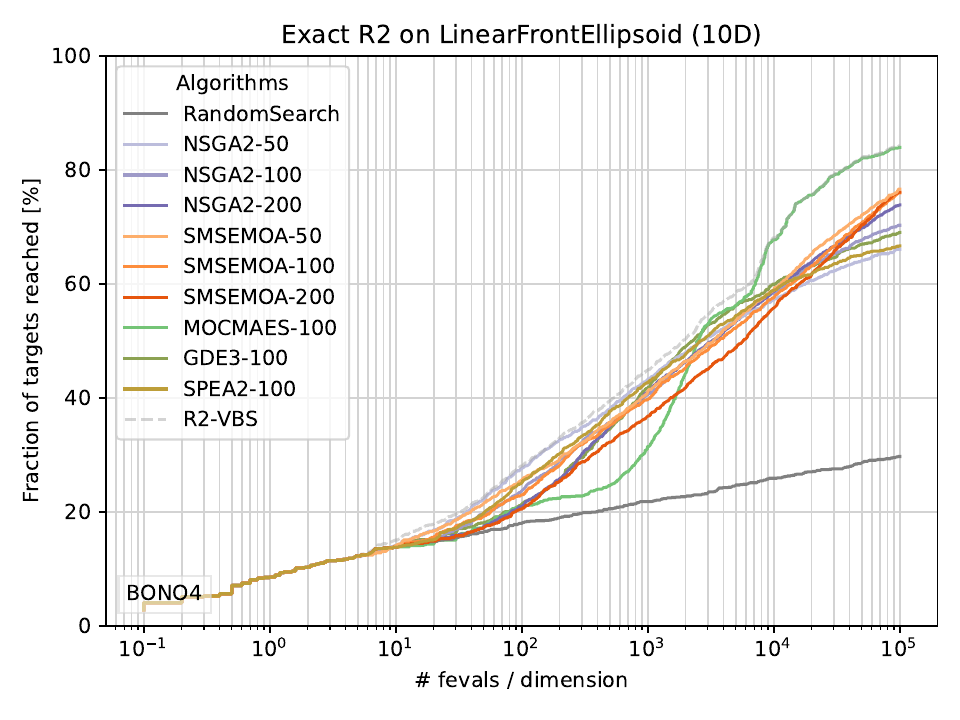}
    \includegraphics[width=0.495\linewidth]{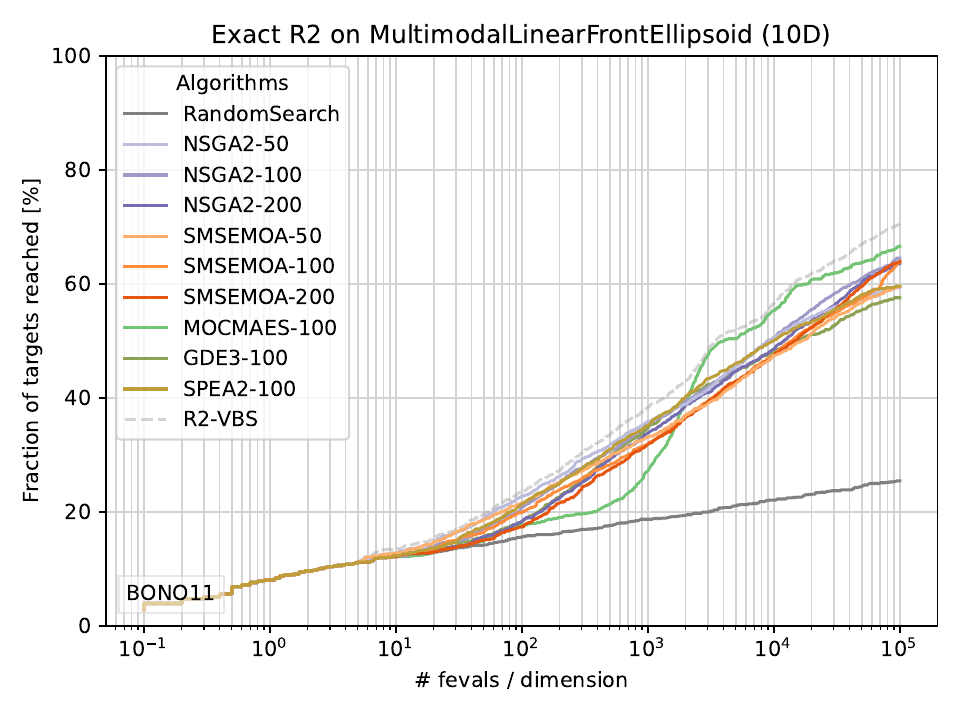}
    \includegraphics[width=0.495\linewidth]{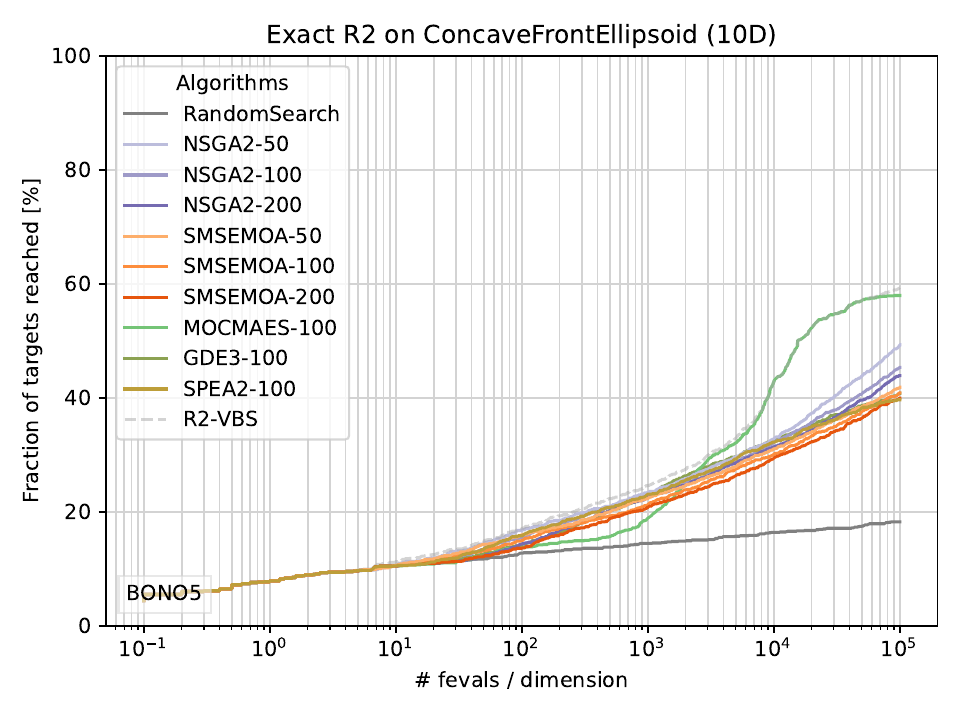}
    \includegraphics[width=0.495\linewidth]{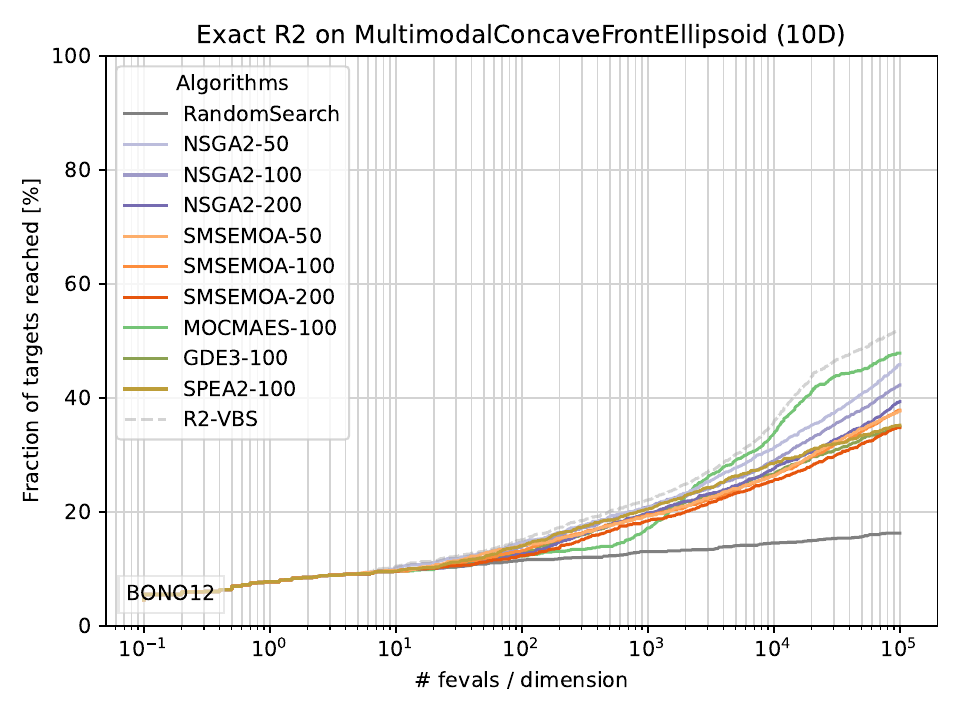}
    \caption{Runtime profiles for convex, linear and concave front ellipsoids in 10D. From the \texttt{pymoo} algorithms, SMS-EMOA has an advantage on the unimodal problems with convex fronts, while performance is more average for linear and deteriorated for concave fronts. Additionally, multimodality has a disproportionate impact on its performance. MO-CMA-ES has the most robust performance of all benchmarked optimizers.}
    \label{fig:rtp-frontshape-10D}
\end{figure}

\Cref{fig:rtp-frontshape-10D} presents the runtime profiles for the problems with linear Pareto set, but different front shapes in both the unimodal as well as their perturbed multimodal variants.
The qualitative pattern is similar to the overall problem set: While MO-CMA-ES is slower in the beginning, it solves the most problems at the end of the runs, however, the distance between MO-CMA-ES and the next best algorithms is larger in most cases.
From the classical MOEAs, SMS-EMOA also reveals an interesting pattern.
Its performance is best on the convex front problem, only slightly better than on the remaining algorithms on the linear fronts, and among the worst on the concave problems.
Beyond this, introducing the perturbation functions also reduces performance for all algorithms, particularly reducing SMS-EMOA's performance as well.
Finally, a general effect of the front shape on the indicator regret becomes evident in this comparison: While random search evaluates the exact same solution trajectories in all problems, a decreasing number of targets is solved from convex to linear and then concave problems.
A reduction in the number of solved targets can also be observed for the other algorithms.
This emphasizes the need for a neutral baseline in order to establish the effectiveness of more sophisticated approaches.

\subsubsection{Discretization}

\begin{figure}
    \centering
    \includegraphics[width=0.495\linewidth]{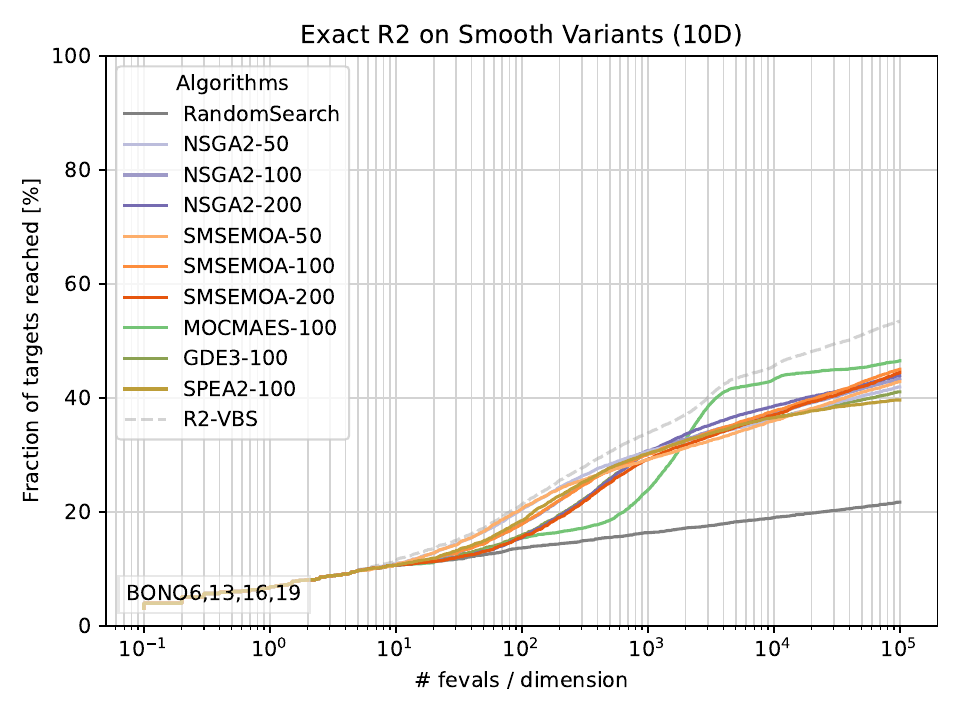}
    \includegraphics[width=0.495\linewidth]{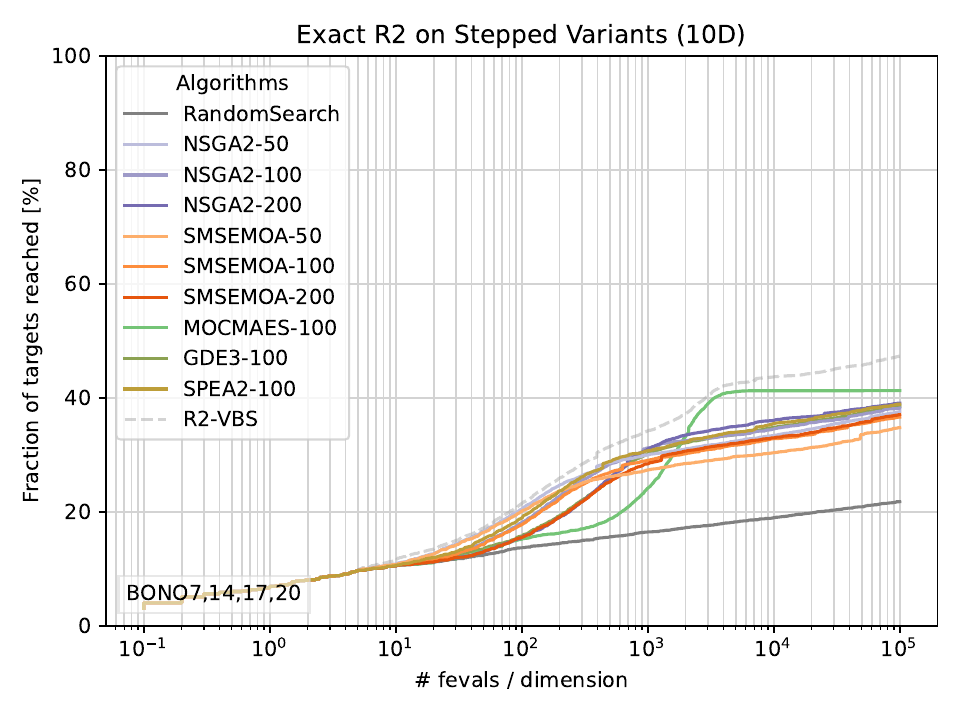}
    \caption{Impact of objective space discretizations on algorithm performance in 10D. The underlying optimization problems are identical except for the discretization. While all optimizers show degenerated performance, SMS-EMOA is particularly sensitive to this change.}
    \label{fig:rtp-stepping-10D}
\end{figure}

Finally, \Cref{fig:rtp-stepping-10D} presents the combined runtime profiles for all stepped problems (BONO7,14,17,20) compared to their counterparts without discretization (BONO6,13,16,19).
The overall performance level is reduced for all algorithms, but SMS-EMOA's performance seems to suffer disproportionately, going from performing among the best (classical) MOEAs to the worst.
Also, MO-CMA-ES hits a performance plateau earlier on the discretized rather than the smooth problems, though it still manages to solve more targets than any other algorithm before plateauing.

\section{Conclusions} \label{sec:conclusions}

Benchmarking is the foundation for progress in many empirical sciences, including multi-objective optimization.
It supports researchers in designing and refining multi-objective optimizers and assists practitioners in deciding on effective optimizers for application-specific problems.
Moreover, benchmarking creates the preliminaries for a large-scale data set, facilitating automated algorithm configuration and selection studies.

This paper contributes to the state of the art in benchmarking bi-objective numerical black-box optimizers by introducing a problem generation approach with configurable single- and bi-objective challenges, capable of defining uni- and multimodal problems with or without discretization in the objective space.
Utilizing the generator, we transform convex-quadratic functions into a test suite of 20 configurable bi-objective problems: the BONO-Bench.
For each problem, we can approximate the Pareto front for the exact R2 and hypervolume indicators up to any given precision, which enables efficient benchmarking of optimizers with precise target values.
This is a trade-off that -- to the authors' best knowledge -- is not yet achieved by any other test suite for bi-objective optimization.
We produce problems with complex optimization challenges, controlled by a generator and with approximation guarantees regarding their Pareto-optimal front to create useful reference values for benchmarking optimizers.

BONO-Bench builds upon and beneficially combines multiple previous works, which have pioneered the usage of convex-quadratic problems for benchmarking \cite{glasmachers2019challenges,toure2019bi}, their combination into multimodal multi-objective problems with tractable Pareto fronts and their indicator approximation \cite{wessing2015multiple,kerschke2019search,schaepermeier2023peak}, as well as problem generators with identifiable properties that are evaluated using target-based runtime profiles \cite{hansen2021coco,brockhoff2022using}.
BONO-Bench comes with the accompanying Python package \texttt{bonobench} that facilitates the benchmarking process, including a flexible problem generator that allows a user to easily define problems beyond the $20$ exemplary classes presented in our test suite.
New problem instances can be generated almost instantaneously, with a one time indicator approximation that finishes in seconds for unimodal problems, with runtimes up to a few minutes for most of the more complex problems.

Despite these advances, there are multiple possible directions to extend this work: First, regarding the problem generator, further placement strategies for the convex-quadratic subproblems or strict and non-strict monotone transformations could be introduced.
One pragmatic limitation of the test problems presented in this paper is the focus on multi-objective problems where the individual objectives have homogeneous challenges.
The creation of heterogeneous problems, e.g., the combination of a discretized unimodal objective with a smooth multimodal one, is entirely possible with the approach as presented in this paper, and should be evaluated in future studies, along with evaluations focusing on the detailed impact of certain problem properties, e.g., degree of multimodality or Pareto front shape.
In particular, future studies could investigate the impact of scaling the degree of multimodality with the dimensionality of the search space.

Another limitation is the usage of a purely numerical approximation to the indicator targets.
A similar target may require a completely different number of points to be solved, e.g., to reach a hypervolume precision of $0.01$ on a linear front, one would need ca. $100$ Pareto-optimal points, while a strongly convex front may be optimized to this degree by just one point.
We believe that revisiting other target-setting approaches, e.g., by performing a subsample of the approximated Pareto front to generate problems that can be reliably solved with a similar number of points, may improve the reliability of runtime profiles.
Extensions and generator configurations should also be guided by the motivation to create problems with properties resembling real-world problems.

Finally, our approach may be extended to problems with more than two objectives.
While convex-quadratic functions can still be analyzed analogously with a higher number of objectives, the approximation to a specific target precision is expected to require more evaluated points in comparison to the bi-objective problems presented here.
Additionally, the degree of multimodality is going to be a limiting factor: $1\,000$ peaks per objective in a bi-objective setting produce the same number of peak combinations as $100$ peaks each with just three objectives.
To create a similarly sophisticated test suite for higher objective space dimensions, additional refinements of the approach may be required.

We plan to extend the experimental studies presented in this paper to eventually have a good representation of the state of the art in bi-objective numerical black-box optimization.
This may also support the development of exploratory landscape analysis features specifically for numerical bi-objective optimization, leading the way to algorithm selection and configuration studies \cite{kerschke2019automated} for this domain.


\begin{acks}
This work was supported by the \href{https://scads.ai}{\em Center for Scalable Data Analytics and Artificial Intelligence (ScaDS.AI) Dresden/Leipzig}.
\end{acks}

\expandafter\def\expandafter\UrlBreaks\expandafter{\UrlBreaks%
  \do\a\do\b\do\c\do\d\do\e\do\f\do\g\do\h\do\i\do\j%
  \do\k\do\l\do\m\do\n\do\o\do\p\do\q\do\r\do\s\do\t%
  \do\u\do\v\do\w\do\x\do\y\do\z\do\A\do\B\do\C\do\D%
  \do\E\do\F\do\G\do\H\do\I\do\J\do\K\do\L\do\M\do\N%
  \do\O\do\P\do\Q\do\R\do\S\do\T\do\U\do\V\do\W\do\X%
  \do\Y\do\Z}

\bibliographystyle{ACM-Reference-Format}
\bibliography{library}










\end{document}